\documentclass[a4paper,12pt,reqno]{amsart}
\pdfoutput=1
\usepackage{amsfonts, amsmath, amssymb, amstext} 
\usepackage{tikz}
\usepackage[normalem]{ulem}
\usepackage{mathrsfs}
\usepackage{hyperref}
\usepackage{enumitem}
\usepackage{tikz-cd, float}
\usetikzlibrary{decorations.markings}

\allowdisplaybreaks[2]

\numberwithin{equation}{section}

\theoremstyle{plain}
\newtheorem{thm}{Theorem}[section] 
\newtheorem{lem}[thm]{Lemma}
\newtheorem{prop}[thm]{Proposition}
\newtheorem{cor}[thm]{Corollary}

\theoremstyle{definition}
\newtheorem{defn}[thm]{Definition}
\newtheorem{exam}[thm]{Example}

\theoremstyle{remark}
\newtheorem{rem}[thm]{Remark}

\newcommand{\norm}[1]{\left\|#1\right\|}
\newcommand{\cspan}{\overline{\mathrm{span}}}
\def\C{{\mathbb C}}
\def\N{{\mathbb N}}

\def\Z{{\mathbb Z}}
\def\T{{\mathbb T}}

\usepackage{braket}

\usepackage{caption, subcaption}
\begin{document}

\title{Higher-rank graph algebras are iterated Cuntz--Pimsner algebras}

\author{James Fletcher}
\address{School of Mathematics and Statistics \\ Victoria University of Wellington \\ Wellington, New Zealand}
\email{james.fletcher@vuw.ac.nz}

\date{\today}
\thanks{This research was supported by an Australian Government Research Training Program (RTP) Scholarship and by the Marsden grant 15-UOO-071 from the 
Royal Society of New Zealand.}



\begin{abstract}
Given a finitely aligned $k$-graph $\Lambda$, we let $\Lambda^i$ denote the $(k-1)$-graph formed by removing all edges of degree $e_i$ from $\Lambda$. We show that the Toeplitz--Cuntz--Krieger algebra of $\Lambda$, denoted by $\mathcal{T}C^*(\Lambda)$, may be realised as the Toeplitz algebra of a Hilbert  $\mathcal{T}C^*(\Lambda^i)$-bimodule. When $\Lambda$ is locally-convex, we show that the Cuntz--Krieger algebra of $\Lambda$, which we denote by $C^*(\Lambda)$, may be realised as the Cuntz--Pimsner algebra of a Hilbert  $C^*(\Lambda^i)$-bimodule. Consequently, $\mathcal{T}C^*(\Lambda)$ and $C^*(\Lambda)$ may be viewed as iterated Toeplitz and iterated Cuntz--Pimsner algebras over $C_0(\Lambda^0)$ respectively.
\end{abstract}

\maketitle


\section{Introduction}

Higher-rank graphs were first introduced by Kumjian and Pask as a generalisation of directed graphs \cite{MR1745529}. Loosely speaking, a higher-rank graph of rank $k$ (or simply a $k$-graph) is a countable small category $\Lambda$ together with a functor $d:\Lambda\rightarrow \N^k$ satisfying the following factorisation property: for any $\lambda\in \Lambda$ and $m,n\in \N^k$ with $d(\lambda)=m+n$, there exist unique $\mu,\nu\in \Lambda$ with $d(\mu)=m$ and $d(\nu)=n$ such that $\lambda=\mu\nu$. In the same paper, Kumjian and Pask showed how to associate a $C^*$-algebra to each row finite higher-rank graph $\Lambda$ with no sources, which we call the Cuntz--Krieger algebra of $\Lambda$. Subsequently Raeburn, Sims, and Yeend \cite{MR2069786}, showed how to relax the hypotheses of \cite{MR1745529}, and defined Cuntz--Krieger algebras for arbitrary finitely aligned higher-rank graphs. Sims subsequently defined relative Cuntz--Krieger algebras for finitely aligned higher-rank graphs, which includes the class of Toeplitz--Cuntz--Krieger algebras as a special case \cite{sims, MR2225454}.

In this article, we show how the Toeplitz--Cuntz--Krieger algebra and Cuntz--Krieger algebra of a finitely aligned higher-rank graph $\Lambda$ may be viewed as iterated Toeplitz and iterated Cuntz--Pimsner algebras over $C_0(\Lambda^0)$ (the space of functions on the graph's vertex set that vanish at infinity) respectively. Writing $e_1,\ldots, e_k$ for the standard generators of $\N^k$, we let $\Lambda^i$ denote the higher-rank graph formed by removing all edges of degree $e_i$ from $\Lambda$. In Theorem~\ref{Toeplitz algebra iterated construction} we show that the Toeplitz--Cuntz--Krieger algebra of $\Lambda$ may be realised as the Toeplitz algebra of a Hilbert bimodule whose coefficient algebra is the Toeplitz--Cuntz--Krieger algebra of $\Lambda^i$. In Theorem~\ref{main theorem} we show that, provided $\Lambda$ is locally-convex, the Cuntz--Krieger algebra of $\Lambda$ may be realised as the Cuntz--Pimsner algebra of a Hilbert bimodule whose coefficient algebra is the Cuntz--Krieger algebra of $\Lambda^i$. Repeatedly removing all edges of a fixed degree from $\Lambda$ eventually leaves a graph consisting solely of vertices, whose Toeplitz--Cuntz--Krieger and Cuntz--Krieger algebras are both isomorphic to $C_0(\Lambda^0)$. When $k=1$ the bimodule we construct is equivalent to the graph correspondence associated to a directed graph \cite[Example~8.3]{MR2135030}, and so we like to think of our procedure as a higher-rank graph correspondence. We also point out that our procedure is similar to the work of Kumjian, Pask, and Sims on $k$-morphs \cite{MR2763728} (introduced as a systematic way of extending a $k$-graph to a $(k+1)$-graph by inserting a collection of edges of degree $e_{k+1}$ between the vertices of the original graph). In \cite[Remark~6.9]{MR2763728} Kumjian, Pask, and Sims show that $C^*(\Lambda)$ may be realised as the Cuntz--Pimsner algebra of a Hilbert $C^*(\Lambda^i)$-bimodule, provided $\Lambda$ is row finite, and has no sources and no sinks (in contrast to our procedure, which only requires local-convexity and finite alignment). 

Our main motivation for wanting to view Cuntz--Krieger algebras associated to higher-rank graphs as iterated Cuntz--Pimsner algebras is to try and determine their $K$-theory. It is well-known that the $K$-theory of a directed graph algebra (equivalently a $1$-graph algebra) can be readily extracted from the graph's adjacency matrix \cite[Theorem~6.1]{MR1988256}. Using a homological spectral sequence, Evans derived expressions for the $K$-theory of Cuntz--Krieger algebras associated to row finite $2$-graphs with no sources, again in terms of the graph's adjacency matrices \cite{evans, MR2383584}. Unfortunately, Evans' techniques do not generalise to $k\geq 3$, and it remains an open problem to find nice formulae for the $K$-groups of higher-rank graph algebras in terms of just their graphical data. In the future, we hope to be able to combine Theorem~\ref{main theorem} and the Pimsner--Voiculescu exact sequence \cite[Theorem~8.6]{MR2102572} (a result that relates the $K$-theory of a Cuntz--Pimsner algebra associated to a Hilbert bimodule and the $K$-theory of the bimodule's coefficient algebra) to do this. As an immediate consequence of combining Theorem~\ref{Toeplitz algebra iterated construction} with \cite[Theorem~4.4]{MR1426840}, we are able to conclude that the Toeplitz--Cuntz--Krieger algebra of a finitely aligned higher-rank graph $\Lambda$ is $KK$-equivalent to $C_0(\Lambda^0)$, generalising an earlier result of Burgstaller \cite[Theorem~1.1]{MR2498556}. Consequently, $K_0(\mathcal{T}C^*(\Lambda))\cong \oplus_{v\in \Lambda^0}\Z$ and $K_1(\mathcal{T}C^*(\Lambda))\cong 0$.

The inspiration for our attempts to realise the Toeplitz--Cuntz--Krieger and Cuntz--Krieger algebras of a finitely aligned higher-rank graph as iterated Toeplitz and Cuntz--Pimsner algebras was Deaconu's work on iterating the Pimsner construction \cite{MR2336239}. Unfortunately, some of Deaconu's proofs lack detail, and it is not clear which of his various hypotheses are necessary to make the procedure work. Motivated by this lack of clarity, as well as the results in \cite[Chapter~2]{fletcherphd}, we explained in \cite{2017arXiv170608626F} how Deaconu's iterative procedure can be extended to quasi-lattice ordered groups that are more general than $(\Z^2,\N^2)$. In particular, \cite[Theorem~4.17]{2017arXiv170608626F} shows that the Nica--Toeplitz algebra of a compactly aligned product system over $\N^k$ can be realised as a $k$-fold iterated Toeplitz algebra. Furthermore, \cite[Theorem~5.20]{2017arXiv170608626F} shows that the Cuntz--Nica--Pimsner algebra can be realised as a $k$-fold iterated Cuntz--Pimsner algebra, provided the action on each fibre of the product system is faithful and by compacts.

In \cite[Section~5.3]{MR2718947} Sims and Yeend show that the Cuntz--Krieger algebra of a finitely aligned $k$-graph may be realised as the Cuntz--Nica--Pimsner algebra of a compactly aligned product system over $\N^k$. It is routine to show that the the action on each fibre of this product system is faithful if and only if the graph has no sources, and by compacts if and only if the graph is row finite. In Subsection~\ref{Iterating the Nica--Toeplitz and Cuntz--Nica--Pimsner construction} we discuss how, in the situation where the graph is row finite and has no sources, Theorem~\ref{main theorem} can be deduced from \cite[Proposition~5.4]{MR2718947} and \cite[Theorem~5.20]{2017arXiv170608626F}. The main purpose of this paper is thus to show that our iterative procedure still works if we drop the hypothesis of row finiteness and the hypothesis of no sources is relaxed to local-convexity. The construction presented in this paper is also significantly simpler and easier to understand than the construction found in \cite{2017arXiv170608626F}, which could be of use to those interested specifically in the $C^*$-algebras associated to higher-rank graphs, and do not want to delve into the theory of product systems. Furthermore, the isomorphisms given by Theorem~\ref{Toeplitz algebra iterated construction} and Theorem~\ref{main theorem} are more explicit than those given by combining the results of \cite{2017arXiv170608626F} and \cite[Section~5.3]{MR2718947}.

In the analysis of \cite[Section~2.6]{fletcherphd}, the assumption that $\Lambda$ has no sources serves two key purposes. Firstly, it ensures that the inclusion of $\Lambda^i$ in $\Lambda$ induces an (injective) $*$-homomorphism from $C^*(\Lambda^i)$ to $C^*(\Lambda)$ (which we use to construct our bimodule), and secondly, it implies that $C^*(\Lambda^i)$ acts faithfully on our bimodule. The assumption that $\Lambda$ is row finite is used to ensure that $C^*(\Lambda^i)$ acts compactly on our bimodule. Combining these two hypotheses, we concluded in \cite[Theorem~2.6.12]{fletcherphd} that the Katsura ideal of our bimodule was all of $C^*(\Lambda^i)$, which made it relatively easy to determine the structure of the bimodule's Cuntz--Pimsner algebra. As shown in \cite[Example~5.4]{2017arXiv170608626F} and Remark~\ref{what goes wrong without local convexity}, if $\Lambda$ has sources, then $C^*(\Lambda)$ need not contain a copy of $C^*(\Lambda^i)$. In Proposition~\ref{injective map for Cuntz--Krieger algebras finitely aligned, locally convex}, we show that this issue can be avoided, provided we restrict our attention to locally-convex graphs. Allowing $\Lambda$ to have sources and/or infinite-receivers can also result in the Katsura ideal being a proper ideal of $C^*(\Lambda^i)$, and the majority of Section~\ref{section: realising as Cuntz--Pimsner algebra} is spent determining what the ideal looks like in this situation. 

Our strategy is to show that the Katsura ideal is gauge-invariant (see  Proposition~\ref{checking gauge invariant ideals}), and then use the results of \cite{MR3262073} to determine its generators. Given a finitely aligned $k$-graph $\Sigma$, it follows from \cite[Theorem~4.6]{MR3262073} that if $I$ is a gauge-invariant ideal of $C^*(\Sigma)$, then $I$ is generated as an ideal by its vertex projections and a collection of projections corresponding to certain finite exhaustive subsets of a subgraph of $\Sigma$. In Proposition~\ref{sufficient to consider exhaustive sets of edges} we show that it suffices to consider only those finite exhaustive sets consisting of edges. We present this result separately in an appendix since it may may be of general interest to those investigating gauge-invariant ideals of higher-rank graph algebras. 

Applying these results to our bimodule, we show in Proposition~\ref{putting everything together} that the Katsura ideal is generated as an ideal by the vertex projections corresponding to vertices admitting a finite and non-zero number of edges of degree $e_i$ (see Proposition~\ref{vertex projection in Katsura ideal}) and a collection of projections corresponding to finite exhaustive subsets of a subgraph of $\Lambda^i$ that can be extended to finite exhaustive subsets of $\Lambda$ (see Lemma~\ref{extending to finite exhaustive set implies compact} for the precise description). With this description of the Katsura ideal it is then relatively straightforward to check that the Cuntz--Pimsner algebra of our bimodule coincides with the Cuntz--Krieger algebra of our original graph. 

Finally, we point out that the results in Section~\ref{section: realising as Cuntz--Pimsner algebra} suggest that the hypothesis of faithful and compact actions present in the author's work on iterating the Cuntz--Nica--Pimsner construction for compactly aligned product systems (see \cite[Theorem~5.20]{2017arXiv170608626F}) can be relaxed (at least for product systems over $\N^k$). The idea would be to develop a suitable notion of local-convexity for product systems (see the discussion before and after \cite[Example~5.4]{2017arXiv170608626F}), and then make use of Katsura's work on gauge-invariant ideals of Cuntz--Pimsner algebras \cite[Theorem~8.6]{MR2413377}. 

\section{Preliminaries}

\subsection{Hilbert bimodules and their associated $C^*$-algebras}

Let $A$ be a $C^*$-algebra. An inner product $A$-module is a complex vector space $X$ equipped with a right action of $A$ and a map $\langle \cdot, \cdot \rangle_A:X\times X\rightarrow A$, linear in its second argument, such that for any $x,y\in X$ and $a\in A$, we have
\begin{enumerate}[label=\upshape(\roman*)]
\item $\langle x,y\rangle_A=\langle y,x\rangle_A^*$;
\item $\langle x,y\cdot a \rangle_A=\langle x,y\rangle_Aa$;
\item $\langle x,x\rangle_A\geq 0$ in $A$; and
\item $\langle x,x\rangle_A=0$ if and only if $x=0$.
\end{enumerate}
It follows from \cite[Proposition~1.1]{lance} that the formula $\norm{x}_X:=\norm{\langle x,x\rangle_A}_A^{1/2}$ defines a norm on $X$. If $X$ is complete with respect to this norm, we say that $X$ is a Hilbert $A$-module. 

We say that a map $T:X\rightarrow X$ is adjointable if there exists a map $T^*:X\rightarrow X$ such that $\langle Tx, y\rangle_A=\langle x, T^*y\rangle_A$ for each $x, y\in X$. Every adjointable  operator $T$ is automatically linear and continuous, and the adjoint $T^*$ is unique. The collection of adjointable operators on $X$, denoted by $\mathcal{L}_A(X)$, equipped with the operator norm is a $C^*$-algebra. For each $x,y\in X$ there is an adjointable operator $\Theta_{x,y}\in \mathcal{L}_A(X)$ defined by $\Theta_{x,y}(z)=x\cdot \langle y, z\rangle_A$. We call operators of this form generalised rank-one operators. The closed subspace $\mathcal{K}_A(X):=\cspan\{\Theta_{x,y}:x,y\in X\}$ is an essential ideal of $\mathcal{L}_A(X)$, whose elements we refer to as generalised compact operators. 

A Hilbert $A$-bimodule consists of a Hilbert $A$-module $X$ together with a $*$-homomorphism $\phi:A\rightarrow \mathcal{L}_A(X)$. We think of $\phi$ as implementing a left action of $A$ on $X$, and frequently write $a\cdot x$ for $\phi(a)(x)$. Since each $\phi(a)\in \mathcal{L}_A(X)$ is $A$-linear, we have that $a\cdot (x\cdot b)=(a\cdot x)\cdot b$ for each $a,b\in A$ and $x\in X$. If we let $A$ act on itself by left and right multiplication, and define an $A$-valued inner product on $A$ by $\langle a,b\rangle_A:=a^*b$, we get a Hilbert $A$-bimodule, which we denote by ${}_A A_A$. We say that a map between two Hilbert $A$-bimodules is a Hilbert $A$-bimodule isomorphism if it is left $A$-linear, surjective, and preserves the $A$-valued inner product (this last condition implies that the map is right $A$-linear and injective). 


The balanced tensor product of two Hilbert $A$-bimodules $X$ and $Y$, denoted by $X\otimes_A Y$, is the completion of the complex vector space spanned by elements $x\otimes_A y$ where $x\in X$ and $y\in Y$, subject to the relation $(x\cdot a) \otimes_A y=x\otimes_A (a\cdot y)$, in the norm determined by the $A$-valued inner product $\langle x\otimes_A y, w\otimes z\rangle_A=\langle y,\langle x,w\rangle_A\cdot z\rangle_A$. There are right and left actions of $A$ on $X\otimes_A Y$ determined by $a\cdot(x\otimes_A y)\cdot b = (a\cdot x)\otimes_A (y\cdot b)$, which gives $X\otimes_A Y$ the structure of a Hilbert $A$-bimodule. We define the balanced tensor powers of $X$ as follows: $X^{\otimes 0}:={}_A A_A$, $X^{\otimes 1}:=X$, and $X^{\otimes n}:=X\otimes_A X^{\otimes n-1}$ for $n\geq 2$. 

A Toeplitz representation of a Hilbert $A$-bimodule $X$ in a $C^*$-algebra $B$ consists of a pair of maps $(\psi,\pi)$, where $\psi:X\rightarrow B$ is linear and $\pi:A\rightarrow B$ is a $*$-homomorphism, satisfying the following relations
\begin{enumerate}
\item[(T1)] $\psi(a\cdot x)=\pi(a)\psi(x)$ for each $a\in A$, $x\in X$;
\item[(T2)] $\psi(x\cdot a)=\psi(x)\pi(a)$ for each $a\in A$, $x\in X$;
\item[(T3)] $\psi(x)^*\psi(y)=\pi(\langle x,y\rangle_A)$ for each $x,y\in X$. 
\end{enumerate}

Given a Hilbert $A$-bimodule $X$, we define the Fock space $\mathcal{F}_X$ to be the set of sequences $(x_n)_{n=0}^\infty$ such that $x_n\in X^{\otimes n}$ for each $n\geq 0$ and $\sum_{n\geq 0} \langle x_n, x_n\rangle_A$ converges in $A$. One can then show that $\langle (x_n)_{n=0}^\infty, (y_n)_{n=0}^\infty \rangle_A: =\sum_{n\geq 0}\langle x_n, y_n\rangle_A$ converges in $A$ for each $(x_n)_{n=0}^\infty, (y_n)_{n=0}^\infty \in \mathcal{F}_X$. Letting $A$ act on $\mathcal{F}_X$ component-wise gives $\mathcal{F}_X$ the structure of a Hilbert $A$-bimodule \cite[p.6]{lance}. It follows that there exists a $*$-homomorphism $\pi:A\rightarrow \mathcal{L}_A(\mathcal{F}_X)$ such that $\pi(a)((x_n)_{n=0}^\infty)=(a\cdot x_n)_{n=0}^\infty$, as well as a linear map $\psi:X\rightarrow \mathcal{L}_A(\mathcal{F}_X)$ such that
\[
\big(\psi(x)((x_n)_{n=0}^\infty)\big)_m
=\begin{cases}
0 & \text{if $m=0$}\\
x\cdot x_0 & \text{if $m=1$}\\
x\otimes_A x_{m-1} & \text{if $m\geq 2$}.
\end{cases}
\]
Routine calculations show that the pair $(\psi,\pi)$ is a Toeplitz representation of $X$ in $\mathcal{L}_A(\mathcal{F}_X)$, which we call the Fock representation of $X$. 

Proposition~1.8 of \cite{MR1722197} shows that a Toeplitz representation $(\psi,\pi)$ of a Hilbert $A$-bimodule $X$ gives rise to Toeplitz representations of the tensor powers of $X$. If we define $\psi^{\otimes 0}:=\pi$, $\psi^{\otimes 1}:=\psi$, and, for $n\geq 2$, let $\psi^{\otimes n}$ be the linear map determined inductively by $\psi^{\otimes n}(x\otimes_A y)=\psi(x)\psi^{\otimes n-1}(y)$ for each $x\in X$ and $y\in X^{\otimes n-1}$, then $(\psi^{\otimes n},\pi)$ is a Toeplitz representation of $X^{\otimes n}$ for each $n\in \N\cup \{0\}$. Using relations (T1)--(T3) and the Hewitt--Cohen--Blanchard factorisation theorem \cite[Proposition~2.31]{MR1634408}, it can be shown that the $C^*$-subalgebra generated by $\psi(X)\cup \pi(A)$ is $\cspan\left\{\psi^{\otimes m}(x)\psi^{\otimes n}(y)^*:m,n\geq 0, x\in X^{\otimes m}, y\in X^{\otimes n} \right\}$.

Theorem~2.10 of \cite{MR2679392} can be used to show that there exists a $C^*$-algebra $\mathcal{T}_X$, which we call the Toeplitz algebra of $X$, and a Toeplitz representation $(i_X,i_A)$ of $X$ in $\mathcal{T}_X$, that are universal in the following sense:
\begin{enumerate}[label=\upshape(\roman*)]
\item $\mathcal{T}_X$ is generated by $i_X(X)\cup i_A(A)$;
\item given any Toeplitz representation $(\psi,\pi)$ of $X$ in a $C^*$-algebra $B$, there exists a $*$-homomorphism $\psi\times_\mathcal{T} \pi:\mathcal{T}_X\rightarrow B$ such that 
\[
(\psi\times_\mathcal{T} \pi) \circ i_X=\psi
\quad \text{and} \quad
(\psi\times_\mathcal{T} \pi) \circ i_A=\pi.
\]
\end{enumerate}
It follows that $\mathcal{T}_X$ is equal to $\cspan\{i_X^{\otimes m}(x)i_X^{\otimes n}(y)^*:m,n\geq 0, x\in X^{\otimes m}, y\in X^{\otimes n}\}$.

The universal property of the Toeplitz algebra ensures it carries a strongly continuous action of the circle group $\gamma:\T\rightarrow \mathrm{Aut}(\mathcal{T}_X)$, which we call the gauge action. The action is determined by $\gamma_z(i_X(x))=zi_X(x)$ and $\gamma_z(i_A(a))=i_A(a)$ for each $z\in \T$, $x\in X$, and $a\in A$.

In \cite{MR2102572} Katsura defined what has come to be accepted as the correct notion of a Cuntz--Pimsner algebra for a Hilbert bimodule with a non-faithful left action. Given a Toeplitz representation $(\psi,\pi)$ of a Hilbert $A$-bimodule $X$ in a $C^*$-algebra $B$, by \cite[Proposition~8.11]{MR2135030} there exists a $*$-homomorphism $(\psi,\pi)^{(1)}:\mathcal{K}_A(X)\rightarrow B$ such that $(\psi,\pi)^{(1)}\left(\Theta_{x,y}\right)=\psi(x)\psi(y)^*$ for each $x,y\in X$. We also let $\ker(\phi)^\perp:=\{a\in A:ab=0 \text{ for all $b\in \ker(\phi)$}\}$. We then say that $(\psi,\pi)$ is Cuntz--Pimsner covariant if $(\psi,\pi)^{(1)}(\phi(a))=\pi(a)$ for every $a\in J_X:=\phi^{-1}(\mathcal{K}_A(X))\cap \ker(\phi)^\perp$.

Theorem~2.10 of \cite{MR2679392} can again be used to show that there exists a $C^*$-algebra $\mathcal{O}_X$, which we call the Cuntz--Pimsner algebra of $X$, and a Cuntz--Pimsner covariant Toeplitz representation $(j_X,j_A)$ of $X$ in $\mathcal{O}_X$ that are universal in the following sense:
\begin{enumerate}[label=\upshape(\roman*)]
\item $\mathcal{O}_X$ is generated by $j_X(X)\cup j_A(A)$;
\item given any Cuntz--Pimsner covariant Toeplitz representation $(\psi,\pi)$ of $X$ in a $C^*$-algebra $B$, there exists a $*$-homomorphism $\psi\times_\mathcal{O} \pi:\mathcal{O}_X\rightarrow B$ such that 
\[
(\psi\times_\mathcal{O} \pi) \circ j_X=\psi
\quad \text{and} \quad
(\psi\times_\mathcal{O} \pi) \circ j_A=\pi.
\]
\end{enumerate}
It follows that $\mathcal{O}_X$ is a quotient of $\mathcal{T}_X$, and routine calculations show that the gauge action on the Toeplitz algebra descends to this quotient. 

\subsection{Higher-rank graphs and their associated $C^*$-algebras}

A higher-rank graph of rank $k$ (also known simply as a $k$-graph) consists of a countable small category $\Lambda$ and a functor $d:\Lambda\rightarrow \mathbb{N}^k$, called the degree map, satisfying the following factorisation property: for any $m,n\in \mathbb{N}^k$ and $\lambda\in \Lambda$ with $d(\lambda)=m+n$, there exist unique $\mu,\nu\in \Lambda$ with $d(\mu)=m$ and $d(\nu)=n$ such that $\lambda=\mu\circ\nu$. Since we think of the morphisms in the category as paths in a graph, we write $\lambda\mu$ for $\lambda \circ \mu$ whenever $\lambda,\mu \in \Lambda$ with $\mathrm{dom}(\lambda)=\mathrm{cod}(\mu)$. 

The factorisation property has some important consequences. Firstly, it follows easily that $d^{-1}(0)=\{\mathrm{id}_o:o\in \mathrm{Obj}(\Lambda)\}$. Secondly, if $\lambda\in \Lambda$ and $m,n\in\mathbb{N}^k$ are such that $m\leq n\leq d(\lambda)$, then two applications of the factorisation property shows that there exist unique $\mu,\nu,\eta\in \Lambda$ with $\lambda=\mu\nu\eta$ and $d(\mu)=m$, $d(\nu)=n-m$, $d(\eta)=d(\lambda)-n$. We write $\lambda(0,m)$ for $\mu$, $\lambda(m,n)$ for $\nu$, and $\lambda(n,d(\lambda))$ for $\eta$. 

The following notation and terminology is standard when working with higher-rank graphs. We write $e_i$ for the $i$th generator of $\mathbb{N}^k$, and $n_i$ for the $i$th component of $n\in \mathbb{N}^k$. We define a partial order on $\mathbb{N}^k$ by $m\leq n \iff m_i\leq n_i$ for all $i$. For any nonempty finite set $E:=\{m_1,\ldots, m_n\}\subseteq \mathbb{N}^k$, we write $\bigvee E$ and $\bigwedge E$ for the component-wise maximum and minimum of $m_1,\ldots, m_n$ (and define $\bigvee \emptyset$ and $\bigwedge \emptyset$ to be zero). For simplicity's sake, we write $m\vee n$ for $\bigvee \{m,n\}$, and $m\wedge n$ for $\bigwedge \{m,n\}$. For each $n\in \mathbb{N}^k$, we define $\Lambda^n:=\{\lambda\in \Lambda:d(\lambda)=n\}$. For each $\lambda\in \Lambda$, we define $r(\lambda):=\mathrm{id}(\mathrm{cod}(\lambda))\in \Lambda^0$ and $s(\lambda):=\mathrm{id}(\mathrm{dom}(\lambda))\in \Lambda^0$. The maps $r,s:\Lambda \rightarrow \Lambda^0$ are called the range and source maps of $\Lambda$. Given a subset $E\subseteq \Lambda$ and $\lambda\in \Lambda$, we define $\lambda E:=\{\lambda\mu:\mu\in E, s(\lambda)=r(\mu)\}$ and $E\lambda :=\{\mu\lambda:\mu\in E, \ r(\lambda)=s(\mu)\}$. We say that a $k$-graph $\Lambda$ has no sources if for every $v\in \Lambda^0$ and every $n\in \mathbb{N}^k$, the set $v\Lambda^n$ is nonempty. For $n\in \mathbb{N}^k$, we define $\Lambda^{\leq n}:=\{\lambda\in \Lambda: d(\lambda)\leq n \text{ and } d(\lambda)_i<n_i \implies s(\lambda)\Lambda^{e_i}=\emptyset\}$ (a simple induction argument shows that each $v\Lambda^{\leq n}$ is always nonempty). We say that a $k$-graph $\Lambda$ is locally-convex if whenever $\lambda\in \Lambda^{e_i}$ and $\mu\in \Lambda^{e_j}$ with $i\neq j$ and $r(\lambda)=r(\mu)$, we have $s(\lambda)\Lambda^{e_j}\neq \emptyset$ and $s(\mu)\Lambda^{e_i}\neq \emptyset$.

Before we look at associating $C^*$-algebras to higher-rank graphs, we need to discuss the concept of (minimal) common extensions. For $\mu,\nu\in \Lambda$ we set
\begin{align*}
\mathrm{CE}(\mu,\nu)&:=\mu\Lambda\cap \nu\Lambda \quad \text{and}\\
\mathrm{MCE}(\mu,\nu)&:=\mathrm{CE}(\mu,\nu)\cap \Lambda^{d(\mu)\vee d(\nu)}.
\end{align*}
We call elements of $\mathrm{CE}(\mu,\nu)$ common extensions of $\mu$ and $\nu$, and elements of $\mathrm{MCE}(\mu,\nu)$ minimal common extensions of $\mu$ and $\nu$. We also define
\[
\Lambda^{\min}(\mu,\nu):=\{(\alpha,\beta): \mu\alpha=\nu\beta \in \mathrm{MCE}(\mu,\nu)\}.
\]
That is, a common extension of $\mu,\nu\in \Lambda$ is a path that ends with both $\mu$ and $\nu$, and a minimal common extension is a common extension that has minimal degree (i.e. $d(\mu)\vee d(\nu)$). Elements of $\Lambda^{\min}(\mu,\nu)$ are then ordered pairs of paths that when prepended to $\mu$ and $\nu$ respectively give a minimal common extension. The factorisation property implies that if $\lambda$ is a common extension of $\mu$ and $\nu$, then $\lambda(0,d(\mu)\vee d(\nu))$ is a minimal common extension and 
$\left(\lambda\left(d(\mu),d(\mu)\vee d(\nu)\right), \lambda\left(d(\nu),d(\mu)\vee d(\nu)\right)\right)\in \Lambda^{\min}(\mu,\nu)$. We can also extend the notion of minimal common extensions from pairs of paths to arbitrary nonempty finite subsets $G\subseteq \Lambda$ by setting $\mathrm{CE}(G):=\bigcap_{\nu\in G}\nu\Lambda$ and $\mathrm{MCE}(G):=\mathrm{CE}(G)\cap \Lambda^{\bigvee d(G)}$. We say that a higher-rank graph $\Lambda$ is finitely aligned if $\Lambda^{\min}(\mu,\nu)$ is finite (possibly empty) for every $\mu,\nu \in \Lambda$ (equivalently $\mathrm{MCE}(\mu,\nu)$ is finite for every $\mu,\nu \in \Lambda$).

Given $v\in \Lambda^0$, we say that a set $E\subseteq v\Lambda$ is exhaustive in $\Lambda$ if for each $\mu\in v\Lambda$ there exists $\nu\in E$ such that $\Lambda^{\min}(\mu,\nu)$ is nonempty. We point out that if $v\in E$, then $E$ is automatically exhaustive. We write 
\[
\mathrm{FE}(\Lambda):=\bigcup_{v\in \Lambda^0} \{E\subseteq v\Lambda\setminus \{v\}: \text{$E$ is finite and exhaustive in $\Lambda$}\}.
\]
For $E\in \mathrm{FE}(\Lambda)$, we write $r(E)$ for the vertex $v\in \Lambda^0$ such that $E\subseteq v\Lambda$. We also define $v\mathrm{FE}(\Lambda):=\{E\in \mathrm{FE}(\Lambda):r(E)=v\}$. 

We now define Toeplitz--Cuntz--Krieger families for finitely aligned $k$-graphs. We say that a collection $\{q_\lambda:\lambda\in \Lambda\}$ of elements in a $C^*$-algebra is a Toeplitz--Cuntz--Krieger $\Lambda$-family if
\begin{enumerate}[leftmargin=1.9cm]
\item[$\mathrm{(TCK1)}$] $\{q_v:v\in \Lambda^0\}$ is a set of mutually orthogonal projections;
\item[$\mathrm{(TCK2)}$] $q_\mu q_\nu=q_{\mu\nu}$ for all $\mu,\nu \in \Lambda$ with $s(\mu)=r(\nu)$; 
\item[$\mathrm{(TCK3)}$] $q_\mu^* q_\nu =\sum_{(\alpha,\beta)\in \Lambda^{\min}(\mu,\nu)} q_\alpha q_\beta^*$ for all $\mu,\nu \in \Lambda$, where the empty sum is interpreted as zero. 
\end{enumerate}
It follows from relation (TCK3) that $q_\lambda^* q_\mu=\delta_{\lambda,\mu}q_{s(\lambda)}$ for $\lambda,\mu\in \Lambda$ with $d(\lambda)=d(\mu)$, and so by (TCK1) Toeplitz--Cuntz--Krieger families consist of partial isometries. Furthermore, relations (TCK1)--(TCK3) imply that $C^*(\{q_\lambda:\lambda\in \Lambda\})=\cspan\{q_\lambda q_\mu^*:\lambda,\mu\in \Lambda\}$. Given a vertex $v\in \Lambda^0$ and a finite set $E\subseteq v\Lambda$, we fix the following notation
\[
\Delta(q)^E:=\prod_{\lambda\in E}(q_v-q_\lambda q_\lambda^*).
\]

Using Theorem~2.10 of \cite{MR2679392}, it can be shown that there exists a $C^*$-algebra $\mathcal{T}C^*(\Lambda)$, which we call the Toeplitz--Cuntz--Krieger algebra of $\Lambda$, and a Toeplitz--Cuntz--Krieger $\Lambda$-family $\{t_\lambda^\Lambda:\lambda \in \Lambda\}$ in $\mathcal{T}C^*(\Lambda)$, that are universal in the following sense
\begin{enumerate}[label=\upshape(\roman*)]
\item 
$\mathcal{T}C^*(\Lambda)$ is generated by $\{t_\lambda^\Lambda:\lambda \in \Lambda\}$;
\item 
if $\{q_\lambda:\lambda\in \Lambda\}$ is a Toeplitz--Cuntz--Krieger $\Lambda$-family in a $C^*$-algebra $B$, then there exists a $*$-homomorphism $\pi_q:\mathcal{T}C^*(\Lambda)\rightarrow B$ that carries $t_\lambda^\Lambda$ to $q_\lambda$ for each $\lambda\in \Lambda$.
\end{enumerate}

It is useful to know when the $*$-homomorphism induced by the universal property of $\mathcal{T}C^*(\Lambda)$ is faithful: by \cite[Theorem~3.15]{MR3262073}, if $\{q_\lambda:\lambda \in \Lambda\}$ is a  Toeplitz--Cuntz--Krieger $\Lambda$-family, then $\pi_q$ is faithful provided each vertex projection $q_v$ is nonzero and $\Delta(q)^E\neq 0$ for each $E\in \mathrm{FE}(\Lambda)$.

We say that a Toeplitz--Cuntz--Krieger $\Lambda$-family $\{q_\lambda:\lambda\in \Lambda\}$ is a Cuntz--Krieger $\Lambda$-family if 
\begin{enumerate}[leftmargin=1.9cm]
\item[(CK)] $\Delta(q)^E=0$ for each $E\in \mathrm{FE}(\Lambda)$. 
\end{enumerate}
It follows from Theorem~2.10 of \cite{MR2679392} that there exists a $C^*$-algebra $C^*(\Lambda)$, which we call the Cuntz--Krieger algebra of $\Lambda$, and a Cuntz--Krieger $\Lambda$-family $\{s_\lambda^\Lambda:\lambda \in \Lambda\}$ in $C^*(\Lambda)$, that are universal in the following sense
\begin{enumerate}[label=\upshape(\roman*)]
\item 
$C^*(\Lambda)$ is generated by $\{s_\lambda^\Lambda:\lambda \in \Lambda\}$;
\item 
if $\{q_\lambda:\lambda\in \Lambda\}$ is a Cuntz--Krieger $\Lambda$-family in a $C^*$-algebra $B$, then there exists a $*$-homomorphism $\pi_q:C^*(\Lambda)\rightarrow B$ that carries $s_\lambda^\Lambda$ to $q_\lambda$ for each $\lambda\in \Lambda$.
\end{enumerate}

The universal property of the Cuntz--Krieger algebra gives the existence of an action $\gamma^\Lambda:\T^k\rightarrow \mathrm{Aut}(C^*(\Lambda))$, which we call the gauge action, such that $\gamma_z^\Lambda(s_\lambda^\Lambda)=z^{d(\lambda)}s_\lambda^\Lambda$ for each $\lambda\in \Lambda$ and $z\in \T^k$ (where $z^m:=\prod_{i=1}^k z_i^{m_i}$ for each $m\in \mathbb{N}^k$). An $\varepsilon/3$ argument shows that $\gamma^\Lambda$ is strongly continuous.

The gauge action can be used to determine when representations of Cuntz--Krieger algebras are faithful (see \cite[Theorem~4.2]{MR2069786}). If $\pi:C^*(\Lambda)\rightarrow B$ is a representation in a $C^*$-algebra $B$, then $\pi$ is injective provided $\pi(s_v^\Lambda)\neq 0$ for each $v\in \Lambda^0$ and there exists a strongly continuous action $\theta:\mathbb{T}^k\rightarrow \mathrm{Aut}(C^*(\{\pi(s_\lambda^\Lambda):\lambda\in \Lambda\}))$
such that $\theta_z\circ \pi=\pi\circ \gamma^\Lambda_z$ for each $z\in \T^k$.

\section{Realising $\mathcal{T}C^*(\Lambda)$ as a Toeplitz algebra}
\label{section: realising as Toeplitz algebra}

Given a $k$-graph $\Lambda$ (with $k\geq 1$), we fix $i\in \{1,\ldots, k\}$ and let $\Lambda^i:=\{\lambda\in \Lambda:d(\lambda)_i=0\}$ (i.e. we remove all edges of degree $e_i$ from $\Lambda$). Restricting the degree functor gives $\Lambda^i$ the structure of a $(k-1)$-graph. In this section we show how the Toeplitz--Cuntz--Krieger algebra of $\Lambda$ may be realised as the Toeplitz algebra of a Hilbert $\mathcal{T}C^*(\Lambda^i)$-bimodule. We will define the Hilbert $\mathcal{T}C^*(\Lambda^i)$-bimodule that we are interested in to be a certain closed subspace of $\mathcal{T}C^*(\Lambda)$. To equip this set with left and right actions of $\mathcal{T}C^*(\Lambda^i)$ we want a $*$-homomorphism from $\mathcal{T}C^*(\Lambda^i)$ to $\mathcal{T}C^*(\Lambda)$. Moreover, to ensure that we have a $\mathcal{T}C^*(\Lambda^i)$-valued inner product, we need to know that this $*$-homomorphism is injective. 

\begin{prop}
\label{existence of phi for TCK algebra}
Let $\Lambda$ be a finitely aligned $k$-graph. Then there exists an injective $*$-homomorphism $\phi:\mathcal{T}C^*(\Lambda^i)\rightarrow \mathcal{T}C^*(\Lambda)$ such that $\phi(t_\lambda^{\Lambda^i})=t_\lambda^{\Lambda}$ for each $\lambda\in \Lambda^i$. 
\begin{proof}
Clearly, the collection $\{t_\lambda^\Lambda:\lambda\in \Lambda^i\}$ satisfies (TCK1) and (TCK2). To see that $\{t_\lambda^\Lambda:\lambda\in \Lambda^i\}$ also satisfies (TCK3), it suffices to show that $\Lambda^{\min}(\mu,\nu)=(\Lambda^i)^{\min}(\mu,\nu)$ for any $\mu,\nu\in \Lambda^i$. To see this, observe that for any $(\alpha,\beta)\in\Lambda^{\min}(\mu,\nu)$ we have 
\[
d(\alpha)_i=\left(d(\mu)\vee d(\nu)-d(\mu)\right)_i=\max\{d(\mu)_i, d(\nu)_i\}-d(\mu)_i=0
\] 
and 
\[
d(\beta)_i=\left(d(\mu)\vee d(\nu)-d(\nu)\right)_i=\max\{d(\mu)_i, d(\nu)_i\}-d(\nu)_i=0,
\]
and so $(\alpha,\beta)\in(\Lambda^i)^{\min}(\mu,\nu)$. Thus, $\{t_\lambda^\Lambda:\lambda\in \Lambda^i\}$ is a Toeplitz--Cuntz--Krieger $\Lambda^i$-family in $\mathcal{T}C^*(\Lambda)$, and so by the universal property of $\mathcal{T}C^*(\Lambda^i)$, there exists a $*$-homomorphism $\phi:\mathcal{T}C^*(\Lambda^i)\rightarrow \mathcal{T}C^*(\Lambda)$ such that $\phi(t_\lambda^{\Lambda^i})=t_\lambda^{\Lambda}$ for each $\lambda\in \Lambda^i$. It remains to check that $\phi$ is injective. 

Routine calculations show that for each $\lambda\in \Lambda$, there exists $w_\lambda\in \mathcal{B}(\ell^2(\Lambda))$ such that $w_\lambda\xi_\mu=\delta_{s(\lambda),r(\mu)}\xi_{\lambda\mu}$ for each $\mu \in \Lambda$ (where $\{\xi_\lambda:\lambda\in \Lambda\}$ is the canonical orthonormal basis for $\ell^2(\Lambda)$). Further straightforward calculations show that the adjoint of $w_\lambda$ is determined by the formula 
\[
w_\lambda^*\xi_\nu=
\begin{cases}
\xi_\eta &\text{if $\nu=\lambda\eta$ for some $\eta\in \Lambda$}\\
0 & \text{otherwise,}
\end{cases}
\]
and that the collection $\{w_\lambda:\lambda\in \Lambda\}$ is a Toeplitz--Cuntz--Krieger $\Lambda$-family in $\mathcal{B}(\ell^2(\Lambda))$. Since the $*$-homomorphism $\pi_w:\mathcal{T}C^*(\Lambda)\rightarrow \mathcal{B}(\ell^2(\Lambda))$ that sends $t_\lambda^\Lambda$ to $w_\lambda$ is necessarily norm-decreasing, and $w_\lambda \xi_{s(\lambda)}=\xi_\lambda\neq 0$ for each $\lambda\in \Lambda$, we conclude that each $t_\lambda^\Lambda$ in the universal Toeplitz--Cuntz--Krieger $\Lambda$-family is nonzero. In particular, $t_v^\Lambda\neq 0$ for each $v\in \Lambda^0$. Thus, to prove that $\phi$ is injective, by \cite[Theorem~3.15]{MR3262073}, it remains to show that $\Delta(t^\Lambda)^E\neq 0$ for each $E\in \mathrm{FE}(\Lambda^i)$. A simple calculation shows that for each $\mu\in \Lambda$,
\[
\Delta(w)^E\xi_\mu=
\begin{cases}
\xi_\mu & \text{if }r(\mu)=r(E) \ \text{and } \mu \not \in \lambda \Lambda \ \text{for all } \lambda\in E \\
0 & \text{otherwise.}
\end{cases}
\]
Since $r(E)\not\in E$, we have that $\pi_w(\Delta(t^\Lambda)^E)\xi_{r(E)}=\Delta(w)^E\xi_{r(E)}=\xi_{r(E)}\neq 0$, and so $\Delta(t^\Lambda)^E\neq 0$. We conclude that $\phi$ is injective. 
\end{proof}
\end{prop}

Using the injective $*$-homomorphism from the previous proposition, we define a collection of Hilbert $\mathcal{T}C^*(\Lambda^i)$-bimodules. 

\begin{prop}
\label{module for TCK algebra}
Let $\Lambda$ be a finitely aligned $k$-graph. For each $n\geq 0$, define 
\[
X_n:=\cspan\{t_\lambda^\Lambda {t_\mu^\Lambda}^*:\lambda,\mu\in \Lambda, d(\lambda)_i=n, d(\mu)_i=0\}\subseteq \mathcal{T}C^*(\Lambda),
\] 
taking the closure with respect to the norm on $\mathcal{T}C^*(\Lambda)$. Then $X_n$ is a Hilbert $\mathcal{T}C^*(\Lambda^i)$-module with inner product and right action given by 
\begin{equation}
\label{inner product and right action}
\langle x,y \rangle^n_{\mathcal{T}C^*(\Lambda^i)}=\phi^{-1}(x^*y) \quad \text{and} 
\quad 
x\cdot a=x\phi(a)
\quad \text{for $x,y\in X_n$, $a\in \mathcal{T}C^*(\Lambda^i)$.}
\end{equation}
Moreover, the norm on $X_n$ induced by $\langle \cdot, \cdot \rangle^n_{\mathcal{T}C^*(\Lambda^i)}$ agrees with the norm on $\mathcal{T}C^*(\Lambda)$. Additionally, there exists a $*$-homomorphism $\psi_n: \mathcal{T}C^*(\Lambda^i)\rightarrow \mathcal{L}_{\mathcal{T}C^*(\Lambda^i)}(X_n)$ such that $\psi_n(a)(x)=\phi(a)x$ for each $a\in \mathcal{T}C^*(\Lambda^i)$ and $x\in X_n$, giving $X_n$ the structure of a Hilbert $\mathcal{T}C^*(\Lambda^i)$-bimodule.
\begin{proof}
By \cite[Lemma~3.2(1)]{MR1426840} (see also \cite{MR1638139} for a more categorical approach), if
\begin{enumerate}[label=\upshape(\roman*)]
\item $X_n^*X_n\subseteq \phi(\mathcal{T}C^*(\Lambda^i))$; and 
\item $X_n\phi(\mathcal{T}C^*(\Lambda^i))\subseteq X_n$,
\end{enumerate}
then $X_n$ is a Hilbert $\mathcal{T}C^*(\Lambda^i)$-module with inner product and right action given by \eqref{inner product and right action}, and the norm on $X_n$ agrees with the norm on $\mathcal{T}C^*(\Lambda)$.

Let us check that (i) holds. Fix paths $\lambda,\lambda',\mu,\mu'\in \Lambda$ with $d(\lambda)_i=d(\lambda')_i=n$ and $d(\mu)_i=d(\mu')_i=0$. If $(\alpha,\beta)\in \Lambda^{\min}(\lambda,\lambda')$, then 
\begin{align*}
&d(\mu\alpha)_i=d(\mu)_i+d(\alpha)_i=d(\alpha)_i=\max\{d(\lambda)_i,d(\lambda')_i\}-d(\lambda)_i=n-n=0,\\
&d(\mu'\beta)_i=d(\mu')_i+d(\beta)_i=d(\beta)_i=\max\{d(\lambda)_i,d(\lambda')_i\}-d(\lambda')_i=n-n=0.
\end{align*}
Hence, making use of relation (TCK3), we see that
\begin{align*}
\big(t_\lambda^\Lambda {t_\mu^\Lambda}^*\big)^*\big(t_{\lambda'}^\Lambda {t_{\mu'}^\Lambda}^*\big)=t_\mu^\Lambda {t_\lambda^\Lambda}^* t_{\lambda'}^\Lambda{t_{\mu'}^\Lambda}^*=\sum_{(\alpha,\beta)\in \Lambda^{\min}(\lambda,\lambda')}t_{\mu\alpha}^\Lambda{t_{\mu'\beta}^\Lambda}^*
\in \phi(\mathcal{T}C^*(\Lambda^i)).
\end{align*}
Since the adjoint and multiplication are both continuous on $\mathcal{T}C^*(\Lambda)$ and $\phi(\mathcal{T}C^*(\Lambda^i))$ is a $*$-subalgebra of $\mathcal{T}C^*(\Lambda)$, we conclude that $X_n^*X_n\subseteq \phi(\mathcal{T}C^*(\Lambda^i))$. 

Now we check that (ii) holds. Again, by linearity and continuity, it suffices to show that if $\lambda,\eta, \rho, \mu \in \Lambda$ with $d(\lambda)_i=n$ and $d(\eta)_i=d(\rho)_i=d(\mu)_i=0$, then 
$
t_\lambda^\Lambda {t_\mu^\Lambda}^*\phi(t_\eta^{\Lambda^i}  {t_\rho^{\Lambda^i}}^*)\in X_n. 
$
Observe that if $(\alpha,\beta)\in \Lambda^{\min}(\mu,\eta)$, then
\[
d(\lambda\alpha)_i=d(\lambda)_i+d(\alpha)_i=n+d(\alpha)_i=n+\max\{d(\mu)_i,d(\eta)_i\}-d(\mu)_i=n
\]
and
\[
d(\rho\beta)_i=d(\rho)_i+d(\beta)_i=\max\{d(\mu)_i,d(\eta)_i\}-d(\eta)_i=0.
\]
Hence,
\[
t_\lambda^\Lambda {t_\mu^\Lambda}^*\phi\big(t_\eta^{\Lambda^i}  {t_\rho^{\Lambda^i}}^*\big)=t_\lambda^\Lambda {t_\mu^\Lambda}^*t_\eta^\Lambda {t_\rho^\Lambda}^* =\sum_{(\alpha,\beta)\in \Lambda^{\min}(\mu,\eta)}t_{\lambda\alpha}^\Lambda {t_{\rho\beta}^\Lambda}^*\in X_n
\]
as required.

It remains to show that there exists a $*$-homomorphism $\psi_n:\mathcal{T}C^*(\Lambda^i)\rightarrow \mathcal{L}_{\mathcal{T}C^*(\Lambda^i)}(X_n)$ such that $\psi_n(a)(x)=\phi(a)x$ for $a\in \mathcal{T}C^*(\Lambda^i)$ and $x\in X_n$. Observe that if $\lambda, \eta, \rho, \mu \in \Lambda$ with $d(\lambda)_i=n$ and $d(\eta)_i=d(\rho)_i=d(\mu)_i=0$, and $(\alpha,\beta)\in \Lambda^{\min}(\rho,\lambda)$, then
\[
d(\eta\alpha)_i=d(\eta)_i+d(\alpha)_i=d(\alpha)_i=\max\{d(\rho)_i,d(\lambda)_i\}-d(\rho)_i=d(\lambda)_i-0=n
\]
and
\[
d(\mu\beta)_i=d(\mu)_i+d(\beta)_i=d(\beta)_i=\max\{d(\rho)_i,d(\lambda)_i\}-d(\lambda)_i=d(\lambda)_i-d(\lambda)_i=0.
\] 
Thus, an application of relation (TCK3) shows that
\[
\phi\big(t_\eta^{\Lambda^i} {t_\rho^{\Lambda^i}}^*\big)t_\lambda^\Lambda {t_\mu^\Lambda}^*=t_\eta^\Lambda {t_\rho^\Lambda}^* t_\lambda^\Lambda {t_\mu^\Lambda}^*=\sum_{(\alpha,\beta)\in \Lambda^{\min}(\rho,\lambda)}t_{\eta\alpha}^\Lambda {t_{\mu\beta}^\Lambda}^*\in X_n. 
\]
By linearity and continuity it follows that $\phi(\mathcal{T}C^*(\Lambda^i))X_n\subseteq X_n$. It follows from \cite[Lemma~3.2(2)]{MR1426840} that for each $a\in \mathcal{T}C^*(\Lambda^i)$ the map $\psi_n(a):X_n\rightarrow X_n$ defined by $\psi_n(a)(x):=\phi(a)x$ is adjointable. Since $\phi$ is a $*$-homomorphism, it follows that the map $\psi_n:\mathcal{T}C^*(\Lambda^i)\rightarrow  \mathcal{L}_{\mathcal{T}C^*(\Lambda^i)}(X_n)$ is also a $*$-homomorphism.
\end{proof}
\end{prop}

Our aim is to show that the Toeplitz algebra of the Hilbert $\mathcal{T}C^*(\Lambda^i)$-bimodule $X:=X_1$ is isomorphic to the Toeplitz--Cuntz--Krieger algebra of $\Lambda$. Before we do this we need to analyse the tensor powers of $X$. Firstly, we need a lemma telling us how, given paths $\eta,\rho\in \Lambda$, we can factorise elements of $\Lambda^{\min}(\eta,\rho)$. 

\begin{lem}
\label{factorising minimal common extensions}
Let $\Lambda$ be a finitely aligned $k$-graph. For each $\eta,\rho \in \Lambda$ and $m\in \mathbb{N}^k$ with $m\leq d(\rho)$, we have
\[
\Lambda^{\min}(\eta,\rho)=\{(\alpha\gamma,\delta):(\alpha,\beta)\in \Lambda^{\min}(\eta,\rho(0,m)), \ (\gamma,\delta)\in \Lambda^{\min}(\beta,\rho(m,d(\rho)))\}.
\]
\begin{proof}
To start we prove that 
\[
\{(\alpha\gamma,\delta):(\alpha,\beta)\in \Lambda^{\min}(\eta,\rho(0,m)), \ (\gamma,\delta)\in \Lambda^{\min}(\beta,\rho(m,d(\rho)))\}\subseteq \Lambda^{\min}(\eta,\rho).
\]
Fix $(\alpha,\beta)\in \Lambda^{\min}(\eta,\rho(0,m))$ and $(\gamma,\delta)\in \Lambda^{\min}(\beta,\rho(m,d(\rho)))$. Then
\begin{align*}
\eta\alpha\gamma=\rho(0,m)\beta\gamma=\rho(0,m)\rho(m,d(\rho))\delta=\rho\delta,
\end{align*}
which shows that $\eta\alpha\gamma=\rho\delta\in \mathrm{CE}(\eta,\rho)$. We show that the common extension $\eta\alpha\gamma=\rho\delta$ of $\eta$ and $\rho$ is minimal by computing the degree of $\rho\delta$. Since $(\gamma,\delta)\in \Lambda^{\min}(\beta,\rho(m,d(\rho)))$, we see that
\begin{align*}
d(\rho\delta)
=d(\rho(0,m))+d(\rho(m,d(\rho))\delta)
&=m+d(\beta)\vee d(\rho(m,d(\rho))\\
&=m+d(\beta)\vee\left(d(\rho)-m\right).
\end{align*}
Since $(\alpha,\beta)\in \Lambda^{\min}(\eta,\rho(0,m))$, this must be the same as
\begin{align*}
m+\left(d(\eta)\vee d(\rho(0,m))-d(\rho(0,m))\right)\vee \left(d(\rho)-m\right)
=m+\left(d(\eta)\vee m-m\right)\vee \left(d(\rho)-m\right).
\end{align*}
Fix $i\in\{1,\ldots, k\}$. If $d(\eta)_i\geq m_i$, then
\begin{align*}
\left(m+\left(d(\eta)\vee m-m\right)\vee \left(d(\rho)-m\right)\right)_i
&=m_i+\max\{d(\eta)_i-m_i,d(\rho)_i-m_i\}\\
&=\max\{d(\eta)_i,d(\rho)_i\}
=\left(d(\eta)\vee d(\rho)\right)_i.
\end{align*}
On the other hand, suppose that $d(\eta)_i< m_i$. Using that fact that $d(\eta)_i<m_i\leq d(\rho)_i$ for the penultimate equality, we see that
\begin{align*}
\left(m+\left(d(\eta)\vee m-m\right)\vee \left(d(\rho)-m\right)\right)_i
&=m_i+\max\{0,d(\rho)_i-m_i\}\\
&=d(\rho)_i
=\max\{d(\rho)_i,d(\eta)_i\}
=\left(d(\eta)\vee d(\rho)\right)_i.
\end{align*}
Thus, $d(\rho\delta)=d(\eta)\vee d(\rho)$, and we conclude that $(\alpha\gamma,\delta)\in \Lambda^{\min}(\eta,\rho)$.  

Next we check that 
\[
\Lambda^{\min}(\eta,\rho)\subseteq\{(\alpha\gamma,\delta):(\alpha,\beta)\in \Lambda^{\min}(\eta,\rho(0,m)), \ (\gamma,\delta)\in \Lambda^{\min}(\beta,\rho(m,d(\rho)))\}.
\]
Suppose that $(\lambda,\tau)\in \Lambda^{\min}(\eta,\rho)$ and define paths $\alpha:=\lambda(0,d(\eta)\vee m - d(\eta))$, $\beta:=(\rho\tau)(m,d(\eta)\vee m)$, $\gamma:=\lambda(d(\eta)\vee m-d(\eta)$, $d(\lambda))$, and $\delta:=\tau$. By construction, $(\alpha\gamma,\delta)=(\lambda,\tau)$. Thus it remains to show that $(\alpha,\beta)\in \Lambda^{\min}(\eta, \rho(0,m))$ and $(\gamma,\delta)\in \Lambda^{\min}(\beta, \rho(m,d(\rho)))$. Since $\eta\lambda=\rho\tau$, we see that
\begin{align*}
\eta\alpha
&=\eta\lambda(0,d(\eta)\vee m - d(\eta))
=(\eta\lambda)(0,d(\eta)\vee m)\\
&=(\rho\tau)(0,d(\eta)\vee m) 
=(\rho\tau)(0,m)\,(\rho\tau)(m,d(\eta\vee m)).
\end{align*}
As $m\leq d(\rho)$, this must be the same as
\begin{align*}
\rho(0,m)\,(\rho\tau)(m,d(\eta\vee m))
=\rho(0,m)\beta.
\end{align*}
Hence, $\eta\alpha=\rho(0,m)\beta\in \mathrm{CE}(\eta,\rho(0,m))$. Since
\begin{align*}
d(\eta\alpha)=d(\eta)+d(\alpha)=d(\eta)+d(\eta)\vee m -d(\eta)=d(\eta)\vee m=d(\eta)\vee d(\rho(0,m)),
\end{align*}
we conclude that $(\alpha,\beta)\in \Lambda^{\min}(\eta, \rho(0,m))$. Since $d(\eta)\leq d(\eta)\vee m$, we see that 
\begin{align*}
\beta\gamma
&=(\rho\tau)(m,d(\eta)\vee m)\,\lambda(d(\eta)\vee m-d(\eta),d(\lambda))\\
&=(\eta\lambda)(m,d(\eta)\vee m)\,\lambda(d(\eta)\vee m-d(\eta),d(\lambda))\\
&=(\eta\lambda)(m,d(\eta)\vee m)\,(\eta\lambda)(d(\eta)\vee m,d(\eta\lambda))\\
&=(\eta\lambda)(m,d(\eta\lambda)).
\end{align*}
As $\eta\lambda=\rho\tau$ and $m\leq d(\rho)$ this is equal to
\begin{align*}
(\rho\tau)(m,d(\rho\tau))
=\rho(m,d(\rho))\tau
=\rho(m,d(\rho))\delta.
\end{align*}
Thus, $\beta\gamma=\rho(m,d(\rho))\delta\in \mathrm{CE}(\beta,\rho(m,d(\rho)))$. Since  $(\lambda,\tau)\in \Lambda^{\min}(\eta,\rho)$, we have
\begin{align*}
d(\beta\gamma)
&=d(\eta)\vee m -m+d(\lambda)-d(\eta)\vee m +d(\eta)
=d(\lambda)+d(\eta)-m\\
&=d(\eta)\vee d(\rho)-d(\eta)+d(\eta)-m
=\left(d(\eta)\vee d(\rho)\right)-m.
\end{align*}
As $m\leq d(\rho)$, this is the same as
\begin{align*}
d(\eta)\vee m\vee d(\rho)-m
&=(d(\eta)\vee m-m)\vee (d(\rho)-m)\\
&=d(\beta)\vee (d(\rho)-m)
=d(\beta)\vee d(\rho(m,d(\rho))),
\end{align*}
which shows that $(\gamma, \delta)\in \Lambda^{\min}(\beta,\rho(m,d(\rho)))$.
\end{proof}
\end{lem}

\begin{prop}
\label{omega isomorphisms}
Let $\Lambda$ be a finitely aligned $k$-graph. Define $X_n$ as in Proposition~\ref{module for TCK algebra} and set $X:=X_1$. Then for each $n\in \mathbb{N}\cup \{0\}$,  there exists a Hilbert $\mathcal{T}C^*(\Lambda^i)$-bimodule isomorphism $\Omega_n:X_n\rightarrow X^{\otimes n}$ such that $\Omega_0=\phi^{-1}$ and, for $n\geq 1$,
\begin{equation}
\label{inductively defined omega maps}
\begin{aligned}
\Omega_n\big(t_\lambda^\Lambda{t_\mu^\Lambda}^*\big)&=t_{\lambda(0,e_i)}^\Lambda\otimes_{\mathcal{T}C^*(\Lambda^i)}\Omega_{n-1}\big(t_{\lambda(e_i,d(\lambda))}^\Lambda {t_\mu^\Lambda}^*\big) \\
&\text{for each $\lambda,\mu\in \Lambda$ with $d(\lambda)_i=n$ and $d(\mu)_i=0$.}
\end{aligned}
\end{equation}
\begin{proof}
Define $\Omega_0:X_0\rightarrow X^{\otimes 0}=\mathcal{T}C^*(\Lambda^i)$ to be $\phi^{-1}$. Clearly, $\Omega_0$ is a Hilbert $\mathcal{T}C^*(\Lambda^i)$-bimodule isomorphism. For $n\geq 1$, we claim that there exists a Hilbert $\mathcal{T}C^*(\Lambda^i)$-bimodule isomorphism $\Omega_n:X_n\rightarrow X^{\otimes n}$ satisfying \eqref{inductively defined omega maps}. We will define this collection of maps inductively. 

Fix $n\geq 0$ and suppose that $\Omega_n:X_n\rightarrow X^{\otimes n}$ is a Hilbert $\mathcal{T}C^*(\Lambda^i)$-bimodule isomorphism satisfying \eqref{inductively defined omega maps}. Let $\lambda,\mu, \nu,\eta\in \Lambda$ with $d(\lambda)_i=d(\nu)_i=n+1$ and $d(\mu)_i=d(\eta)_i=0$. Using the fact that $\Omega_n$ is left $\mathcal{T}C^*(\Lambda^i)$-linear for the second equality, we see that
\begin{align*}
\big\langle t_{\lambda(0,e_i)}^\Lambda&\otimes_{\mathcal{T}C^*(\Lambda^i)}\Omega_n\big(t_{\lambda(e_i,d(\lambda))}^\Lambda {t_\mu^\Lambda}^*\big), t_{\nu(0,e_i)}^\Lambda\otimes_{\mathcal{T}C^*(\Lambda^i)}\Omega_n\big(t_{\nu(e_i,d(\lambda))}^\Lambda {t_\eta^\Lambda}^*\big)\big\rangle\\
&=
\big\langle \Omega_n\big(t_{\lambda(e_i,d(\lambda))}^\Lambda {t_\mu^\Lambda}^*\big), \langle  t_{\lambda(0,e_i)}^\Lambda, t_{\nu(0,e_i)}^\Lambda\rangle^1_{\mathcal{T}C^*(\Lambda^i)} \cdot \Omega_n\big(t_{\nu(e_i,d(\lambda))}^\Lambda {t_\eta^\Lambda}^*\big)\big\rangle\\
&=
\big\langle \Omega_n\big(t_{\lambda(e_i,d(\lambda))}^\Lambda {t_\mu^\Lambda}^*\big),  \Omega_n\big(\langle  t_{\lambda(0,e_i)}^\Lambda, t_{\nu(0,e_i)}^\Lambda\rangle^1_{\mathcal{T}C^*(\Lambda^i)} \cdot t_{\nu(e_i,d(\lambda))}^\Lambda {t_\eta^\Lambda}^*\big)\big\rangle. 
\end{align*}
Since $\Omega_n$ is inner product preserving, this is equal to 
\begin{align*}
\big\langle t_{\lambda(e_i,d(\lambda))}^\Lambda {t_\mu^\Lambda}^*,  {t_{\lambda(0,e_i)}^\Lambda}^*t_{\nu(0,e_i)}^\Lambda t_{\nu(e_i,d(\lambda))}^\Lambda {t_\eta^\Lambda}^*\big\rangle^n_{\mathcal{T}C^*(\Lambda^i)}
&=
\phi^{-1}\big(t_\mu^\Lambda{t_{\lambda(e_i,d(\lambda))}^\Lambda}^* {t_{\lambda(0,e_i)}^\Lambda}^*t_{\nu(0,e_i)}^\Lambda t_{\nu(e_i,d(\lambda))}^\Lambda {t_\eta^\Lambda}^*\big)\\
&=
\phi^{-1}\big(t_\mu^\Lambda{t_\lambda^\Lambda}^* t_\nu^\Lambda {t_\eta^\Lambda}^*\big)\\
&=\big\langle t_\lambda^\Lambda{t_\mu^\Lambda}^*, t_\nu^\Lambda {t_\eta^\Lambda}^* \big\rangle^{n+1}_{\mathcal{T}C^*(\Lambda^i)}.
\end{align*}
Thus, there exists a well-defined norm-decreasing map 
\[
\sum c_{(\lambda,\mu)} t_\lambda^\Lambda{t_\mu^\Lambda}^*
\mapsto 
\sum c_{(\lambda,\mu)} t_{\lambda(0,e_i)}^\Lambda\otimes_{\mathcal{T}C^*(\Lambda^i)}\Omega_n\big(t_{\lambda(e_i,d(\lambda))}^\Lambda {t_\mu^\Lambda}^*\big)
\]
on $\mathrm{span}\big\{t_\lambda^\Lambda{t_\mu^\Lambda}^*:\lambda,\mu\in \Lambda, d(\lambda)_i=n+1, d(\mu)_i=0\big\}$, which extends to $X_{n+1}$ by continuity. We denote this extension by $\Omega_{n+1}$. The previous calculation then shows that $\Omega_{n+1}$ is inner product preserving. 

We now show that $\Omega_{n+1}$ is left $\mathcal{T}C^*(\Lambda^i)$-linear. For any $\lambda,\mu, \nu,\eta\in \Lambda$ with $d(\lambda)_i=n+1$ and $d(\nu)_i=d(\mu)_i=d(\eta)_i=0$, we have
\begin{equation}
\label{calculations for left linearity of omega maps}
\begin{aligned}
t_\nu^{\Lambda^i} {t_\eta^{\Lambda^i}}^*\cdot \Omega_{n+1}\big(t_\lambda^\Lambda{t_\mu^\Lambda}^*\big)
&=
t_\nu^{\Lambda^i} {t_\eta^{\Lambda^i}}^*\cdot \big( t_{\lambda(0,e_i)}^\Lambda\otimes_{\mathcal{T}C^*(\Lambda^i)}\Omega_n\big(t_{\lambda(e_i,d(\lambda))}^\Lambda {t_\mu^\Lambda}^*\big)\big)\\
&=
t_\nu^\Lambda {t_\eta^\Lambda}^* t_{\lambda(0,e_i)}^\Lambda\otimes_{\mathcal{T}C^*(\Lambda^i)}\Omega_n\big(t_{\lambda(e_i,d(\lambda))}^\Lambda {t_\mu^\Lambda}^*\big)\\
&=
\sum_{(\alpha,\beta)\in \Lambda^{\min}(\eta,\lambda(0,e_i))}t_{\nu\alpha}^\Lambda {t_\beta^\Lambda}^* \otimes_{\mathcal{T}C^*(\Lambda^i)}\Omega_n\big(t_{\lambda(e_i,d(\lambda))}^\Lambda {t_\mu^\Lambda}^*\big).
\end{aligned}
\end{equation}
To simplify this expression, observe that if $(\alpha,\beta)\in \Lambda^{\min}(\eta,\lambda(0,e_i))$, then 
\[
d(\nu\alpha)_i=d(\nu)_i+\max\{d(\eta)_i,d(\lambda(0,e_i))_i\}-d(\eta)_i=1
\]
and 
\[
d(\beta)_i=\max\{d(\eta)_i,d(\lambda(0,e_i))_i\}-d(\lambda(0,e_i))_i=0.
\]
Thus, since $\Omega_n$ is left $\mathcal{T}C^*(\Lambda^i)$-linear, we see that \eqref{calculations for left linearity of omega maps} is equal to
\begin{align*}
\sum_{(\alpha,\beta)\in \Lambda^{\min}(\eta,\lambda(0,e_i))}&t_{(\nu\alpha)(0,e_i)}^\Lambda \otimes_{\mathcal{T}C^*(\Lambda^i)} t_{(\nu\alpha)(e_i,d(\nu\alpha))}^{\Lambda^i} {t_\beta^{\Lambda^i}}^*\cdot \Omega_n\big(t_{\lambda(e_i,d(\lambda))}^\Lambda {t_\mu^\Lambda}^*\big)\\
&=
\sum_{(\alpha,\beta)\in \Lambda^{\min}(\eta,\lambda(0,e_i))}t_{(\nu\alpha)(0,e_i)}^\Lambda \otimes_{\mathcal{T}C^*(\Lambda^i)} \Omega_n\big(t_{(\nu\alpha)(e_i,d(\nu\alpha))}^\Lambda {t_\beta^\Lambda}^*t_{\lambda(e_i,d(\lambda))}^\Lambda {t_\mu^\Lambda}^*\big)\\
&=
\sum_{\substack{(\alpha,\beta)\in \Lambda^{\min}(\eta,\lambda(0,e_i))\\(\gamma,\delta)\in \Lambda^{\min}(\beta,\lambda(e_i,d(\lambda)))}}t_{(\nu\alpha)(0,e_i)}^\Lambda \otimes_{\mathcal{T}C^*(\Lambda^i)} \Omega_n\big(t_{(\nu\alpha)(e_i,d(\nu\alpha))\gamma}^\Lambda {t_{\mu\delta}^\Lambda}^*\big).
\end{align*}
Since $d(\nu\alpha)\geq e_i$, the factorisation property gives that
$
(\nu\alpha)(0,e_i)=(\nu\alpha\gamma)(0,e_i)
$
and
$
(\nu\alpha)(e_i,d(\nu\alpha))\gamma=(\nu\alpha\gamma)(e_i,d(\nu\alpha\gamma)).
$
Assembling these arguments and using Lemma~\ref{factorising minimal common extensions} for the third equality, we have
\begin{align*}
t_\nu^{\Lambda^i} {t_\eta^{\Lambda^i}}^*\cdot \Omega_{n+1}\big(t_\lambda^\Lambda{t_\mu^\Lambda}^*\big)
&=
\sum_{\substack{(\alpha,\beta)\in \Lambda^{\min}(\eta,\lambda(0,e_i))\\(\gamma,\delta)\in \Lambda^{\min}(\beta,\lambda(e_i,d(\lambda)))}}t_{(\nu\alpha\gamma)(0,e_i)}^\Lambda \otimes_{\mathcal{T}C^*(\Lambda^i)} \Omega_n\big(t_{(\nu\alpha\gamma)(e_i,d(\nu\alpha\gamma))}^\Lambda {t_{\mu\delta}^\Lambda}^*\big)\\
&=
\sum_{\substack{(\alpha,\beta)\in \Lambda^{\min}(\eta,\lambda(0,e_i))\\(\gamma,\delta)\in \Lambda^{\min}(\beta,\lambda(e_i,d(\lambda)))}}\Omega_{n+1}\big(t_{\nu\alpha\gamma}^\Lambda {t_{\mu\delta}^\Lambda}^*\big)\\
&=
\sum_{(\tau,\sigma)\in \Lambda^{\min}(\eta,\lambda)}\Omega_{n+1}\big(t_{\nu\tau}^\Lambda {t_{\mu\sigma}^\Lambda}^*\big)\\
&=
\Omega_{n+1}\big(t_\nu^{\Lambda^i} {t_\eta^{\Lambda^i}}^*\cdot t_\lambda^\Lambda{t_\mu^\Lambda}^*\big).
\end{align*}
Using the fact that
$
X_{n+1}=\cspan\{t_\lambda^\Lambda{t_\mu^\Lambda}^*:\lambda,\mu\in \Lambda, d(\lambda)_i=n+1,  d(\mu)_i=0\}
$
and 
$
\mathcal{T}C^*(\Lambda^i)=\cspan\{t_\nu^{\Lambda^i} {t_\eta^{\Lambda^i}}^*:\nu,\eta\in \Lambda^i\},
$
we conclude, by linearity and continuity, that $\Omega_{n+1}$ is left $\mathcal{T}C^*(\Lambda^i)$-linear.

Next, we show that $\Omega_{n+1}:X_{n+1}\rightarrow X^{\otimes n+1}$ is surjective. Fix $\lambda,\mu, \nu,\eta\in \Lambda$ with $d(\lambda)_i=1$, $d(\nu)_i=n$, and $d(\mu)_i=d(\eta)_i=0$. Using the left $\mathcal{T}C^*(\Lambda^i)$-linearity of $\Omega_n$ for the last equality, we see that
\begin{align*}
\Omega_{n+1}\bigg(\sum_{(\alpha,\beta)\in \Lambda^{\min}(\mu,\nu)}t_{\lambda\alpha}^\Lambda {t_{\eta\beta}^\Lambda}^*\bigg)
&=
\sum_{(\alpha,\beta)\in \Lambda^{\min}(\mu,\nu)}t_{(\lambda\alpha)(0,e_i)}^\Lambda\otimes_{\mathcal{T}C^*(\Lambda^i)} \Omega_n\big(t_{(\lambda\alpha)(e_i,d(\lambda\alpha))}^\Lambda {t_{\eta\beta}^\Lambda}^*\big)\\
&=
\sum_{(\alpha,\beta)\in \Lambda^{\min}(\mu,\nu)}t_{\lambda(0,e_i)}^\Lambda\otimes_{\mathcal{T}C^*(\Lambda^i)} \Omega_n\big(t_{\lambda(e_i,d(\lambda))\alpha}^\Lambda {t_{\eta\beta}^\Lambda}^*\big)\\
&=
t_{\lambda(0,e_i)}^\Lambda\otimes_{\mathcal{T}C^*(\Lambda^i)} \Omega_n\big(t_{\lambda(e_i,d(\lambda))}^\Lambda {t_\mu^\Lambda}^* t_\nu^\Lambda {t_\eta^\Lambda}^*\big)\\
&=
t_\lambda^\Lambda {t_\mu^\Lambda}^* \otimes_{\mathcal{T}C^*(\Lambda^i)} \Omega_n\big( t_\nu^\Lambda {t_\eta^\Lambda}^*\big)\\
&\in X\otimes_{\mathcal{T}C^*(\Lambda^i)} X^{\otimes n}.
\end{align*}
Since $X_m=\cspan\{t_\lambda^\Lambda{t_\mu^\Lambda}^*:\lambda,\mu\in \Lambda,  d(\lambda)_i=m,  d(\mu)_i=0\}$ for each $m\geq 0$ and the map $\Omega_n:X_n\rightarrow X^{\otimes n}$ is surjective, we conclude that $\Omega_{n+1}$ is surjective. 

We have now shown that $\Omega_{n+1}$ is inner product preserving and surjective. Thus, $\Omega_{n+1}$ is adjointable (with adjoint $\Omega_{n+1}^{-1}$). Since $\Omega_{n+1}$ is also left $\mathcal{T}C^*(\Lambda^i)$-linear, we conclude that $\Omega_{n+1}$ is a $\mathcal{T}C^*(\Lambda^i)$-bimodule isomorphism from $X_{n+1}$ to $X^{\otimes n+1}$ as required. 
\end{proof}
\end{prop}

We now work towards showing that the Toeplitz algebra of the Hilbert $\mathcal{T}C^*(\Lambda^i)$-bimodule $X$ is isomorphic to the Toeplitz--Cuntz--Krieger algebra of $\Lambda$. The idea is to use the universal properties of $\mathcal{T}_X$ and $\mathcal{T}C^*(\Lambda)$ to get $*$-homomorphisms between the two $C^*$-algebras, and then argue that these maps are mutually inverse. Firstly, we need a result telling us how the Hilbert $\mathcal{T}C^*(\Lambda^i)$-bimodule isomorphisms from Proposition~\ref{omega isomorphisms} interact with the tensor product. 

\begin{lem}
\label{bimodule isomorphisms and tensor products}
Let $\{\Omega_n:n\geq 0\}$ be the collection of Hilbert $\mathcal{T}C^*(\Lambda^i)$-bimodule isomorphisms defined in Proposition~\ref{omega isomorphisms}. Then for any $m,n\geq 0$ and $x\in X_m$, $y\in X_n$,
\begin{align}
\label{compatibility of Omega maps}
\Omega_{m}\left(x\right)\otimes_{\mathcal{T}C^*(\Lambda^i)} \Omega_{n}\left(y\right)
= \Omega_{m+n}\left(xy\right).
\end{align}
In particular, if $\lambda,\mu\in \Lambda$ with $r(\mu)=s(\lambda)$, then
\[
\Omega_{d(\lambda)_i}\big(t_\lambda^\Lambda\big)\otimes_{\mathcal{T}C^*(\Lambda^i)} \Omega_{d(\mu)_i}\big(t_\mu^\Lambda\big)
= \Omega_{d(\lambda\mu)_i}\big(t_{\lambda\mu}^\Lambda\big).
\]
\begin{proof}
We will use induction on $m$. The $m=0$ case is equivalent to the left $\mathcal{T}C^*(\Lambda^i)$-linearity of $\Omega_n$, which we proved in Proposition~\ref{omega isomorphisms}. Now suppose that \eqref{compatibility of Omega maps} holds for some $m\geq 0$. Let $n\geq 0$ and fix $\lambda,\mu,\nu,\tau\in \Lambda$ with $d(\lambda)_i=m+1$, $d(\nu)_i=n$, and $d(\mu)_i=d(\tau)_i=0$. Applying the inductive hypothesis, and using relation (TCK3) in $\mathcal{T}C^*(\Lambda)$ for the final equality, we see that
\begin{equation}
\label{omega map induction}
\begin{aligned}
\Omega_{m+1}\big(t_\lambda^\Lambda{t_\mu^\Lambda}^*\big)&\otimes_{\mathcal{T}C^*(\Lambda^i)} \Omega_{n}\big(t_\nu^\Lambda{t_\tau^\Lambda}^*\big)\\
&=
t_{\lambda(0,e_i)}^\Lambda\otimes_{\mathcal{T}C^*(\Lambda^i)}\Omega_{m}\big(t_{\lambda(e_i,d(\lambda))}^\Lambda{t_\mu^\Lambda}^*\big)\otimes_{\mathcal{T}C^*(\Lambda^i)} \Omega_{n}\big(t_\nu^\Lambda{t_\tau^\Lambda}^*\big)\\
&=
t_{\lambda(0,e_i)}^\Lambda\otimes_{\mathcal{T}C^*(\Lambda^i)} \Omega_{m+n}\big(t_{\lambda(e_i,d(\lambda))}^\Lambda{t_\mu^\Lambda}^*t_\nu^\Lambda{t_\tau^\Lambda}^*\big)\\
&=
\sum_{(\alpha,\beta)\in \Lambda^{\min}(\mu,\nu)}t_{\lambda(0,e_i)}^\Lambda\otimes_{\mathcal{T}C^*(\Lambda^i)} \Omega_{m+n}\big(t_{\lambda(e_i,d(\lambda))\alpha}^\Lambda{t_{\tau\beta}^\Lambda}^*\big).
\end{aligned}
\end{equation}
Since $d(\lambda)\geq (m+1)e_i\geq e_i$, \eqref{omega map induction} must be equal to
\begin{align*}
\sum_{(\alpha,\beta)\in \Lambda^{\min}(\mu,\nu)}t_{(\lambda\alpha)(0,e_i)}^\Lambda\otimes_{\mathcal{T}C^*(\Lambda^i)}\Omega_{m+n}\big(t_{(\lambda\alpha)(e_i,d(\lambda\alpha))}^\Lambda{t_{\tau\beta}^\Lambda}^*\big)
&=
\sum_{(\alpha,\beta)\in \Lambda^{\min}(\mu,\nu)}\Omega_{m+n+1}\big(t_{\lambda\alpha}^\Lambda{t_{\tau\beta}^\Lambda}^*\big)\\
&=
\Omega_{m+n+1}\big(t_\lambda^\Lambda{t_\mu^\Lambda}^*t_\nu^\Lambda{t_\tau^\Lambda}^*\big).
\end{align*}
Since $X_j=\cspan\{t_\lambda^\Lambda{t_\mu^\Lambda}^*:\lambda,\mu\in\Lambda, d(\lambda)_i=j, d(\mu)_i=0\}$ for each $j\geq 0$, we conclude that \eqref{compatibility of Omega maps} holds for $m+1$ as well. 
\end{proof}
\end{lem}

We now get a $*$-homomorphism from $\mathcal{T}C^*(\Lambda)$ to $\mathcal{T}_X$ by exhibiting a Toeplitz--Cuntz--Krieger $\Lambda$-family in $\mathcal{T}_X$.

\begin{prop}
\label{hom from TCK algebra to Toeplitz algebra}
Let $\Lambda$ be a finitely aligned $k$-graph. Define $X_n$ as in Proposition~\ref{module for TCK algebra} and set $X:=X_1$. Consider the collection of Hilbert $\mathcal{T}C^*(\Lambda^i)$-bimodule isomorphisms $\{\Omega_n:n\geq 0\}$ defined in Proposition~\ref{omega isomorphisms}. For each $\lambda\in \Lambda$, define $u_\lambda \in \mathcal{T}_X$ by
\[
u_\lambda:=i_X^{\otimes d(\lambda)_i}\left(\Omega_{d(\lambda)_i}\left(t_\lambda^\Lambda\right)\right).
\]
Then $\{u_\lambda:\lambda \in \Lambda\}$ is a Toeplitz--Cuntz--Krieger $\Lambda$-family in $\mathcal{T}_X$. Hence, there exists a $*$-homomorphism $\pi_u:\mathcal{T}C^*(\Lambda)\rightarrow \mathcal{T}_X$ such that $\pi_u\big(t_\lambda^\Lambda\big)=u_\lambda$ for each $\lambda\in \Lambda$.
\begin{proof}
Firstly, we check that $\{u_\lambda:\lambda \in \Lambda\}$ satisfies (TCK1). For any $v\in \Lambda^0$, we see that
\begin{align*}
u_v=i_X^{\otimes d(v)_i}\big(\Omega_{d(v)_i}\big(t_v^\Lambda\big)\big)=i_X^{\otimes 0}\big(\Omega_0\big(t_v^\Lambda\big)\big)
=i_{\mathcal{T}C^*(\Lambda^i)}\big(\phi^{-1}\big(t_v^\Lambda\big)\big)
=i_{\mathcal{T}C^*(\Lambda^i)}\big(t_v^{\Lambda^i}\big).
\end{align*}
Since $i_{\mathcal{T}C^*(\Lambda^i)}$ is a $*$-homomorphism and $\{t_v^{\Lambda^i}:v\in \Lambda^0\}$ is a collection of mutually orthogonal projections, it follows that the set $\{u_v:v\in \Lambda^0\}$ also consists of mutually orthogonal projections. 

Next we check that $\{u_\lambda:\lambda \in \Lambda\}$ satisfies (TCK2). Fix $\lambda,\mu\in \Lambda$ with $r(\mu)=s(\lambda)$. Making use of Lemma~\ref{bimodule isomorphisms and tensor products}, we see that
\begin{align*}
u_\lambda u_\mu
=i_X^{\otimes d(\lambda)_i}\big(\Omega_{d(\lambda)_i}\big(t_\lambda^\Lambda\big)\big)i_X^{\otimes d(\mu)_i}\big(\Omega_{d(\mu)_i}\big(t_\mu^\Lambda\big)\big)
&=i_X^{\otimes \left(d(\lambda)_i+d(\mu)_i\right)}\big(\Omega_{d(\lambda)_i}\big(t_\lambda^\Lambda\big)\hspace{-0.2em}\otimes_{\mathcal{T}C^*(\Lambda^i)}\hspace{-0.2em}\Omega_{d(\mu)_i}\big(t_\mu^\Lambda\big)\big)\\
&=i_X^{\otimes d(\lambda\mu)_i}\big(\Omega_{d(\lambda\mu)_i}\big(t_{\lambda\mu}^\Lambda\big)\big)\\
&=u_{\lambda\mu}.
\end{align*}

Finally we check that $\{u_\lambda:\lambda \in \Lambda\}$ satisfies (TCK3). Let $\lambda,\mu\in \Lambda$. Suppose that $d(\mu)_i\geq d(\lambda)_i$. By Lemma~\ref{bimodule isomorphisms and tensor products}
\begin{align*}
\Omega_{d(\mu)_i}\big(t_\mu^\Lambda\big)
=\Omega_{d(\lambda)_i}\big(t_{\mu(0,d(\lambda)_i e_i)}^\Lambda\big)\otimes_{\mathcal{T}C^*(\Lambda^i)}\Omega_{d(\mu)_i-d(\lambda)_i}\big(t_{\mu(d(\lambda)_i e_i,d(\mu))}^\Lambda\big),
\end{align*}
and it follows that
\begin{align*}
u_\lambda^* u_\mu
&=i_X^{\otimes d(\lambda)_i}\big(\Omega_{d(\lambda)_i}\big(t_\lambda^\Lambda\big)\big)^*i_X^{\otimes d(\mu)_i}\big(\Omega_{d(\mu)_i}\big(t_\mu^\Lambda\big)\big)\\
&=i_X^{\otimes d(\mu)_i-d(\lambda)_i}\big(\big\langle \Omega_{d(\lambda)_i}\big(t_\lambda^\Lambda\big),\Omega_{d(\lambda)_i}\big(t_{\mu(0,d(\lambda)_i e_i)}^\Lambda\big)\big\rangle \cdot\Omega_{d(\mu)_i-d(\lambda)_i}\big(t_{\mu(d(\lambda)_i e_i,d(\mu))}^\Lambda\big)\big).
\end{align*}
As $\Omega_{d(\lambda)_i}$ preserves inner products and $\Omega_{d(\mu)_i-d(\lambda)_i}$ is left $\mathcal{T}C^*(\Lambda^i)$-linear, this must be the same as
\begin{align*}
i_X^{\otimes d(\mu)_i-d(\lambda)_i}&\big(\Omega_{d(\mu)_i-d(\lambda)_i}\big(\big\langle t_\lambda^\Lambda,t_{\mu(0,d(\lambda)_i e_i)}^\Lambda\big\rangle^{d(\lambda)_i}_{\mathcal{T}C^*(\Lambda^i)}\cdot t_{\mu(d(\lambda)_i e_i,d(\mu))}^\Lambda\big)\big)\\
&=i_X^{\otimes d(\mu)_i-d(\lambda)_i}\big(\Omega_{d(\mu)_i-d(\lambda)_i}\big({t_\lambda^\Lambda}^*t_{\mu(0,d(\lambda)_i e_i)}^\Lambda t_{\mu(d(\lambda)_i e_i,d(\mu))}^\Lambda\big)\big)\\
&=i_X^{\otimes d(\mu)_i-d(\lambda)_i}\big(\Omega_{d(\mu)_i-d(\lambda)_i}\big({t_\lambda^\Lambda}^*t_\mu^\Lambda\big)\big)\\
&=i_X^{\otimes d(\mu)_i-d(\lambda)_i}\bigg(\Omega_{d(\mu)_i-d(\lambda)_i}\bigg(\sum_{(\alpha,\beta)\in \Lambda^{\min}(\lambda,\mu)}t_\alpha^\Lambda{t_\beta^\Lambda}^*\bigg)\bigg),
\end{align*}
where the last equality comes from the fact that $\{t_\lambda^\Lambda:\lambda\in \Lambda\}$ satisfies (TCK3). Moreover, if $(\alpha,\beta)\in \Lambda^{\min}(\lambda,\mu)$, then
\begin{align*}
d(\alpha)_i=\max\{d(\lambda)_i,d(\mu)i\}-d(\lambda)_i=d(\mu)_i-d(\lambda)_i
\end{align*}
and
\begin{align*}
d(\alpha)_i=\max\{d(\lambda)_i,d(\mu)i\}-d(\mu)_i=0.
\end{align*}
Thus, as $\Omega_{d(\mu)_i-d(\lambda)_i}$ is right $\mathcal{T}C^*(\Lambda^i)$-linear, 
\begin{align*}
u_\lambda^* u_\mu
&=\sum_{(\alpha,\beta)\in \Lambda^{\min}(\lambda,\mu)}i_X^{\otimes d(\mu)_i-d(\lambda)_i}\big(\Omega_{d(\mu)_i-d(\lambda)_i}\left(t_\alpha^\Lambda\right)\cdot {t_\beta^{\Lambda^i}}^*\big)\\
&=\sum_{(\alpha,\beta)\in \Lambda^{\min}(\lambda,\mu)}i_X^{\otimes d(\mu)_i-d(\lambda)_i}\big(\Omega_{d(\mu)_i-d(\lambda)_i}\big(t_\alpha^\Lambda\big)\big)i_X^{\otimes 0}\big({t_\beta^{\Lambda^i}}^*\big)\\
&=\sum_{(\alpha,\beta)\in \Lambda^{\min}(\lambda,\mu)}i_X^{\otimes d(\alpha)_i}\big(\Omega_{d(\alpha)_i)}\big(t_\alpha^\Lambda\big)\big)i_X^{\otimes d(\beta)_i}\big(\Omega_{d(\beta)_i}\big(t_\beta^\Lambda\big)\big)^*\\
&=\sum_{(\alpha,\beta)\in \Lambda^{\min}(\lambda,\mu)} u_\alpha u_\beta^*.
\end{align*}
If $d(\lambda)_i\geq d(\mu)_i$, we can apply the previous working to $(u_\lambda^* u_\mu)^*=u_\mu^* u_\lambda$. This completes the proof that $\{u_\lambda:\lambda\in \Lambda\}$ satisfies (TCK3). Hence, $\{u_\lambda:\lambda \in \Lambda\}$ is a Toeplitz--Cuntz--Krieger $\Lambda$-family in $\mathcal{T}_X$. The universal property of $\mathcal{T}C^*(\Lambda)$ then induces a $*$-homomorphism $\pi_u:\mathcal{T}C^*(\Lambda)\rightarrow \mathcal{T}_X$ such that $\pi_u\big(t_\lambda^\Lambda\big)=u_\lambda$ for each $\lambda\in \Lambda$.
\end{proof}
\end{prop}

It is considerably easier to get a $*$-homomorphism from $\mathcal{T}_X$ to $\mathcal{T}C^*(\Lambda)$. Once we have it, there is still some work left to show that it is the inverse of the $*$-homomorphism $\pi_u:\mathcal{T}C^*(\Lambda)\rightarrow \mathcal{T}_X$ from Proposition~\ref{hom from TCK algebra to Toeplitz algebra}. 

\begin{thm}
\label{Toeplitz algebra iterated construction}
Let $\Lambda$ be a finitely aligned $k$-graph. Define $X_n$ as in Proposition~\ref{module for TCK algebra} and set $X:=X_1$. Let  $\iota:X\rightarrow \mathcal{T}C^*(\Lambda)$ denote the inclusion map. Then $(\iota, \phi)$ is a Toeplitz representation of $X$ in $\mathcal{T}C^*(\Lambda)$, and hence, by the universal property of $\mathcal{T}_X$, there exists a $*$-homomorphism $\iota \times_\mathcal{T} \phi:\mathcal{T}_X\rightarrow \mathcal{T}C^*(\Lambda)$ such that $(\iota \times_\mathcal{T} \phi)\circ i_X=\iota$ and $(\iota \times_\mathcal{T} \phi)\circ i_{\mathcal{T}C^*(\Lambda^i)}=\phi$. Moreover, $\pi_u$ and $\iota \times_\mathcal{T} \phi$ are mutually inverse. Thus, $\mathcal{T}C^*(\Lambda)\cong \mathcal{T}_X$. 
\begin{proof}
It is elementary to check that $(\iota, \phi)$ is a Toeplitz representation of $X$ in $\mathcal{T}C^*(\Lambda)$. Observe that for any $x\in X$ and $a\in A$, we have
$
\iota(a\cdot x)=a\cdot x=\phi(a)x=\phi(a)\iota(x)
$
and 
$
\iota(x\cdot a)=x\cdot a=x\phi(a)=\iota(x)\phi(a),
$
which proves that $(\iota, \phi)$ satisfies (T1) and (T2). If $x,y\in X$, then
$
\iota(x)^*\iota(y)=x^*y=\phi\big(\phi^{-1}(x^*y)\big)=\phi\big(\langle x,y\rangle^1_{\mathcal{T}C^*(\Lambda^i)}\big),
$ 
and so $(\iota, \phi)$ satisfies (T3). 

It remains to check that $\iota \times_\mathcal{T} \phi$ and $\pi_u$ are mutually inverse. Fix $\lambda \in \Lambda$. If $d(\lambda)_i=0$, then
\begin{align*}
\left((\iota \times_\mathcal{T} \phi)\circ \pi_u\right) \big(t_\lambda^\Lambda\big)
&=(\iota \times_\mathcal{T} \phi)(u_\lambda)
=(\iota \times_\mathcal{T} \phi)\big(i_{X}^{\otimes d(\lambda)_i}\big(\Omega_{d(\lambda)_i}\big(t_\lambda^\Lambda\big)\big)\big)\\
&=(\iota \times_\mathcal{T} \phi)\big(i_{\mathcal{T}C^*(\Lambda^i)}\big(t_\lambda^{\Lambda^i}\big)\big)
=\phi\big(t_\lambda^{\Lambda^i}\big)
=t_\lambda^\Lambda.
\end{align*}
If $d(\lambda)_i=1$, then
\begin{align*}
\left((\iota \times_\mathcal{T} \phi)\circ \pi_u\right) \big(t_\lambda^\Lambda\big)
&=(\iota \times_\mathcal{T} \phi)(u_\lambda)
=(\iota \times_\mathcal{T} \phi)\big(i_{X}^{\otimes d(\lambda)_i}\big(\Omega_{d(\lambda)_i}\big(t_\lambda^\Lambda\big)\big)\big)\\
&=(\iota \times_\mathcal{T} \phi)\big(i_X\big(t_\lambda^\Lambda\big)\big)
=\iota\big(t_\lambda^\Lambda\big)
=t_\lambda^\Lambda.
\end{align*}
If $d(\lambda)_i\geq 2$, then
\begin{align*}
\left((\iota \times_\mathcal{T} \phi)\circ \pi_u\right) \big(t_\lambda^\Lambda\big)
=(\iota \times_\mathcal{T} \phi)(u_\lambda)
&=(\iota \times_\mathcal{T} \phi)\big(i_{X}^{\otimes d(\lambda)_i}\big(\Omega_{d(\lambda)_i}\big(t_\lambda^\Lambda\big)\big)\big)\\
&=(\iota \times_\mathcal{T} \phi)\big(i_X\big(t_{\lambda(0,e_i)}^\Lambda\big)\cdots i_X\big(t_{\lambda((d(\lambda)_i-1)e_i,d(\lambda))}^\Lambda\big)\big)\\
&=t_{\lambda(0,e_i)}^\Lambda\cdots t_{\lambda((d(\lambda)_i-1)e_i,d(\lambda))}^\Lambda
=t_\lambda^\Lambda.
\end{align*}
Since $\mathcal{T}C^*(\Lambda)$ is generated by $\big\{t_\lambda^\Lambda:\lambda\in \Lambda\big\}$, we conclude that $(\iota \times_\mathcal{T} \phi)\circ \pi_u = \text{id}_{\mathcal{T}C^*(\Lambda)}$. 

We now show that $\pi_u\circ (\iota \times_\mathcal{T} \phi)=\text{id}_{\mathcal{T}_X}$. If $\mu\in \Lambda^i$, then
\begin{align*}
\left(\pi_u\circ (\iota \times_\mathcal{T} \phi)\right)\big(i_{\mathcal{T}C^*(\Lambda^i)}\big(t_\mu^{\Lambda^i}\big)\big)
&=\pi_u\big(\phi\big(t_\mu^{\Lambda^i} \big)\big)
=\pi_u\big(t_\mu^\Lambda\big)
=u_\mu\\
&=i_X^{\otimes d(\lambda)_i}\big(\Omega_{d(\mu)_i}\big(t_\mu^\Lambda\big)\big)
=i_{\mathcal{T}C^*(\Lambda^i)}\big(t_\mu^{\Lambda^i}\big).
\end{align*}
For any $\lambda\in\Lambda$ with $d(\lambda)_i=1$ and $\mu\in \Lambda^i$, we see that
\begin{align*}
\left(\pi_u\circ (\iota \times_\mathcal{T} \phi)\right)\big(i_X\big(t_\lambda^\Lambda {t_\mu^\Lambda}^*\big)\big)
=\pi_u\big(t_\lambda^\Lambda {t_\mu^\Lambda}^*\big)
=u_\lambda u_\mu^*
&=i_X^{\otimes d(\lambda)_i}\big(\Omega_{d(\lambda)_i}\big(t_\lambda^\Lambda\big)\big) i_X^{\otimes d(\mu)_i}\big(\Omega_{d(\mu)_i}\big(t_\mu^\Lambda\big)\big)^*\\
&=i_X\big(t_\lambda^\Lambda\big) i_{\mathcal{T}C^*(\Lambda^i)}\big(t_\mu^{\Lambda^i}\big)^*
=i_X\big(t_\lambda^\Lambda {t_\mu^\Lambda}^*\big).
\end{align*}
Since $\mathcal{T}_X$ is generated by $i_X(X)\cup i_{\mathcal{T}C^*(\Lambda^i)}(\mathcal{T}C^*(\Lambda^i))$, whilst $\mathcal{T}C^*(\Lambda^i)$ is generated by $\{t_\mu^{\Lambda^i}:\mu\in\Lambda^i\}$ and $X=\cspan\{t_\lambda^\Lambda {t_{\mu}^\Lambda}^*:\lambda,\mu\in\Lambda, d(\lambda)_i=1, d(\mu)_i=0\}$, we conclude that $\pi_u\circ (\iota \times_\mathcal{T} \phi)=\text{id}_{\mathcal{T}_X}$. Thus, $\iota \times_\mathcal{T} \phi$ and $\pi_u$ are mutually inverse. 
\end{proof}
\end{thm}

\begin{cor}
Let $\Lambda$ be a finitely aligned $k$-graph. Then the $*$-homomorphism $\Phi:C_0(\Lambda)\rightarrow\mathcal{T}C^*(\Lambda)$ that sends $\delta_v$ to $t_v^\Lambda$ for each $v\in \Lambda^0$ induces a $KK$-equivalence between $C_0(\Lambda^0)$ and $\mathcal{T}C^*(\Lambda)$. 
\begin{proof}
We will use induction on $k$. If $k=0$, then the map $\Phi$ is an isomorphism between $C_0(\Lambda)$ and $\mathcal{T}C^*(\Lambda)$, and so of course gives a $KK$-equivalence. Suppose that the result holds for some $k\geq 0$ and let $\Lambda$ be a finitely aligned $(k+1)$-graph. Fix $i\in \{1,\ldots, k+1\}$ and let $\Phi':C_0(\Lambda)\rightarrow\mathcal{T}C^*(\Lambda^i)$ denote the $*$-homomorphism that sends $\delta_v$ to $t_v^{\Lambda^i}$ for each $v\in \Lambda^0$. By the inductive hypothesis, $\Phi'$ induces a $KK$-equivalence. Proposition~\ref{existence of phi for TCK algebra} gives a $*$-homomorphism $\phi:\mathcal{T}C^*(\Lambda^i)\rightarrow \mathcal{T}C^*(\Lambda)$ such that $\phi\big(t_v^{\Lambda^i}\big)=t_v^\Lambda$ for $v\in \Lambda^0$. By Proposition~\ref{module for TCK algebra} and Theorem~\ref{Toeplitz algebra iterated construction} there exists a Hilbert $\mathcal{T}C^*(\Lambda^i)$-bimodule $X$ and an isomorphism $\iota \times_\mathcal{T} \phi:\mathcal{T}_X\rightarrow \mathcal{T}C^*(\Lambda)$ such that $(\iota \times_\mathcal{T} \phi)\circ i_{\mathcal{T}C^*(\Lambda^i)}=\phi$. Since higher-rank graphs are, by definition, countable categories, $\mathcal{T}C^*(\Lambda^i)$ is separable and $X$ is countably generated as a right $\mathcal{T}C^*(\Lambda^i)$-module. Thus, by \cite[Theorem~4.4]{MR1426840}, the $*$-homomorphism $i_{\mathcal{T}C^*(\Lambda^i)}$ induces a $KK$-equivalence between $\mathcal{T}C^*(\Lambda^i)$ and $\mathcal{T}_X$. Hence, the $*$-homomorphism $\Phi:=\phi\circ \Phi'=(\iota \times_\mathcal{T} \phi)\circ i_{\mathcal{T}C^*(\Lambda^i)}\circ \Phi'$ induces a $KK$-equivalence between $C_0(\Lambda)$ and $\mathcal{T}C^*(\Lambda)$ and sends $\delta_v$ to $t_v^\Lambda$ for each $v\in \Lambda^0$.
\end{proof}
\end{cor}

\begin{rem}
Since $KK$-equivalent $C^*$-algebras have the same $K$-theory, we have an alternative proof of \cite[Theorem~1.1]{MR2498556} that $K_0(\mathcal{T}C^*(\Lambda))\cong \bigoplus_{v\in \Lambda^0}\mathbb{Z}$ and $K_1(\mathcal{T}C^*(\Lambda))\cong 0$ for any finitely aligned $k$-graph. 
\end{rem}

\section{Realising $C^*(\Lambda)$ as a Cuntz--Pimsner algebra}
\label{section: realising as Cuntz--Pimsner algebra}

In Section~\ref{section: realising as Toeplitz algebra} we showed how the Toeplitz--Cuntz--Krieger algebra of a finitely aligned $k$-graph $\Lambda$ can be realised as the Toeplitz algebra of a Hilbert $\mathcal{T}C^*(\Lambda^i)$-bimodule. In this section we prove an analogous result for Cuntz--Krieger algebras, defining a Hilbert $C^*(\Lambda^i)$-bimodule (which, for simplicity, we also denote by $X$) and showing that the Cuntz--Pimsner algebra of this bimodule is isomorphic to $C^*(\Lambda)$. 

Our methodology is very similar to that of Section~\ref{section: realising as Toeplitz algebra}.  Similar to Proposition~\ref{existence of phi for TCK algebra}, the first step is to show that the inclusion of $\Lambda^i$ in $\Lambda$ induces an (injective) $*$-homomorphism from $C^*(\Lambda^i)$ to $C^*(\Lambda)$.  However, unlike in Proposition~\ref{existence of phi for TCK algebra}, such a $*$-homomorphism need not exist unless we place additional constraints on the graph (see Remark~\ref{what goes wrong without local convexity} for an example of what can go wrong). In the analysis of \cite[Chapter~2]{fletcherphd}, to get around this problem we assumed that $\Lambda$ had no sources (see \cite[Proposition~2.6.4]{fletcherphd}). In Proposition~\ref{injective map for Cuntz--Krieger algebras finitely aligned, locally convex} we improve the situation, by showing that local-convexity of $\Lambda$ is sufficient. Before we prove Proposition~\ref{injective map for Cuntz--Krieger algebras finitely aligned, locally convex}, we prove some (probably) well-known properties of locally-convex higher-rank graphs that we will need. 

\begin{lem}
\label{extending local convexity}
Let $\Lambda$ be a locally-convex $k$-graph. If $\mu \in \Lambda^{e_i}$ and $\nu\in r(\mu)\Lambda$ with $d(\nu)_i=0$, then $s(\nu)\Lambda^{e_i}\neq \emptyset$. 
\begin{proof}
We use induction on the quantity $L(\nu):=\sum_{j=1}^k d(\nu)_j$. If $L(\nu)=0$, then $\nu\in \Lambda^0$ and so $\mu \in r(\nu)\Lambda^{e_i}=s(\nu)\Lambda^{e_i}$. Suppose $M\in \N\cup \{0\}$ and the result holds whenever $L(\nu)=M$. Fix $\nu'\in r(\mu)\Lambda$ with $d(\nu')_i=0$ and $L(\nu')=M+1$. Thus for some $l\in \{1,\ldots, k\}\setminus \{i\}$, we have $d(\nu')_l\geq 1$. Then $L(\nu'(0,d(\nu')-e_l))=M$, and so $s(\nu'(0,d(\nu')-e_l))\Lambda^{e_i}$ is nonempty by the inductive hypothesis. As $\nu'(d(\nu')-e_l, d(\nu'))\in s(\nu'(0,d(\nu')-e_l))\Lambda^{e_l}$ and $\Lambda$ is locally-convex, we have that $s(\nu')\Lambda^{e_i}=s(\nu'(d(\nu')-e_l, d(\nu')))\Lambda^{e_i}\neq \emptyset$ as required. 
\end{proof}
\end{lem}

\begin{lem}
\label{factorisation in locally convex graphs}
Let $\Lambda$ be a $k$-graph. Then $\Lambda^{\leq m}\Lambda^{\leq n}\subseteq \Lambda^{\leq m+n}$ for each $m,n\in \N^k$. If $\Lambda$ is locally-convex then $\Lambda^{\leq m}\Lambda^{\leq n}= \Lambda^{\leq m+n}$.
\begin{proof}
Firstly, suppose that $\mu\in \Lambda^{\leq m}$ and $\nu\in \Lambda^{\leq n}$ with $s(\mu)=r(\nu)$. Clearly, $d(\mu\nu)=d(\mu)+d(\nu)\leq m+n$. Suppose that $d(\mu\nu)_i< (m+n)_i=m_i+n_i$ for some $i$. If $d(\nu)_i<n_i$, then $s(\mu\nu)\Lambda^{e_i}=s(\nu)\Lambda^{e_i}=\emptyset$ since $\nu\in \Lambda^{\leq n}$. Thus, $\mu\nu\in \Lambda^{\leq m+n}$. Alternatively, $d(\mu)_i<m_i$, and so $s(\mu)\Lambda^{e_i}=\emptyset$ since $\mu\in \Lambda^{\leq m}$. Then by the factorisation property, $s(\mu\nu)\Lambda^{e_i}=s(\nu)\Lambda^{e_i}=\emptyset$. Thus, $\mu\nu\in \Lambda^{\leq m+n}$.

Now suppose that $\Lambda$ is locally-convex. We need to show that $\Lambda^{\leq m+n}\subseteq \Lambda^{\leq m}\Lambda^{\leq n}$. Fix $\lambda\in \Lambda^{\leq m+n}$. Let $m':=m\wedge d(\lambda)$ be the component-wise minimum of $m$ and $d(\lambda)$, and set $\mu:=\lambda(0,m')$ and $\nu:=\lambda(m', d(\lambda))$. Clearly, $\lambda=\mu\nu$. We claim that $\mu\in \Lambda^{\leq m}$ and $\nu\in \Lambda^{\leq n}$. Obviously, $d(\mu)\leq m$, and routine calculations show that $d(\nu)\leq n$. Suppose that $d(\nu)_i<n_i$ for some $i$. Then $d(\lambda)_i<(m'+n)_i\leq (m+n)_i$ and so $s(\nu)\Lambda^{e_i}=s(\lambda)\Lambda^{e_i}=\emptyset$. Thus, $\nu\in \Lambda^{\leq n}$. Now suppose that $d(\mu)_i< m_i$ for some $i$. Hence, $m'_i<m_i$, and so $m'_i=d(\lambda)_i=d(\mu)_i$ and $d(\nu)_i=0$. Also, $d(\lambda)_i=m'_i<m_i\leq (m+n)_i$, and so $s(\nu)\Lambda^{e_i}=s(\lambda)\Lambda^{e_i}= \emptyset$. By Lemma~\ref{extending local convexity}, this forces $s(\mu)\Lambda^{e_i}=r(\nu)\Lambda^{e_i}=\emptyset$. Thus, $\mu\in \Lambda^{\leq m}$.
\end{proof}
\end{lem}

We now work towards showing that the inclusion of $\Lambda^i$ in $\Lambda$ induces a $*$-homomorphism from $C^*(\Lambda^i)$ to $C^*(\Lambda)$. The key point is that when $\Lambda$ is locally-convex, finite-exhaustive subsets of $\Lambda^i$ are also exhaustive in $\Lambda$. 

\begin{defn}
Let $\Lambda$ be a $k$-graph. For any $E\subseteq \Lambda$ and $\mu\in \Lambda$, we define
\[
\mathrm{Ext}_\Lambda(\mu;E):=\bigcup_{\lambda\in E}\{\alpha\in s(\mu)\Lambda:\mu\alpha\in\mathrm{MCE}(\mu,\lambda)\}.
\]
\end{defn}

Informally speaking, $\mathrm{Ext}_\Lambda(\mu;E)$ is the set of paths in $\Lambda$ that when prepended to $\mu$ give a minimal common extension of $\mu$ with something in $E$. 

\begin{lem}[\cite{MR2069786}, Lemma C.5]
\label{extending finite exhaustive sets}
Let $\Lambda$ be a finitely aligned $k$-graph. Fix $v\in \Lambda^0$ and let $E\subseteq v\Lambda$ be a finite exhaustive set in $\Lambda$. Then for any $\mu\in v\Lambda$, the set $\mathrm{Ext}_\Lambda(\mu;E)\subseteq s(\mu)\Lambda$ is finite and exhaustive in $\Lambda$. 
\begin{proof}
Firstly, we check that $\mathrm{Ext}_\Lambda(\mu;E)$ is finite. For each $\lambda\in E$,  since $\Lambda$ is finitely aligned, the set $\{\alpha\in \Lambda:\mu\alpha\in\mathrm{MCE}(\mu,\lambda)\}$ is finite. As $E$ is finite, $\mathrm{Ext}_\Lambda(\mu;E)$ is the finite union of finite sets, and so finite. It remains to verify that $\mathrm{Ext}_\Lambda(\mu;E)$ is exhaustive in $\Lambda$. Fix $\sigma\in s(\mu)\Lambda$. Since $\mu\sigma \in v\Lambda$ and $E\subseteq v\Lambda$ is exhaustive in $\Lambda$, there exists $\lambda\in E$ and $\alpha,\beta\in \Lambda$ such that $\mu\sigma\alpha=\lambda\beta\in \mathrm{MCE}(\lambda, \mu\sigma)$. Let 
$
\tau:=(\sigma\alpha)(0, d(\lambda)\vee d(\mu)-d(\mu)),
$
which is well-defined because
\[
d(\lambda)\vee d(\mu)-d(\mu)\leq d(\lambda)\vee d(\mu\sigma)-d(\mu)=d(\mu\sigma\alpha)-d(\mu)=d(\sigma\alpha).
\]
Then
\begin{align*}
\mu\tau
&=\mu(\sigma\alpha)(0, d(\lambda)\vee d(\mu)-d(\mu))
=(\mu\sigma\alpha)(0, d(\lambda)\vee d(\mu))\\
&=(\lambda\beta)(0, d(\lambda)\vee d(\mu))
=\lambda\beta(0, d(\lambda)\vee d(\mu)-d(\lambda))
\in \Lambda^{d(\lambda)\vee d(\mu)},
\end{align*}
which shows that $\mu\tau\in \mathrm{MCE}(\mu,\lambda)$. As $\lambda\in E$, we see that $\tau\in \mathrm{Ext}_\Lambda(\mu;E)$. Furthermore,
\begin{align*}
\tau(\sigma\alpha)(d(\lambda)\vee d(\mu)-d(\mu),d(\sigma\mu))=\sigma\alpha,
\end{align*}
which shows that $\mathrm{CE}(\tau,\sigma)\neq \emptyset$, and so $\mathrm{MCE}(\tau,\sigma)\neq \emptyset$. Therefore, $\mathrm{Ext}_\Lambda(\mu;E)$ is exhaustive in $\Lambda$.
\end{proof}
\end{lem}

\begin{lem}
\label{finite exhaustive sets in subgraphs are finite exhaustive sets in whole graph - locally convex version}
Let $\Lambda$ be a locally-convex $k$-graph. Then $\mathrm{FE}(\Lambda^i)\subseteq\mathrm{FE}(\Lambda)$.  
\begin{proof}
We need to show that if $E\in\mathrm{FE}(\Lambda^i)$, then $E$ is exhaustive in $\Lambda$. Fix $\lambda\in r(E)\Lambda$ and write $\lambda=\lambda'\lambda_i$ with $\lambda'\in \Lambda^i$ and $\lambda_i\in \Lambda^{\N e_i}$. Let 
$
N:=\bigvee\{d(\mu):\mu\in \mathrm{Ext}_{\Lambda^i}(\lambda';E)\},
$
which exists since $\mathrm{Ext}_{\Lambda^i}(\lambda';E)$ is finite by Lemma~\ref{extending finite exhaustive sets}. Since $N_i=0$, we can choose $\tau\in s(\lambda_i)\Lambda^{\leq N}\subseteq \Lambda^i$. Thus, $\lambda_i\tau\in \Lambda^{\leq d(\lambda_i)}\Lambda^{\leq N}\subseteq \Lambda^{\leq d(\lambda_i)+N}$. Since $\Lambda$ is locally-convex, Lemma~\ref{factorisation in locally convex graphs} says that we can find $\tau'\in \Lambda^{\leq N}\subseteq \Lambda^i$ and $\lambda_i'\in \Lambda^{\leq  d(\lambda_i)}\subseteq \Lambda^{\N e_i}$ such that $\lambda_i\tau=\tau'\lambda_i'$. Since $r(\tau')=r(\lambda_i)=s(\lambda')$ and $\mathrm{Ext}_{\Lambda^i}(\lambda';E)\subseteq s(\lambda')\Lambda^i$ is exhaustive in $\Lambda^i$ by Lemma~\ref{extending finite exhaustive sets}, there exists $\mu\in \mathrm{Ext}_{\Lambda^i}(\lambda';E)$ such that $\mathrm{MCE}(\mu,\tau')\neq \emptyset$. That is, we can find $\alpha,\beta \in \Lambda^i$ such that $\tau'\alpha=\mu\beta \in \Lambda^{d(\mu)\vee d(\tau')}$. As $N$ is maximal, $N\geq d(\mu)$, and so $d(\mu)\vee d(\tau')\leq N$.  Since $\tau'\in\Lambda^{\leq N}$, this forces $\alpha=s(\tau')$, and so $\tau'=\mu\beta$. Moreover, since $\mu\in \mathrm{Ext}_{\Lambda^i}(\lambda';E)$, we know that $\lambda'\mu=\sigma\xi\in \mathrm{MCE}(\lambda',\sigma)$ for some $\sigma\in E$ and $\xi\in \Lambda^i$. Therefore,
\[
\sigma\xi\beta\lambda_i'=\lambda'\mu\beta\lambda_i'=\lambda'\tau'\lambda_i'=\lambda'\lambda_i\tau=\lambda\tau.
\]
Thus, $\mathrm{CE}(\sigma, \lambda)\neq \emptyset$, and so $\mathrm{MCE}(\sigma, \lambda)\neq \emptyset$. As $\sigma\in E$, we conclude that $E$ is exhaustive in $\Lambda$.
\end{proof}
\end{lem}

\begin{prop}
\label{injective map for Cuntz--Krieger algebras finitely aligned, locally convex}
Let $\Lambda$ be a finitely aligned locally-convex $k$-graph. Then there exists an injective $*$-homomorphism $\phi:C^*(\Lambda^i)\rightarrow C^*(\Lambda)$ carrying $s_\lambda^{\Lambda^i}$ to $s_\lambda^\Lambda$ for each $\lambda\in \Lambda^i$.
\begin{proof}
We claim that $\{s_\lambda^\Lambda:\lambda\in \Lambda^i\}\subseteq C^*(\Lambda)$ is a Cuntz--Krieger $\Lambda^i$-family. The same argument as in the proof of Proposition~\ref{existence of phi for TCK algebra} shows that $\{s_\lambda^\Lambda:\lambda\in \Lambda^i\}$ satisfies (TCK1), (TCK2), and (TCK3), so we need only worry about checking that relation (CK) holds. With this in mind, fix $v\in (\Lambda^i)^0=\Lambda^0$ and suppose that $E\in v\mathrm{FE}(\Lambda^i)$ . By Lemma~\ref{finite exhaustive sets in subgraphs are finite exhaustive sets in whole graph - locally convex version}, $E$ is exhaustive in $\Lambda$. As $\{s_\lambda^\Lambda:\lambda\in \Lambda\}$ satisfies relation (CK), we conclude that $\{s_\lambda^\Lambda:\lambda\in \Lambda^i\}$ does as well. The universal property of $C^*(\Lambda^i)$ then induces a $*$-homomorphism $\phi$ from $C^*(\Lambda^i)$ to $C^*(\Lambda)$ such that $\phi(s_\lambda^{\Lambda^i})=s_\lambda^\Lambda$ for each $\lambda\in \Lambda^i$.
\\ \indent
The injectivity of $\phi$ follows from an application of \cite[Theorem~4.2]{MR2069786}. For each $v\in \Lambda^0$, we have $\phi(s_v^{\Lambda^i})=s_v^\Lambda$ which is nonzero by \cite[Proposition~2.12]{MR2069786}. Restricting the gauge action $\gamma^\Lambda$ of $\mathbb{T}^k$ on $C^*(\Lambda)$ to $\mathbb{T}^{k-1}$ gives an action of $\mathbb{T}^{k-1}$ on $C^*(\{\phi(s_\lambda^{\Lambda^i}):\lambda\in \Lambda^i\})=C^*(\{s_\lambda^\Lambda:\lambda\in \Lambda^i\})\subseteq C^*(\Lambda)$ that intertwines $\phi$ and the gauge action $\gamma^{\Lambda^i}$ of $\T^{k-1}$ on $C^*(\Lambda^i)$.  
\end{proof}
\end{prop}

\begin{rem}
\label{what goes wrong without local convexity}
There are simple examples to show what can go wrong if we do not have a locally-convex graph. 
Consider the $2$-graph $\Lambda$ consisting of just two edges $\lambda\in \Lambda^{e_1}$ and $\mu \in \Lambda^{e_2}$ with common range $v$ and distinct sources. In this situation the second part of Lemma~\ref{factorisation in locally convex graphs} is false: the path $\lambda\in \Lambda^{\leq e_1+e_2}$ cannot be written in the form $\eta\nu$ where $\eta \in \Lambda^{\leq e_2}$ and $\nu \in \Lambda^{\leq e_1}$ (due to the presence of the edge $\mu$, the vertex $v$ is not in $\Lambda^{\leq e_2}$). Furthermore, $\{\lambda\}$ is exhaustive in $\Lambda^2$, but not exhaustive in $\Lambda$ since $\mathrm{MCE}(\lambda,\mu)=\emptyset$. Thus the conclusion of Lemma~\ref{finite exhaustive sets in subgraphs are finite exhaustive sets in whole graph - locally convex version} need not hold if we drop the local-convexity hypothesis. This example also shows that the conclusion of Proposition~\ref{injective map for Cuntz--Krieger algebras finitely aligned, locally convex} is false if we drop the local-convexity hypothesis. The Cuntz--Krieger relation in $C^*(\Lambda^i)$ says that $s_v^{\Lambda^i}=s_\lambda^{\Lambda^i}{s_\lambda^{\Lambda^i}}^*$. On the other hand, the Cuntz--Krieger relation in $C^*(\Lambda)$ (applied to the finite exhaustive set $\{\lambda, \mu\}$) gives
\[
0=\big(s_v^\Lambda-s_\lambda^\Lambda{s_\lambda^\Lambda}^*\big)\big(s_v^\Lambda-s_\mu^\Lambda{s_\mu^\Lambda}^*\big)
=s_v^\Lambda-s_\lambda^\Lambda{s_\lambda^\Lambda}^*-s_\mu^\Lambda{s_\mu^\Lambda}^*
\]
since $\lambda$ and $\mu$ have no common extensions. Hence, if there existed a $*$-homomorphism $\phi$ from $C^*(\Lambda^i)$ to $C^*(\Lambda)$ induced by the inclusion of $\Lambda^i$ in $\Lambda$, we would have that
\[
0=\phi(0)=\phi\big(s_v^{\Lambda^i}-s_\lambda^{\Lambda^i}{s_\lambda^{\Lambda^i}}^*\big)
=s_v^\Lambda-s_\lambda^\Lambda{s_\lambda^\Lambda}^*=s_\mu^\Lambda{s_\mu^\Lambda}^*.
\]
Thus, $s_\mu^\Lambda=0$, which is impossible since universal Cuntz--Krieger families always consist of non-zero partial isometries \cite[Proposition~2.12]{MR2069786}. 
\end{rem}

We are now ready to define the collection of Hilbert $C^*(\Lambda^i)$-bimodules that we are interested in. Suppose that $\Lambda$ is locally-convex so that the injective $*$-homomorphism $\phi$ from Proposition~\ref{injective map for Cuntz--Krieger algebras finitely aligned, locally convex} exists. The same working as in Proposition~\ref{module for TCK algebra} shows that 
\[
X_n:=\cspan\{s_\lambda^\Lambda {s_\mu^\Lambda}^*:\lambda,\mu\in \Lambda, d(\lambda)_i=n, d(\mu)_i=0\}\subseteq C^*(\Lambda)
\]
has the structure of a Hilbert $C^*(\Lambda^i)$-bimodule for each $n\in \N\cup \{0\}$, with actions and inner product given by $a\cdot x\cdot b:=\phi(a)x\phi(b)$ and $\langle x,y \rangle^n_{C^*(\Lambda^i)}:=\phi^{-1}(x^*y)$ for each $x, y\in X_n$ and $a,b\in C^*(\Lambda^i)$. For notational convenience, we set $X:=X_1$, and write $\langle \cdot, \cdot  \rangle_{C^*(\Lambda^i)}$ for $\langle \cdot, \cdot  \rangle^1_{C^*(\Lambda^i)}$. We again write $\psi$ for the $*$-homomorphism that implements the left action of $C^*(\Lambda^i)$ on $X$. For each $n\in \N\cup \{0\}$, the same working as in Proposition~\ref{omega isomorphisms} gives a Hilbert $C^*(\Lambda^i)$-bimodule isomorphism $\Omega_n:X_n\rightarrow X^{\otimes n}$, where, in particular, $\Omega_0=\phi^{-1}$ and $\Omega_1$ is the identity map. 

Our goal is to now analyse the Cuntz--Pimsner algebra of $X$. In order to do this, we need to get a grip on the Katsura ideal $J_X:=\psi^{-1}(\mathcal{K}_{C^*(\Lambda^i)}(X))\cap \ker(\psi)^\perp$. In \cite[Lemma 2.6.7]{fletcherphd} we showed that if $\Lambda$ has no sources, then $C^*(\Lambda^i)$ acts faithfully on $X$.  Lemma~2.6.8 of \cite{fletcherphd} also shows that if $v\Lambda^{e_i}$ is finite for each $v\in \Lambda^0$, then $C^*(\Lambda^i)$ acts compactly on $X$. Thus, when $\Lambda$ is row finite and has no sources the Katsura ideal of $X$ is all of $C^*(\Lambda^i)$. Consequently, to determine whether a Toeplitz representation of $X$ was Cuntz--Pimsner covariant, we only needed to check that the covariance relation held on the generating set $\{s_\lambda^{\Lambda^i}:\lambda\in \Lambda\}$ of $C^*(\Lambda^i)$ (see \cite[Theorem~2.6.12]{fletcherphd}). In this paper we are not assuming that the graph $\Lambda$ is source free and row finite (recall, our only assumption so far is that $\Lambda$ is locally-convex), and so it is not immediately obvious what the Katsura ideal looks like, and whether it has a `nice' generating set that is easy to work with. Our strategy is to show that $J_X$ is gauge-invariant and calculate its generators using \cite[Theorem~4.6]{MR3262073}. We begin in Proposition~\ref{checking gauge invariant ideals} by showing that the ideals $\mathrm{ker}(\psi)$, $\mathrm{ker}(\psi)^\perp$, and $\psi^{-1}(\mathcal{K}_{C^*(\Lambda^i)}(X))$ are all gauge-invariant. First, we require a lemma. 

\begin{lem}
\label{existence of U_z}
Let $\Lambda$ be a locally-convex finitely aligned $k$-graph. For each $z\in \mathbb{T}^{k-1}$, there exists a unitary $U_z\in \mathcal{L}_{C^*(\Lambda^i)}(X)$ such that 
\[
U_z\big(s_\lambda^\Lambda \phi(a)\big)=s_\lambda^\Lambda\phi\big(\gamma_z^{\Lambda^i}(a)\big)
\]
for each $\lambda\in \Lambda^{e_i}$ and $a\in C^*(\Lambda^i)$. Moreover, for each $a\in C^*(\Lambda^i)$,
\[
\psi\big(\gamma_z^{\Lambda^i}(a)\big)=U_z\psi(a)U_z^*.
\]
\begin{proof}
We show that for each $z\in \mathbb{T}^{k-1}$, the formula $s_\lambda^\Lambda\phi(a)\mapsto s_\lambda^\Lambda\phi\big(\gamma^{\Lambda^i}_z(a)\big)$, where $\lambda\in \Lambda^{e_i}$, $a\in C^*(\Lambda^i)$,
extends by linearity and continuity to a map on $X$. Let $m\in \mathbb{N}$ and fix $\lambda_1,\ldots, \lambda_m\in \Lambda^{e_i}$ and $a_1,\ldots, a_m\in C^*(\Lambda^i)$. Then 
\begin{align*}
\bigg\|\sum_{j=1}^m s_{\lambda_j}^\Lambda\phi\big(\gamma^{\Lambda^i}_z(a_j)\big)\bigg\|_X^2
&=\bigg\|\bigg\langle \sum_{j=1}^m s_{\lambda_j}^\Lambda\phi\big(\gamma^{\Lambda^i}_z(a_j)\big), \sum_{j=1}^m s_{\lambda_j}^\Lambda\phi\big(\gamma^{\Lambda^i}_z(a_j)\big)\bigg\rangle_{C^*(\Lambda^i)}\bigg\|_{C^*(\Lambda^i)}\\
&=\bigg\|\sum_{j,l=1}^m\phi^{-1}\big(\big(s_{\lambda_j}^\Lambda\phi\big(\gamma^{\Lambda^i}_z(a_j)\big)\big)^*\big(s_{\lambda_l}^\Lambda\phi\big(\gamma^{\Lambda^i}_z(a_l)\big)\big)\big)\bigg\|_{C^*(\Lambda^i)}\\
&=\bigg\|\sum_{j,l=1}^m \gamma^{\Lambda^i}_z(a_j)^*\phi^{-1}\big({s_{\lambda_j}^\Lambda}^*s_{\lambda_l}^\Lambda\big)\gamma^{\Lambda^i}_z(a_l)\bigg\|_{C^*(\Lambda^i)}.
\end{align*}
Since $\lambda_1,\ldots, \lambda_m$ all have degree $e_i$, relation (TCK3) implies that ${s_{\lambda_j}^\Lambda}^*s_{\lambda_l}^\Lambda=\delta_{\lambda_j,\lambda_l}s_{s(\lambda_j)}^\Lambda$, and so the previous line is equal to
\begin{align*}
\bigg\|\sum_{j,l=1}^m \delta_{\lambda_j,\lambda_l}\gamma^{\Lambda^i}_z(a_j^*)s_{s(\lambda_j)}^{\Lambda^i}\gamma^{\Lambda^i}_z(a_l) \bigg\|_{C^*(\Lambda^i)}
&=\bigg\|\sum_{j,l=1}^m \delta_{\lambda_j,\lambda_l}\gamma^{\Lambda^i}_z(a_j^*)\gamma^{\Lambda^i}_z\big(s_{s(\lambda_j)}^{\Lambda^i}\big)\gamma^{\Lambda^i}_z(a_l) \bigg\|_{C^*(\Lambda^i)}\\
&=\bigg\|\gamma^{\Lambda^i}_z\bigg(\sum_{j,l=1}^m \delta_{\lambda_j,\lambda_l}a_j^*s_{s(\lambda_j)}^{\Lambda^i}a_l\bigg) \bigg\|_{C^*(\Lambda^i)}\\
&=\bigg\|\sum_{j,l=1}^m \delta_{\lambda_j,\lambda_l} a_j^*s_{s(\lambda_j)}^{\Lambda^i}a_l \bigg\|_{C^*(\Lambda^i)},
\end{align*}
where the last equality follows from the fact that $\gamma_z$ is an automorphism, and hence isometric. Finally, this is the same as
\begin{align*}
\bigg\|\sum_{j,l=1}^m a_j^*\phi^{-1}\big({s_{\lambda_j}^\Lambda}^*s_{\lambda_l}^\Lambda\big)a_l\bigg\|_{C^*(\Lambda^i)}
&=\bigg\|\bigg\langle \sum_{j=1}^m s_{\lambda_j}^\Lambda\phi(a_j), \sum_{j=1}^m s_{\lambda_j}^\Lambda\phi(a_j)\bigg\rangle_{C^*(\Lambda^i)}\bigg\|_{C^*(\Lambda^i)}\\
&=\bigg\|\sum_{j=1}^m s_{\lambda_j}^\Lambda\phi(a_j)\bigg\|_X^2.
\end{align*}
Thus, the formula $s_\lambda^\Lambda\phi(a)\mapsto s_\lambda^\Lambda\phi\big(\gamma^{\Lambda^i}_z(a)\big)$ extends by linearity and continuity to an inner product preserving map on $X$, which we denote by $U_z$. The map $U_z$ is surjective since, for any $\lambda\in \Lambda^{e_i}$ and $a\in C^*(\Lambda^i)$, we have $U_z\big(s_\lambda^\Lambda\phi\big(\gamma^{\Lambda^i}_{\overline{z}}(a)\big)\big)=s_\lambda^\Lambda\phi(a)$. Consequently, $U_z\in \mathcal{L}_{C^*(\Lambda^i)}(X)$, with $U_z^*=U_z^{-1}=U_{\overline{z}}$. 

It remains to check that for each $z\in \mathbb{T}^{k-1}$ and $a\in C^*(\Lambda^i)$, we have
\[
U_z\psi(a)U_z^*=\psi\big(\gamma^{\Lambda^i}_z(a)\big)\in \psi\big(C^*(\Lambda^i)\big)\subseteq \mathcal{L}_{C^*(\Lambda^i)}(X).
\]
To see this, fix $\eta,\rho, \nu, \tau \in \Lambda^i$ and $\lambda\in \Lambda^{e_i}$. Then 
\begin{align*}
\psi\big(s_\eta^{\Lambda^i} {s_\rho^{\Lambda^i}}^*\big)\big(s_\lambda^\Lambda\phi\big(s_\nu^{\Lambda^i} {s_\tau^{\Lambda^i}}^*\big)\big)
&=s_\eta^\Lambda {s_\rho^\Lambda}^*s_\lambda^\Lambda s_\nu^\Lambda {t_\tau^\Lambda}^*
=\sum_{(\alpha,\beta)\in \Lambda^{\min}(\rho,\lambda\nu)}s_{\eta\alpha}^\Lambda {s_{\tau\beta}^\Lambda}^*\\
&=\sum_{(\alpha,\beta)\in \Lambda^{\min}(\rho,\lambda\nu)}s_{(\eta\alpha)(0,e_i)}^\Lambda \phi\big(s_{(\eta\alpha)(e_i,d(\eta\alpha))}^{\Lambda^i}{s_{\tau\beta}^{\Lambda^i}}^*\big).
\end{align*}
Hence,
\begin{align*}
\big(U_z&\psi\big(s_\eta^{\Lambda^i} {s_\rho^{\Lambda^i}}^*\big)U_z^*\big)\big(s_\lambda^\Lambda\phi\big(s_\nu^{\Lambda^i} {s_\tau^{\Lambda^i}}^*\big)\big)\\
&=U_z\psi\big(s_\eta^{\Lambda^i} {s_\rho^{\Lambda^i}}^*\big)\big(s_\lambda^\Lambda\phi\big(\gamma^{\Lambda^i}_{\overline{z}}\big(s_\nu^{\Lambda^i} {s_\tau^{\Lambda^i}}^*\big)\big)\big)\\
&=z^{d(\tau)-d(\nu)}U_z\bigg(\sum_{(\alpha,\beta)\in \Lambda^{\min}(\rho,\lambda\nu)}s_{(\eta\alpha)(0,e_i)}^\Lambda \phi\big(s_{(\eta\alpha)(e_i,d(\eta\alpha))}^{\Lambda^i}{s_{\tau\beta}^{\Lambda^i}}^*\big)\bigg)\\
&=z^{d(\tau)-d(\nu)}\sum_{(\alpha,\beta)\in \Lambda^{\min}(\rho,\lambda\nu)}s_{(\eta\alpha)(0,e_i)}^\Lambda \phi\big(\gamma^{\Lambda^i}_z\big(s_{(\eta\alpha)(e_i,d(\eta\alpha))}^{\Lambda^i}{s_{\tau\beta}^{\Lambda^i}}^*\big)\big)\\
&=\sum_{(\alpha,\beta)\in \Lambda^{\min}(\rho,\lambda\nu)}z^{d(\tau)-d(\nu)+d((\eta\alpha)(e_i,d(\eta\alpha)))-d(\tau\beta)}s_{(\eta\alpha)(0,e_i)}^\Lambda \phi\big(s_{(\eta\alpha)(e_i,d(\eta\alpha))}^{\Lambda^i}{s_{\tau\beta}^{\Lambda^i}}^*\big)
\end{align*}
However, if $(\alpha,\beta)\in \Lambda^{\min}(\rho,\lambda\nu)$, then
\begin{align*}
d(\tau)-d(\nu)+d((\eta\alpha)(e_i,d(\eta\alpha)))-d(\tau\beta)
&=d(\eta)+d(\alpha)-e_i-d(\nu)-d(\beta)\\
&=d(\eta)+d(\rho)\vee d(\lambda\nu)-d(\rho)-e_i\\
&\qquad\qquad -d(\nu)-d(\rho)\vee d(\lambda\nu)+ d(\lambda\nu)\\
&=d(\eta)-d(\rho)-e_i+ d(\lambda)\\
&=d(\eta)-d(\rho). 
\end{align*}
Therefore,
\begin{align*}
\big(U_z\psi\big(s_\eta^{\Lambda^i} {s_\rho^{\Lambda^i}}^*\big)U_z^*\big)\big(s_\lambda^\lambda\phi\big(s_\nu^{\Lambda^i} {s_\tau^{\Lambda^i}}^*\big)\big)
&=\sum_{(\alpha,\beta)\in \Lambda^{\min}(\rho,\lambda\nu)}z^{d(\eta)-d(\rho)}s_{(\eta\alpha)(0,e_i)}^\Lambda \phi\big(s_{(\eta\alpha)(e_i,d(\eta\alpha))}^{\Lambda^i}{s_{\tau\beta}^{\Lambda^i}}^*\big)\\
&=\psi\big(\gamma^{\Lambda^i}_z\big(s_\eta^{\Lambda^i} {s_\rho^{\Lambda^i}}^*\big)\big)\big(s_\lambda^\lambda\phi\big(s_\nu^{\Lambda^i} {s_\tau^{\Lambda^i}}^*\big)\big).
\end{align*}
Since $X=\cspan\big\{s_\lambda^\Lambda\phi\big(s_\nu^{\Lambda^i} {s_\tau^{\Lambda^i}}^*\big):\lambda\in \Lambda^{e_i}, \nu, \tau\in \Lambda^i\big\}$, we see that $U_z\psi(a)U_z^*=\psi\big(\gamma^{\Lambda^i}_z(a)\big)$ for each $a\in C^*(\Lambda^i)$. 
\end{proof}
\end{lem}

\begin{prop}
\label{checking gauge invariant ideals}
Let $\Lambda$ be a locally-convex finitely aligned $k$-graph. Then $\ker(\psi)$, $\ker(\psi)^\perp$, and $\psi^{-1}(\mathcal{K}_{C^*(\Lambda^i)}(X))$ are gauge-invariant ideals of $C^*(\Lambda^i)$. Hence, the Katsura ideal $J_X:=\psi^{-1}(\mathcal{K}_{C^*(\Lambda^i)}(X))\cap \ker(\psi)^\perp$ is a gauge-invariant ideal of $C^*(\Lambda^i)$.
\begin{proof}
We begin by showing that $\ker(\psi)$ is gauge-invariant. For $z\in \mathbb{T}^{k-1}$, let $U_z$ be the map described in Lemma~\ref{existence of U_z}. For $a\in \ker(\psi)$, we have
$
\psi\big(\gamma_z^{\Lambda^i}(a)\big)=U_z\psi(a)U_z^*=0. 
$
Hence, $\gamma_z^{\Lambda^i}(a)\in \ker(\psi)$. From this we also see that $\ker(\psi)^\perp$ is gauge-invariant: if $a\in \ker(\psi)^\perp$, $b\in \ker(\psi)$, and $z\in \mathbb{T}^{k-1}$, then
$
\gamma_z^{\Lambda^i}(a)b=\gamma_z^{\Lambda^i}\big(a\gamma_{\overline{z}}^{\Lambda^i}(b)\big)=\gamma_z^{\Lambda^i}(0)=0.
$
It remains to show that $\psi^{-1}(\mathcal{K}_{C^*(\Lambda^i)}(X))$ is gauge-invariant. This follows from the fact that $\mathcal{K}_{C^*(\Lambda^i)}(X)$ is an ideal of $\mathcal{L}_{C^*(\Lambda^i)}(X)$, $U_z\in\mathcal{L}_{C^*(\Lambda^i)}(X)$ for each $z\in \mathbb{T}^{k-1}$, and $\psi\big(\gamma^{\Lambda^i}_z(a)\big)=U_z\psi(a)U_z^*$ for each $a\in C^*(\Lambda^i)$.
\end{proof}
\end{prop}

Now that we know the Katsura ideal $J_X$ is gauge-invariant, we seek to apply the analysis of \cite[Section~4]{MR3262073} to determine its generators. Note that when $k=2$ we could also apply the somewhat simpler analysis of \cite{MR1988256}, which deals just with directed graphs. Loosely speaking, if $\Sigma$ is a finitely-aligned $k$-graph, then a gauge-invariant ideal of $C^*(\Sigma)$ is generated (as an ideal) by its vertex projections and a collection of projections corresponding to certain finite exhaustive subsets of a subgraph of $\Sigma$. We now summarise the parts of \cite{MR3262073} that we will need.  

Suppose $I$ is a gauge-invariant ideal of $C^*(\Sigma)$. By \cite[Lemma~4.3]{MR3262073} the set 
\[
H_I:=\{v\in \Sigma^0:s_v^\Sigma\in I\}
\]
is saturated and hereditary (in the sense of  \cite[Definition~4.1]{MR3262073}). In particular, 
\[
\Sigma\setminus \Sigma H_I:=\{\lambda\in \Sigma:s(\lambda)\not\in H_I\}
\]
is itself a finitely aligned $k$-graph. If we define
\[
\tilde{B}_{I}:=\Big\{E\in \mathrm{FE}(\Sigma\setminus \Sigma H_I):E\subseteq \bigcup_{j=1}^k \Sigma^{e_j}, \Delta(s^\Sigma)^E\in I
\Big\},
\]
then \cite[Theorem~4.6]{MR3262073} (along with Proposition~\ref{sufficient to consider exhaustive sets of edges}) tells us that $I$ is generated as an ideal of $C^*(\Sigma)$ by the collection of projections
\[
\big\{s_v^{\Sigma}:v\in H_I\big\}\cup \big\{\Delta(s^\Sigma)^E:E\in \tilde{B}_I\big\}.
\]

We now determine the generators of the gauge-invariant ideal $\mathrm{ker}(\psi)$. We will see that due to the local-convexity of $\Lambda$, the ideal $\mathrm{ker}(\psi)$ is generated (as an ideal of $C^*(\Lambda^i)$) precisely by those vertex projections that act trivially on $X$, and that these projections correspond to the vertices that do not admit an edge of degree $e_i$. 

\begin{prop}
\label{generators of ker(psi)}
Let $\Lambda$ be a locally-convex finitely aligned $k$-graph. Consider the gauge-invariant ideal $\mathrm{ker}(\psi)$ of  $C^*(\Lambda^i)$. Then
\begin{enumerate}[label=\upshape(\roman*)]
\item $H_{\mathrm{ker}(\psi)}=\{v\in \Lambda^0:v\Lambda^{e_i}=\emptyset\}$;
\item $\tilde{B}_{\mathrm{ker}(\psi)}=\big\{E\in\mathrm{FE}(\Lambda^i\setminus \Lambda^i H_{\mathrm{ker}(\psi)}):E\subseteq \bigcup_{j=1}^k \Lambda^{e_j}\big\}$ with
\[
\Delta(s^{\Lambda^i})^E=0
\quad 
\text{for any $E\in\mathrm{FE}(\Lambda^i\setminus \Lambda^i H_{\mathrm{ker}(\psi)})$.}
\] 
\end{enumerate}
Hence, $\mathrm{ker}(\psi)$ is generated as an ideal of $C^*(\Lambda^i)$ by the collection of vertex projections $\{s_v^{\Lambda^i}:v\Lambda^{e_i}=\emptyset\}$. In particular, if $v\Lambda^{e_i}$ is nonempty for each $v\in \Lambda^0$, then the left action of $C^*(\Lambda^i)$ on $X$ is faithful.  
\begin{proof}
For any $v\in \Lambda^0$, $\lambda\in \Lambda^{e_i}$, and $a\in C^*(\Lambda^i)$, we see that
\[
\psi(s_v^{\Lambda^i})(s_\lambda^\Lambda\phi(a))
=\phi(s_v^{\Lambda^i})s_\lambda^\Lambda\phi(a)
=s_v^{\Lambda}s_\lambda^\Lambda\phi(a)
=\begin{cases}
s_\lambda^\Lambda\phi(a) & \text{if $r(\lambda)=v$}\\
0 & \text{otherwise.}
\end{cases}
\]
Since $X=\cspan\{s_\lambda^\Lambda\phi(a):\lambda\in \Lambda^{e_i}, a\in C^*(\Lambda^i)\}$, part (i) now follows immediately. 

We now prove part (ii). Suppose $E\in \mathrm{FE}(\Lambda^i\setminus \Lambda^iH_{\mathrm{ker}(\psi)})$. We claim that $r(E)\Lambda^{e_i}$ is nonempty. Since $E$ is nonempty, we can choose $\nu\in E$. Then $s(\nu)\not\in H_{\mathrm{ker}(\psi)}$, and so $s(\nu)\Lambda^{e_i}\neq \emptyset$ by part (i). By the factorisation property, $r(E)\Lambda^{e_i}=r(\nu)\Lambda^{e_i}\neq \emptyset$, which proves the claim. 
Next we show that $E\in \mathrm{FE}(\Lambda^i)$. Fix $\lambda\in r(E)\Lambda^i$. Since $\Lambda$ is locally-convex and $r(E)\Lambda^{e_i}\neq \emptyset$, we must have that $s(\lambda)\Lambda^{e_i}\neq \emptyset$. Thus, $s(\lambda)\not\in H_{\mathrm{ker}(\psi)}$, and we see that $\lambda\in r(E)\Lambda^i\setminus \Lambda^iH_{\mathrm{ker}(\psi)}$. Since $E$ is exhaustive in $\Lambda^i\setminus \Lambda^i H_{\mathrm{ker}(\psi)}$, we can find $\mu\in E$ such that $\mathrm{MCE}(\lambda,\mu)\neq \emptyset$. Hence, $E\in \mathrm{FE}(\Lambda^i)$ as claimed. Applying relation (CK) in $C^*(\Lambda^i)$ gives
$
\Delta(s^{\Lambda^i})^E=0,
$
which is certainly an element of $\mathrm{ker}(\psi)$. This completes the proof of part (ii). 
\end{proof}
\end{prop}

We will use the following product to sum transformation repeatedly, so we state it as a separate result.

\begin{lem}
\label{products to sums general}
Let $\Lambda$ be a finitely aligned $k$-graph and $\{r_\lambda:\lambda\in \Lambda\}$ a Toeplitz--Cuntz--Krieger $\Lambda$-family. Then for each $v\in \Lambda^0$ and each nonempty finite set $F\subseteq v\Lambda$ we have
\begin{equation}
\label{useful product to sum}
\begin{aligned}
\Delta(r)^F=r_v+\sum_{\substack{\emptyset \neq G\subseteq F\\ \lambda\in \mathrm{MCE}(G)}}
(-1)^{|G|}r_\lambda r_\lambda^*.
\end{aligned}
\end{equation}
In particular, if $n\in \mathbb{N}^k$ and $F$ is a nonempty finite subset of $v\Lambda^n$, then 
\begin{equation}
\label{products to sums equation}
\begin{aligned}
\Delta(r)^F=r_v-\sum_{\lambda\in F}r_\lambda r_\lambda^*.
\end{aligned}
\end{equation}
\begin{proof}
To prove \eqref{useful product to sum}, we will use induction on $|F|$. Clearly, when $|F|=1$ (say $F=\{\lambda\}$), both sides of \eqref{useful product to sum} are $r_v-r_\lambda r_\lambda^*$. So suppose that \eqref{useful product to sum} holds whenever $|F|=n$ and fix a set $F'\subseteq v\Lambda$ with $|F'|=n+1$. Then for any $\mu\in F'$, the inductive hypothesis gives
\begin{align*}
\Delta(r)^{F'}
&=(r_v-r_\mu r_\mu^*)\Delta(r)^{F'\setminus \{\mu\}}\\
&=(r_v-r_\mu r_\mu^*)\Bigg(r_v+\sum_{\substack{\emptyset \neq G\subseteq F'\setminus \{\mu\}\\ \lambda\in \mathrm{MCE}(G)}}
(-1)^{|G|}r_\lambda r_\lambda^*\Bigg)\\
&=r_v-r_\mu r_\mu^*+\sum_{\substack{\emptyset \neq G\subseteq F'\setminus \{\mu\}\\ \lambda\in \mathrm{MCE}(G)}}
(-1)^{|G|}r_\lambda r_\lambda^*-\sum_{\substack{\emptyset \neq G\subseteq F'\setminus \{\mu\}\\ \lambda\in \mathrm{MCE}(G)}}
(-1)^{|G|}r_\mu r_\mu^*r_\lambda r_\lambda^*.
\end{align*}
Applying relation (TCK3) to the product in the last sum, shows that this is equal to
\begin{align*}
r_v-r_\mu r_\mu^*+\sum_{\substack{\emptyset \neq G\subseteq F'\setminus \{\mu\}\\ \lambda\in \mathrm{MCE}(G)}}
(-1)^{|G|}r_\lambda r_\lambda^*
&+\sum_{\substack{\emptyset \neq G\subseteq F'\setminus \{\mu\}\\ \lambda\in \mathrm{MCE}(G\cup \{\mu\})}}
(-1)^{(|G|+1)}r_\lambda r_\lambda^*\\
&=r_v+\sum_{\substack{\emptyset \neq G\subseteq F'\\ \lambda\in \mathrm{MCE}(G)}}
(-1)^{|G|}r_\lambda r_\lambda^*.
\end{align*}
Thus, \eqref{useful product to sum} follows by induction. 

To see how \eqref{products to sums equation} follows from \eqref{useful product to sum}, observe that if $F\subseteq v\Lambda^n$, then for any $\emptyset \neq G\subseteq F$ we have $\mathrm{MCE}(G)=G$ if $G$ is a singleton set, whilst $\mathrm{MCE}(G)=\emptyset$ if $|G|\geq 2$. 
\end{proof}
\end{lem}

The next result tells us precisely which vertex projections belong to the gauge-invariant ideals $\mathrm{ker}(\psi)^\perp$ and $\psi^{-1}(\mathcal{K}_{C^*(\Lambda^i)}(X))$ (and so to $J_X$). 

\begin{prop}
\label{vertex projection in Katsura ideal}
Let $\Lambda$ be a locally-convex finitely aligned $k$-graph. Then for any vertex $v\in \Lambda^0$, $s_v^{\Lambda^i}\in \mathrm{ker}(\psi)^\perp$ if and only if $v\Lambda^{e_i}$ is nonempty, and $s_v^{\Lambda^i}\in \psi^{-1}(\mathcal{K}_{C^*(\Lambda^i)}(X))$ if and only if $v\Lambda^{e_i}$ is finite. Hence, $H_{J_X}=\{v\in \Lambda^0:0<|v\Lambda^{e_i}|<\infty\}$. In particular, when $v\Lambda^{e_i}$ is finite
\begin{equation}
\label{simple compact operators}
\psi(s_v^{\Lambda^i})=\sum_{\lambda\in v\Lambda^{e_i}}\Theta_{s_\lambda^\Lambda,s_\lambda^\Lambda}.
\end{equation}
If, in addition, $v\Lambda^{e_i}$ is nonempty, then 
\[
s_v^\Lambda=\sum_{\lambda\in v\Lambda^{e_i}}s_\lambda^\Lambda{s_\lambda^\Lambda}^*
\]
\begin{proof}
We begin by proving that $s_v^{\Lambda^i}\in \mathrm{ker}(\psi)^\perp$ if and only if $v\Lambda^{e_i}$ is nonempty. Fix $v\in \Lambda^0$. If $v\Lambda^{e_i}=\emptyset$, then $s_v^{\Lambda^i}\in \mathrm{ker}(\psi)$ by Proposition~\ref{generators of ker(psi)}. Thus, $s_v^{\Lambda^i}\not\in \mathrm{ker}(\psi)^\perp$ (otherwise we would have $0=(s_v^{\Lambda^i})^2=s_v^{\Lambda^i}$, which is clearly impossible). For the converse, suppose that $v\Lambda^{e_i}\neq\emptyset$. We need to show that $s_v^{\Lambda^i}a=0$ for each $a\in \mathrm{ker}(\psi)$. Since Proposition~\ref{generators of ker(psi)} tells us that $\mathrm{ker}(\psi)$ is generated as an ideal of $C^*(\Lambda^i)$ by the projections $\{s_w^{\Lambda^i}:w\Lambda^{e_i}=\emptyset\}$, it suffices to show that 
\[
s_v^{\Lambda^i}\big(s_\lambda^{\Lambda^i}{s_\nu^{\Lambda^i}}^*s_w^{\Lambda^i}\big)=0
\]
whenever $\lambda,\nu\in \Lambda^i$ and $w\in \Lambda^0$ is such that $w\Lambda^{e_i}=\emptyset$. Now $s_v^{\Lambda^i}\big(s_\lambda^{\Lambda^i}{s_\nu^{\Lambda^i}}^*s_w^{\Lambda^i}\big)$ is certainly zero if $r(\lambda)\neq v$, or $s(\lambda)\neq s(\nu)$, or $r(\nu)\neq w$, so we suppose otherwise. Since $v\Lambda^{e_i}\neq \emptyset$ and $\lambda\in \Lambda^i$, the local convexity of $\Lambda$ forces $s(\lambda)\Lambda^{e_i}\neq \emptyset$. Since $s(\lambda)=s(\nu)$, the factorisation property then implies that $r(\nu)\Lambda^{e_i}\neq \emptyset$. But this is impossible since $r(\nu)=w$. Thus, $s_v^{\Lambda^i}\in \mathrm{ker}(\psi)^\perp$. We conclude that $s_v^{\Lambda^i}\in \mathrm{ker}(\psi)^\perp$ if and only if $v\Lambda^{e_i}$ is nonempty.

Now we move on to proving that $s_v^{\Lambda^i}\in \psi^{-1}(\mathcal{K}_{C^*(\Lambda^i)}(X))$ if and only if $v\Lambda^{e_i}$ is finite. The proof uses the same ideas as \cite[Proposition~4.4]{MR1722197}. 

We claim that for any $v\in \Lambda^0$, the set $v\Lambda^{e_i}$ is exhaustive in $\Lambda$ provided it is nonempty. To see this, suppose that $\lambda\in v\Lambda$. We need to show there exists $\mu\in v\Lambda^{e_i}$ such that $\Lambda^{\min}(\lambda,\mu)\neq \emptyset$. If $\lambda=v$ then, for any $\mu\in v\Lambda^{e_i}$, we have $\{(\mu,s(\mu))\}=\Lambda^{\min}(\lambda,\mu)$. If $d(\lambda)_i\neq 0$, then with $\mu:=\lambda(0,e_i)\in v\Lambda^{e_i}$ we have $\{(s(\lambda),\lambda(e_i,d(\lambda)))\}=\Lambda^{\min}(\lambda,\mu)$. If $d(\lambda)_i=0$, then the  local-convexity of $\Lambda$ allows us to choose an edge $\nu\in s(\lambda)\Lambda^{e_i}$. With $\mu:=(\lambda\nu)(0,e_i)\in v\Lambda^{e_i}$, we have $(\nu,(\lambda\nu)(e_i,d(\lambda\nu)))\in \Lambda^{\min}(\lambda,\mu)$. Thus, $v\Lambda^{e_i}$ is exhaustive in $\Lambda$. Thus, if $v\Lambda^{e_i}$ is finite and nonempty, then relation (CK) tells us that $\Delta(s^\Lambda)^{v\Lambda^{e_i}}=0$. Applying Lemma~\ref{products to sums general} with $F=v\Lambda^{e_i}$, we conclude that $s_v^\Lambda=\sum_{\lambda\in v\Lambda^{e_i}}s_\lambda^\Lambda {s_\lambda^\Lambda}^*$.

We note that for each $\lambda\in v\Lambda^{e_i}$, $s_\lambda^\Lambda\in X$. To show that $s_v^{\Lambda^i}\in \psi^{-1}(\mathcal{K}_{C^*(\Lambda^i)}(X))$ when $v\Lambda^{e_i}$ is finite, it suffices to show that \eqref{simple compact operators} holds.  If $v\Lambda^{e_i}=\emptyset$, then the right-hand side of \eqref{simple compact operators} is the empty sum, and so zero, whilst $\psi(s_v^{\Lambda^i})=0$ by Proposition~\ref{generators of ker(psi)}. On the other hand, if $v\Lambda^{e_i}$ is nonempty and finite, then for any $\mu\in \Lambda^{e_i}$ and $a\in C^*(\Lambda^i)$, we have
\begin{align*}
\psi(s_v^{\Lambda^i})(s_\mu^\Lambda\phi(a))
=s_v^\Lambda s_\mu^\Lambda\phi(a)
=\sum_{\lambda\in v\Lambda^{e_i}}s_\lambda^\Lambda {s_\lambda^\Lambda}^*s_\mu^\Lambda\phi(a)
&=\sum_{\lambda\in v\Lambda^{e_i}}s_\lambda^\Lambda\cdot \langle s_\lambda^\Lambda, s_\mu^\Lambda\phi(a)\rangle_{C^*(\Lambda^i)}\\
&=\sum_{\lambda\in v\Lambda^{e_i}}\Theta_{s_\lambda^\Lambda,s_\lambda^\Lambda}(s_\mu^\Lambda\phi(a)).
\end{align*}
Since $X=\cspan\{s_\mu^\Lambda\phi(a):\mu\in \Lambda^{e_i}, a\in C^*(\Lambda^i)\}$, we conclude that \eqref{simple compact operators} holds whenever $v\in \Lambda^0$ and $v\Lambda^{e_i}$ is finite.

It remains to show that if $v\in \Lambda^0$ and $s_v^{\Lambda^i}\in \psi^{-1}(\mathcal{K}_{C^*(\Lambda^i)}(X))$, then $v\Lambda^{e_i}$ is finite. Looking for a contradiction, suppose that $v\Lambda^{e_i}$ is infinite and $\psi(s_v^{\Lambda^i})$ is compact. Since $X=\cspan\{s_\lambda^\Lambda \phi(a):\lambda\in \Lambda^{e_i}, a\in C^*(\Lambda^i)\}$, there exist finite sets $E,F\subseteq \Lambda^{e_i}$, $G,H\subseteq C^*(\Lambda^i)$ such that
\[
\|K-\psi(s_v^{\Lambda^i})\|_{\mathcal{L}_{C^*(\Lambda^i)}(X)}<1
\]
where 
\[
K:=\sum_{\substack{(\lambda,a)\in E\times G, \\ (\mu,b)\in F\times H}}\Theta_{s_\lambda^\Lambda \phi(a), s_\mu^\Lambda \phi(b)}
\in\mathcal{K}_{C^*(\Lambda^i)}(X).
\]
Since $F$ is finite and $v\Lambda^{e_i}$ is infinite, we can choose $\nu\in v\Lambda^{e_i}\setminus F$. Then $s_\nu^\Lambda\in X$ and 
\[
\psi(s_v^{\Lambda^i})(s_\nu^\Lambda)=\phi(s_v^{\Lambda^i})s_\nu^\Lambda=s_v^\Lambda s_\nu^\Lambda=s_\nu^\Lambda.
\]
However, since $d(\mu)=d(\nu)$ for each $\mu\in F$, we have that $\mathrm{MCE}(\mu,\nu)=\emptyset$ for every $\mu\in F$, and so
\[
K(s_\nu^\Lambda)
=\sum_{\substack{(\lambda,a)\in E\times G, \\ (\mu,b)\in F\times H}}\Theta_{s_\lambda^\Lambda \phi(a), s_\mu^\Lambda \phi(b)}(s_\nu^\Lambda)
=\sum_{\substack{(\lambda,a)\in E\times G, \\ (\mu,b)\in F\times H}} s_\lambda^\Lambda \phi(ab^*){s_\mu^\Lambda}^*s_\nu^\Lambda
=0.
\]
Since the norm on $X$ is the restriction of the norm on $C^*(\Lambda)$, we have that 
\[
\|s_\nu^\Lambda\|_X^2=\|s_\nu^\Lambda\|_{C^*(\Lambda)}^2=\|{s_\nu^\Lambda}^*s_\nu^\Lambda\|_{C^*(\Lambda)}=\|s_{s(\nu)}^\Lambda\|_{C^*(\Lambda)}=1.
\] 
Thus, 
\begin{align*}
\|K-\psi(s_v^{\Lambda^i})\|_{\mathcal{L}_{C^*(\Lambda^i)}(X)}
&=\sup_{\|x\|_X\leq 1}\|(K-\psi(s_v^{\Lambda^i}))(x)\|_X\\
&\geq \|(K-\psi(s_v^{\Lambda^i}))(s_\nu^\Lambda)\|_X
=\|s_\nu^\Lambda\|_X
=1,
\end{align*}
which is a contradiction.
\end{proof}
\end{prop}

Now that we have determined $H_{J_X}$, we move on to $\tilde{B}_{J_X}$. In contrast to $\mathrm{ker}(\psi)$, where the set $\tilde{B}_{\mathrm{ker}(\psi)}$ did not contribute any nonzero generators, the ideal $\psi^{-1}(\mathcal{K}_{C^*(\Lambda^i)}(X))$ is not necessarily generated solely by its vertex projections (i.e. the projections corresponding to vertices that admit finitely many edges of degree $e_i$). The purpose of the next lemma is to determine for which finite exhaustive sets $E\in \mathrm{FE}(\Lambda^i\setminus \Lambda^i H_{J_X})$ with $E\subseteq \bigcup_{j=1}^k \Lambda^{e_j}$ does the projection $\Delta(s^{\Lambda^i})^E$ belong to the Katsura ideal. We were somewhat surprised to discover that this occurs if and only if $E$ can be extended to an exhaustive subset of $\Lambda$ by adding in a finite collection of edges of degree $e_i$.  

\begin{lem}
\label{extending to finite exhaustive set implies compact}
Let $\Lambda$ be a locally-convex finitely aligned $k$-graph. Let $E\in \mathrm{FE}(\Lambda^i\setminus \Lambda^i H_{J_X})$ with $E\subseteq \bigcup_{j=1}^k\Lambda^{e_j}$. Then 
$
\Delta(s^{\Lambda^i})^E\in J_X
$
if and only if there exists a finite set $F\subseteq r(E)\Lambda^{e_i}$ such that
$E\cup F\in \mathrm{FE}(\Lambda)$. In particular, if $E\cup F\in \mathrm{FE}(\Lambda)$ for some $F\subseteq r(E)\Lambda^{e_i}$, then
\[
\psi\big(\Delta(s^{\Lambda^i})^E\big)
=
\sum_{\substack{G\subseteq E\cup F\\G\cap F\neq \emptyset\\\mu\in \mathrm{MCE}(G)}}
(-1)^{(|G|+1)}\Theta_{s_\mu^\Lambda,s_\mu^\Lambda}\in \mathcal{K}_{C^*(\Lambda^i)}(X). 
\]
\begin{proof}
Fix $E\in \mathrm{FE}(\Lambda^i\setminus \Lambda^i H_{J_X})$ with $E\subseteq \bigcup_{j=1}^k\Lambda^{e_j}$. We show that $\Delta(s^{\Lambda^i})^E\in \mathrm{ker}(\psi)^\perp$. Consider the situation where $r(E)\Lambda^{e_i}=\emptyset$. We claim that $r(E)\Lambda^i=r(E)\Lambda^i\setminus \Lambda^i H_{J_X}$. Clearly, $r(E)\Lambda^i\setminus \Lambda^i H_{J_X}\subseteq r(E)\Lambda^i$, and we just need to check the reverse set inclusion. If $\lambda\in r(E)\Lambda^i$, then $s(\lambda)\Lambda^{e_i}=\emptyset$ by the factorisation property, and so $s(\lambda)\not \in H_{J_X}$ by Proposition~\ref{vertex projection in Katsura ideal}. Hence, $\lambda\in r(E)\Lambda^i\setminus \Lambda^i H_{J_X}$, which proves the claim. Consequently, $E\in r(E)\mathrm{FE}(\Lambda^i)$, and relation (CK) in $C^*(\Lambda^i)$ says that $\Delta(s^{\Lambda^i})^E=0$, which is certainly in $\mathrm{ker}(\psi)^\perp$. On the other hand, if $r(E)\Lambda^{e_i}\neq\emptyset$, then $s_{r(E)}^{\Lambda^i}\in\mathrm{ker}(\psi)^\perp$ by Proposition~\ref{vertex projection in Katsura ideal}. Since $\mathrm{ker}(\psi)^\perp$ is an ideal of $C^*(\Lambda^i)$, we see that $\Delta(s^{\Lambda^i})^E=s_{r(E)}^{\Lambda^i}\Delta(s^{\Lambda^i})^E\in \mathrm{ker}(\psi)^\perp$.

Now suppose that there exists a set $F\subseteq r(E)\Lambda^{e_i}$ such that
$E\cup F\in \mathrm{FE}(\Lambda)$. Applying the Cuntz--Krieger relation in $C^*(\Lambda)$, Lemma~\ref{products to sums general} gives
\[
0=\Delta(s^\Lambda)^{E\cup F}=s_v^\Lambda+\sum_{\substack{\emptyset \neq G\subseteq E\cup F\\ \mu\in \mathrm{MCE}(G)}}(-1)^{|G|}s_\mu^{\Lambda}{s_\mu^{\Lambda}}^*. 
\]
Splitting this sum, we get that
\begin{equation}
\label{splitting sum}
s_v^\Lambda+\sum_{\substack{\emptyset \neq G\subseteq E\\ \mu\in \mathrm{MCE}(G)}}(-1)^{|G|}s_\mu^{\Lambda}{s_\mu^{\Lambda}}^*
=
\sum_{\substack{G\subseteq E\cup F\\ G\cap F\neq \emptyset\\\mu\in \mathrm{MCE}(G)}}(-1)^{(|G|+1)}s_\mu^{\Lambda}{s_\mu^{\Lambda}}^*. 
\end{equation}
Since $E\subseteq \Lambda^i$, we can use Lemma~\ref{products to sums general} again to see that 
\begin{equation}
\label{applying lemma 5.12 to E}
\begin{aligned}
s_v^\Lambda+\sum_{\substack{\emptyset \neq G\subseteq E\\ \mu\in \mathrm{MCE}(G)}}(-1)^{|G|}s_\mu^{\Lambda}{s_\mu^{\Lambda}}^*
=
\phi\bigg(
s_v^{\Lambda^i}+\sum_{\substack{\emptyset \neq G\subseteq E\\ \mu\in \mathrm{MCE}(G)}}(-1)^{|G|}s_\mu^{\Lambda^i}{s_\mu^{\Lambda^i}}^*
\bigg)
=
\phi\big(\Delta(s^{\Lambda^i})^E\big).
\end{aligned}
\end{equation}
Next, observe that if $G\subseteq E\cup F$ and $G\cap F\neq \emptyset$, then $\max\{d(\nu)_i:\nu\in G\}=1$. Hence, if $\mu \in \mathrm{MCE}(G)$, then $d(\mu)_i=1$, and so $s_\mu^\Lambda\in X$.  Since $\psi\big(\Delta(s^{\Lambda^i})^E\big)$ is multiplication by $\phi\big(\Delta(s^{\Lambda^i})^E\big)$ on $X\subseteq C^*(\Lambda)$, and $\Theta_{s_\mu^\Lambda,s_\mu^\Lambda}$ is multiplication by $s_\mu^{\Lambda}{s_\mu^{\Lambda}}^*$ for each $s_\mu^\Lambda\in X$, \eqref{splitting sum} and \eqref{applying lemma 5.12 to E} imply that
\[
\psi\big(\Delta(s^{\Lambda^i})^E\big)
=
\sum_{\substack{G\subseteq E\cup F\\G\cap F\neq \emptyset\\\mu\in \mathrm{MCE}(G)}}
(-1)^{(|G|+1)}\Theta_{s_\mu^\Lambda,s_\mu^\Lambda}. 
\]
Thus, $\Delta(s^{\Lambda^i})^E\in \psi^{-1}(\mathcal{K}_{C^*(\Lambda^i)}(X))$, and we conclude that $\Delta(s^{\Lambda^i})^E\in J_X$. 

Conversely, suppose that $\Delta(s^{\Lambda^i})^E\in J_X$. Since $\Delta(s^{\Lambda^i})^E$ acts compactly on $X$, we can find finite sets $G,H\subseteq\{\lambda\in \Lambda: d(\lambda)_i=1\}$ such that
\begin{equation}
\label{compact op norm inequality}
\big\|K-\psi\big(\Delta(s^{\Lambda^i})^E\big)\big\|_{\mathcal{L}_{C^*(\Lambda^i)}(X)}<1,
\end{equation}
where 
\[
K:=\sum_{(\mu,\nu)\in G\times H}c_{(\mu,\nu)}\Theta_{s_\mu^\Lambda , s_\nu^\Lambda}\in \mathcal{K}_{C^*(\Lambda^i)}(X).
\]
and the $c_{(\mu,\nu)}$ are constants. 

We define $F:=\{\lambda(0,e_i):\lambda\in r(E)H\}$ and claim that $E\cup F\in \mathrm{FE}(\Lambda)$. Since $E$ and $H$ are finite, so is $E\cup F$. Furthermore, since $r(E)\not\in E$, and $F$ consists of edges, we see that $r(E)\not \in E\cup F$. Thus, it remains to show that $E\cup F$ is exhaustive in $\Lambda$. Looking for a contradiction, suppose that there exists some $\tau\in r(E)\Lambda$ that does not have a minimal common extension with anything in $E\cup F$. We consider the situations where $d(\tau)_i=1$, $d(\tau)_i=0$, and $d(\tau)_i\geq 2$ separately. 

Firstly, we consider the case where $d(\tau)_i=1$ (note: this implies that $s_\tau^\Lambda\in X$). Clearly, if $\nu\in H$ and $r(\nu)\neq r(E)$, then $\mathrm{MCE}(\nu,\tau)=\emptyset$. Since $\mathrm{MCE}(\tau,\nu(0,e_i))$ is, by assumption, empty for each $\nu\in r(E)H$, the factorisation property implies that $\mathrm{MCE}(\tau,\nu)=\emptyset$ for each $\nu\in H$. Hence,
\begin{align*}
K(s_\tau^\Lambda)
=\sum_{(\mu,\nu)\in G\times H}c_{(\mu,\nu)}\Theta_{s_\mu^\Lambda , s_\nu^\Lambda}(s_\tau^\Lambda)
&=\sum_{(\mu,\nu)\in G\times H}c_{(\mu,\nu)}s_\mu^\Lambda {s_\nu^\Lambda}^*s_\tau^\Lambda\\
&=\sum_{\substack{(\mu,\nu)\in G\times H\\ (\alpha,\beta)\in \Lambda^{\min}(\nu,\tau)}}c_{(\mu,\nu)}s_{\mu\alpha}^\Lambda {s_\beta^\Lambda}^*
=0.
\end{align*}
Similarly, since $\mathrm{MCE}(\tau,\lambda)=\emptyset$ for each $\lambda\in E$, we have that 
\[
(s_{r(E)}^{\Lambda}-s_\lambda^{\Lambda}{s_\lambda^{\Lambda}}^*)s_\tau^\Lambda=s_\tau^\Lambda-\sum_{(\alpha,\beta)\in \Lambda^{\min}(\lambda,\tau)}s_{\lambda\alpha}^\Lambda {s_\beta^\Lambda}^*=s_\tau^\Lambda
\]
for each $\lambda\in E$. Thus,
\[
\psi\big(\Delta(s^{\Lambda^i})^E\big)(s_\tau^\Lambda)
=\Big(\prod_{\lambda\in E}\big(s_{r(E)}^{\Lambda}-s_\lambda^{\Lambda}{s_\lambda^{\Lambda}}^*\big)\Big)s_\tau^\Lambda
=s_\tau^\Lambda,
\] 
and so
\begin{equation}
\label{required contradiction}
\big(K-\psi\big(\Delta(s^{\Lambda^i})^E\big)\big)(s_\tau^\Lambda)=-s_\tau^\Lambda. 
\end{equation}
Since
\[
\|s_\tau^\Lambda\|_X=\|s_\tau^\Lambda\|_{C^*(\Lambda)}=\|{s_\tau^\Lambda}^*s_\tau^\Lambda\|_{C^*(\Lambda)}^{1/2}=
\|s_{s(\tau)}^\Lambda\|_{C^*(\Lambda)}^{1/2}=1\neq 0,
\]
\eqref{required contradiction} contradicts \eqref{compact op norm inequality}. Hence, for each $\tau\in r(E)\Lambda$ with $d(\tau)_i=1$, there must exist $\lambda\in E\cup F$ such that $\mathrm{MCE}(\tau,\lambda)\neq\emptyset$. 

Now consider the situation when $d(\tau)_i=0$. Consider the case where $r(E)\Lambda^{e_i}$ is nonempty. The local-convexity of $\Lambda$ allows us to choose $\xi\in s(\tau)\Lambda^{e_i}$, and so by the argument in the previous paragraph we can find $\lambda\in E\cup F$ such that $\mathrm{MCE}(\tau\xi,\lambda)\neq\emptyset$. The factorisation property then implies that $\mathrm{MCE}(\tau,\lambda)\neq\emptyset$. On the other hand, if $r(E)\Lambda^{e_i}$ is empty, then $s_{r(E)}^{\Lambda^i}\in \mathrm{ker}(\psi)$, and so
\[
\Delta(s^{\Lambda^i})^E
=\Delta(s^{\Lambda^i})^E s_{r(E)}^{\Lambda^i}
=0
\]
because $\Delta(s^{\Lambda^i})^E\in J_X\subseteq \mathrm{ker}(\psi)^\perp$. Thus,
\[
\Big(\prod_{\lambda\in E}\big(s_{r(E)}^{\Lambda}-s_\lambda^{\Lambda}{s_\lambda^{\Lambda}}^*\big)\Big)s_\tau^\Lambda=\Delta(s^{\Lambda^i})^Es_\tau^\Lambda=0,
\]
which is impossible if $\mathrm{MCE}(\tau,\lambda)=\emptyset$ for each $\lambda\in E$. Hence, $\mathrm{MCE}(\tau,\lambda)$ is nonempty for some $\lambda\in E$.

It remains to consider the case where $d(\tau)_i\geq 2$. Let $K_\tau:=\{j:r(E)\Lambda^{e_j}\neq \emptyset, d(\tau)_j=0\}$ (this could be the empty set). Since $\Lambda$ is locally-convex, we may choose $\xi\in s(\tau)\Lambda^{\sum_{j\in K_\tau}e_j}$. Define $\tau':=(\tau\xi)(0,d(\tau\xi)-(d(\tau\xi)_i-1)e_i)$. Observe that $d(\tau')_i=1$ and $d(\tau')_j=d(\tau\xi)_j$ for all $j\neq i$. Hence, we can find $\lambda\in E\cup F$ such that $\mathrm{MCE}(\tau',\lambda)\neq\emptyset$. If $\lambda\in F$ (which is a subset of $\Lambda^{e_i}$), then $\tau'$ must extend $\lambda$ because  $d(\tau')\geq e_i$. Hence, $\tau\xi$ also extends $\lambda$, and we have that $\{\tau\xi\}=\mathrm{MCE}(\tau\xi,\lambda)$. By the factorisation property, it follows that $\mathrm{MCE}(\tau,\lambda)\neq \emptyset$. Alternatively, $\lambda\in E$, and so $\lambda\in \Lambda^{e_j}$ for some $j\neq i$. Since $d(\tau')_j=d(\tau\xi)_j\geq 1$ by our choice of $\xi$, we see that $\tau'$ (and so $\tau\xi$) must extend $\lambda$. Thus, $\{\tau\xi\}=\mathrm{MCE}(\tau\xi,\lambda)$. Hence, by the factorisation property $\mathrm{MCE}(\tau,\lambda)\neq \emptyset$. This completes the proof of the claim that $E\cup F$ is exhaustive in $\Lambda$. 
\end{proof}
\end{lem}

We now have enough information to give a complete description of the Katsura ideal. 

\begin{prop}
\label{putting everything together}
Let $\Lambda$ be a locally-convex finitely aligned $k$-graph. Then 
\[
\Lambda^i\setminus \Lambda^i H_{J_X}=\{\lambda\in \Lambda^i:|s(\lambda)\Lambda^{e_i}|\in \{0,\infty\}\}
\]
and
\[
\tilde{B}_{J_X} = 
\Big\{
E\in \mathrm{FE}(\Lambda^i\setminus \Lambda^i H_{J_X}):
E\subseteq \bigcup_{j=1}^k\Lambda^{e_j}, E\cup F\in \mathrm{FE}(\Lambda) 
\text{ for some } F\subseteq r(E)\Lambda^{e_i}
\Big\}.
\]
Furthermore, the Katsura ideal $J_X$ is generated as an ideal of $C^*(\Lambda^i)$ by the collection of projections
\begin{equation}
\label{description of generators in Katsura ideal}
\big\{s_v^{\Lambda^i}:0<|v\Lambda^{e_i}|<\infty\big\}\cup \big\{\Delta(s^{\Lambda^i})^E:E\in \tilde{B}_{J_X}, r(E)\Lambda^{e_i}\neq \emptyset\big\}.
\end{equation}
\begin{proof}
We proved in Lemma~\ref{vertex projection in Katsura ideal} that $H_{J_X}=\{v\in \Lambda^0:0<|v\Lambda^{e_i}|<\infty\}$, from which the description of $\Lambda^i\setminus \Lambda^i H_{J_X}$ follows. Lemma~\ref{extending to finite exhaustive set implies compact} shows that if $E\in \mathrm{FE}(\Lambda^i\setminus \Lambda^i H_{J_X})$ and $E\subseteq \bigcup_{j=1}^k\Lambda^{e_j}$, then $\Delta(s^{\Lambda^i})^E\in J_X$ if and only if there exists $F\subseteq r(E)\Lambda^{e_i}$ such that $E\cup F\in \mathrm{FE}(\Lambda)$. This gives us the description of $\tilde{B}_{J_X}$. The first paragraph of the proof of Lemma~\ref{extending to finite exhaustive set implies compact} shows that if $E\in \mathrm{FE}(\Lambda^i\setminus \Lambda^i H_{J_X})$ and $r(E)\Lambda^{e_i}=\emptyset$, then $E\in \mathrm{FE}(\Lambda^i)$, and so $\Delta(s^{\Lambda^i})^E=0$. Consequently, \cite[Theorem~4.6]{MR3262073} and Proposition~\ref{sufficient to consider exhaustive sets of edges} tell us that the collection of projections in \eqref{description of generators in Katsura ideal} generate the gauge-invariant ideal $J_X$. 
\end{proof}
\end{prop}

The next example illustrates the subtlety addressed by Proposition~\ref{putting everything together} for graphs with infinite receivers, and served as the main motivation for the formulation of Lemma~\ref{extending to finite exhaustive set implies compact}. We thank Aidan Sims for bringing this example to our attention. 

\begin{exam}
Consider the locally-convex finitely aligned $2$-graph $\Lambda$ described in \cite[Example~A.3]{MR2069786} with $1$-skeleton 

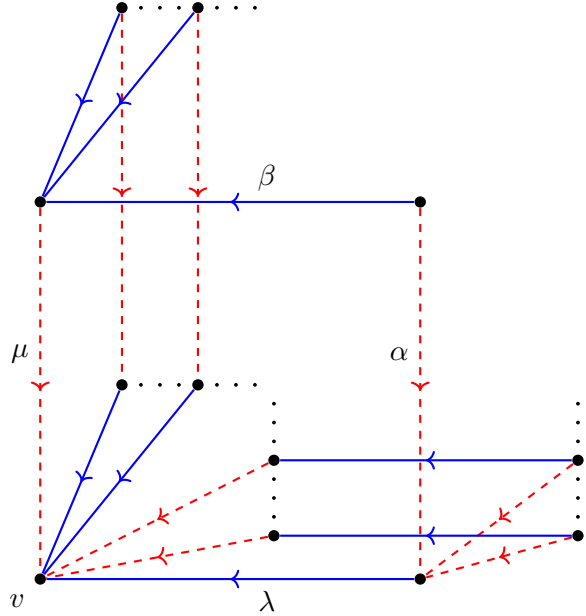
\tikzset{->-/.style={decoration={
  markings,
  mark=at position #1 with {\arrow{>}}},postaction={decorate}}}

\[
\begin{array}{lcr}
&\begin{tikzpicture}[scale=5]
    \node[label=below left:\small$v$,circle,inner sep=1.5pt,fill=black] (001) at (0,0,0) {};
    \node[circle,inner sep=1.5pt,fill=black] (011) at (0,1,0) {};
    \node[circle,inner sep=1.5pt,fill=black] (101) at (1,0,0) {};
    \node[circle,inner sep=1.5pt,fill=black] (111) at (1,1,0) {};
\node[circle,inner sep=1.5pt,fill=black] (a) at (0.6,0.9,1) {};
\node[circle,inner sep=0.5pt,fill=black] (x) at (0.65,0.9,1) {};
\node[circle,inner sep=0.5pt,fill=black] (y) at (0.7,0.9,1) {};
\node[circle,inner sep=0.5pt,fill=black] (z) at (0.75,0.9,1) {};
\node[circle,inner sep=1.5pt,fill=black] (b) at (0.8,0.9,1) {};
\node[circle,inner sep=0.5pt,fill=black] (x) at (0.85,0.9,1) {};
\node[circle,inner sep=0.5pt,fill=black] (y) at (0.9,0.9,1) {};
\node[circle,inner sep=0.5pt,fill=black] (z) at (0.95,0.9,1) {};
\node[circle,inner sep=1.5pt,fill=black] (c) at (0.6,1.9,1) {};
\node[circle,inner sep=0.5pt,fill=black] (x) at (0.65,1.9,1) {};
\node[circle,inner sep=0.5pt,fill=black] (y) at (0.7,1.9,1) {};
\node[circle,inner sep=0.5pt,fill=black] (z) at (0.75,1.9,1) {};
\node[circle,inner sep=1.5pt,fill=black] (d) at (0.8,1.9,1) {};
\node[circle,inner sep=0.5pt,fill=black] (x) at (0.85,1.9,1) {};
\node[circle,inner sep=0.5pt,fill=black] (y) at (0.9,1.9,1) {};
\node[circle,inner sep=0.5pt,fill=black] (z) at (0.95,1.9,1) {};
\node[circle,inner sep=1.5pt,fill=black] (e) at (1.8,0.5,1) {};
\node[circle,inner sep=0.5pt,fill=black] (x) at (1.8,0.55,1) {};
\node[circle,inner sep=0.5pt,fill=black] (y) at (1.8,0.6,1) {};
\node[circle,inner sep=0.5pt,fill=black] (z) at (1.8,0.65,1) {};
\node[circle,inner sep=1.5pt,fill=black] (f) at (1.8,0.7,1) {};
\node[circle,inner sep=0.5pt,fill=black] (x) at (1.8,0.75,1) {};
\node[circle,inner sep=0.5pt,fill=black] (y) at (1.8,0.8,1) {};
\node[circle,inner sep=0.5pt,fill=black] (z) at (1.8,0.85,1) {};
\node[circle,inner sep=1.5pt,fill=black] (g) at (1,0.5,1) {};
\node[circle,inner sep=0.5pt,fill=black] (x) at (1,0.55,1) {};
\node[circle,inner sep=0.5pt,fill=black] (y) at (1,0.6,1) {};
\node[circle,inner sep=0.5pt,fill=black] (z) at (1,0.65,1) {};
\node[circle,inner sep=1.5pt,fill=black] (h) at (1,0.7,1) {};
\node[circle,inner sep=0.5pt,fill=black] (x) at (1,0.75,1) {};
\node[circle,inner sep=0.5pt,fill=black] (y) at (1,0.8,1) {};
\node[circle,inner sep=0.5pt,fill=black] (z) at (1,0.85,1) {};
    \draw[->-=.5, thick, blue] (101)--(001)
        node[below, pos=0.4, black] {\small$\lambda$};
    \draw[->-=.5, thick, blue] (111)--(011)
        node[above, pos=0.4, black] {\small$\beta$};
    \draw[->-=.5, thick, red, dashed] (011)--(001)
        node[left, pos=0.4, black] {\small$\mu$};
    \draw[->-=.5, thick, red, dashed] (111)--(101)
        node[left, pos=0.4, black] {\small$\alpha$};
\draw[->-=.5, thick, blue] (a)--(001)
        node[left, pos=0.4, inner sep=0.5pt, black] {};
\draw[->-=.5, thick, blue] (b)--(001)
        node[left, pos=0.4, inner sep=0.5pt, black] {};
\draw[->-=.5, thick, red, dashed] (c)--(a)
        node[left, pos=0.4, inner sep=0.5pt, black] {};
\draw[->-=.5, thick, blue] (c)--(011)
        node[left, pos=0.4, inner sep=0.5pt, black] {};
\draw[->-=.5, thick, red, dashed] (d)--(b)
        node[left, pos=0.4, inner sep=0.5pt, black] {};
\draw[->-=.5, thick, blue] (d)--(011)
        node[left, pos=0.4, inner sep=0.5pt, black] {};
\draw[->-=.5, thick, red, dashed] (g)--(001)
        node[left, pos=0.4, inner sep=0.5pt, black] {};
\draw[->-=.5, thick, red, dashed] (h)--(001)
        node[left, pos=0.4, inner sep=0.5pt, black] {};
\draw[->-=.5, thick, red, dashed] (e)--(101)
        node[left, pos=0.4, inner sep=0.5pt, black] {};
\draw[->-=.5, thick, red, dashed] (f)--(101)
        node[left, pos=0.4, inner sep=0.5pt, black] {};
\draw[->-=.5, thick, blue] (e)--(g)
        node[left, pos=0.4, inner sep=0.5pt, black] {};
\draw[->-=.5, thick, blue] (f)--(h)
        node[left, pos=0.4, inner sep=0.5pt, black] {};
\end{tikzpicture}
\end{array}
\]
where blue (solid) edges have degree $e_1$ and red (dashed) edges have degree $e_2$. 

Let $i=2$ (i.e. we are removing the red (dashed) edges from the graph). Proposition~\ref{vertex projection in Katsura ideal} tells us that
\[
H_{J_X}=s(v\Lambda^{e_1}\setminus \{\lambda\}).
\]
We now determine the finite exhaustive sets in the graphs $\Lambda^2$ and $\Lambda^2\setminus \Lambda^2 H_{J_X}$.
Since $s(\mu)\Lambda^{e_1}=s(\mu)(\Lambda^2\setminus \Lambda^2 H_{J_X})^{e_1}$ is infinite, we see that
$s(\mu)\mathrm{FE}(\Lambda^2)$ and $s(\mu)\mathrm{FE}(\Lambda^2\setminus \Lambda^2 H_{J_X})$ are empty. Since $v\Lambda^{e_1}$ is infinite, but $v(\Lambda^2\setminus \Lambda^2 H_{J_X})^{e_1}=\{\lambda\}$, we see that $v\mathrm{FE}(\Lambda^2)$ is empty and $v\mathrm{FE}(\Lambda^2\setminus \Lambda^2 H_{J_X})=\{\{\lambda\}\}$. Observe that $\{\lambda\}$ is contained in $\{\lambda,\mu\}\in \mathrm{FE}(\Lambda)$. For each  $\eta\in v\Lambda^{e_2}\setminus \{\mu\}$, we have that $s(\eta)\Lambda^{e_1}=s(\eta)(\Lambda^2\setminus \Lambda^2 H_{J_X})^{e_1}$ is a singleton. Hence, for each $\eta\in v\Lambda^{e_2}\setminus \{\mu\}$, $s(\eta)\mathrm{FE}(\Lambda^2)=s(\eta)\mathrm{FE}(\Lambda^2\setminus \Lambda^2 H_{J_X})=\{s(\eta)\Lambda^{e_1}\}$. Moreover, for each $\eta\in v\Lambda^{e_2}\setminus \{\mu\}$, the singleton set $\{s(\eta)\Lambda^{e_1}\}$ is exhaustive in $\Lambda$. Hence, we conclude that
\[
\mathrm{FE}(\Lambda^2)=\{\{s(\eta)\Lambda^{e_1}\}:\eta\in v\Lambda^{e_2}\setminus \{\mu\}\}
\]
and
\[
\mathrm{FE}(\Lambda^2\setminus \Lambda^2 H_{J_X})=\mathrm{FE}(\Lambda^2)\cup \{\{\lambda\}\}=\tilde{B}_{J_X}.
\]
Since $s(\eta)\Lambda^{e_2}=\emptyset$ for each $\eta\in v\Lambda^{e_2}\setminus \{\mu\}$, Proposition~\ref{putting everything together} tells us that $J_X$ is generated as an ideal of $C^*(\Lambda^2)$ by the projections
\[
\big\{s_w^{\Lambda^2}:w\in s(v\Lambda^{e_1}\setminus \{\lambda\})\big\}\cup \big\{s_v^{\Lambda^2}-s_{\lambda}^{\Lambda^2}{s_{\lambda}^{\Lambda^2}}^*\big\}.
\]
Furthermore, for any $w\in s(v\Lambda^{e_1}\setminus \{\lambda\})$, Proposition~\ref{vertex projection in Katsura ideal} tells us that
\[
\psi\big(s_w^{\Lambda^2}\big)=\sum_{\tau\in w\Lambda^{e_2}}\Theta_{s_\tau^{\Lambda},s_\tau^{\Lambda}},
\]
which is a rank-one operator. By Lemma~\ref{extending to finite exhaustive set implies compact}, we also get that
\[
\psi\big(s_{v}^{\Lambda^2}-s_{\lambda}^{\Lambda^2}{s_{\lambda}^{\Lambda^2}}^*\big)
=\Theta_{s_{\mu}^{\Lambda}, s_{\mu}^{\Lambda}}-\Theta_{s_{\mu\beta}^{\Lambda}, s_{\mu\beta}^{\Lambda}}.
\]
\end{exam}

We now use Proposition~\ref{putting everything together} to prove our main theorem: when $\Lambda$ is locally-convex, $\mathcal{O}_X$ and $C^*(\Lambda)$ are isomorphic.

\begin{thm}
\label{main theorem}
Let $\Lambda$ be a locally-convex finitely aligned $k$-graph. If we let $\iota:X\rightarrow C^*(\Lambda)$ denote the inclusion map, then $(\iota, \phi)$ is a Cuntz--Pimsner covariant Toeplitz representation of $X$ in $C^*(\Lambda)$. For each $\lambda\in \Lambda$, define $u_\lambda \in \mathcal{O}_X$ by
$
u_\lambda:=j_X^{\otimes d(\lambda)_i}\left(\Omega_{d(\lambda)_i}\left(s_\lambda^\Lambda\right)\right)
$.
Then $\{u_\lambda:\lambda \in \Lambda\}$ is a Cuntz--Krieger $\Lambda$-family in $\mathcal{O}_X$. Furthermore, the $*$-homomorphisms $\iota \times_\mathcal{O} \phi:\mathcal{O}_X\rightarrow C^*(\Lambda)$ and $\pi_u:C^*(\Lambda)\rightarrow \mathcal{O}_X$ induced by the universal properties of $\mathcal{O}_X$ and $C^*(\Lambda)$ are mutually inverse. Hence, $C^*(\Lambda)\cong \mathcal{O}_X$. 
\begin{proof}
The proof is very similar to the analogous statement for Toeplitz algebras in Theorem~\ref{Toeplitz algebra iterated construction}. We begin by showing that $(\iota, \phi)$ is a Cuntz--Pimsner covariant Toeplitz representation of $X$ in $C^*(\Lambda)$. Exactly the same argument as in the proof of Theorem~\ref{Toeplitz algebra iterated construction} shows that $(\iota, \phi)$ is a Toeplitz representation.  It remains to check that $(\iota, \phi)$ is Cuntz--Pimsner covariant, i.e. $(\iota, \phi)^{(1)}(\psi(a))=\phi(a)$ for each $a\in J_X=\psi^{-1}(\mathcal{K}_{C^*(\Lambda^i)}(X))\cap \ker(\psi)^\perp$. 
Using the generating set of $J_X$ that we found in Proposition~\ref{putting everything together}, it suffices to show that:
\begin{enumerate}[label=\upshape(\roman*)]
\item
If $v\in \Lambda^0$ and $0<|v\Lambda^{e_i}|<\infty$, then
\[
(\iota,\phi)^{(1)}\big(\psi\big(as_v^{\Lambda^i}b\big)\big)=\phi\big(as_v^{\Lambda^i}b\big)
\]
for each $a,b\in C^*(\Lambda^i)$; and
\item
If $E\in \mathrm{FE}(\Lambda^i\setminus \Lambda^i H_{J_X})$ with $E\subseteq \bigcup_{j=1}^k \Lambda^{e_j}$, and there exists $F\subseteq r(E)\Lambda^{e_i}$ such that $E\cup F\in \mathrm{FE}(\Lambda)$, then
\[
(\iota,\phi)^{(1)}\big(\psi\big(a\Delta(s^{\Lambda^i})^Eb\big)\big)
=\phi\big(a\Delta(s^{\Lambda^i})^Eb\big)
\] 
for each $a,b\in C^*(\Lambda^i)$.
\end{enumerate}

Let us check (i) first. If $v\in \Lambda^0$ with $0<|v\Lambda^{e_i}|<\infty$ and $a,b\in C^*(\Lambda^i)$, making use of Proposition~\ref{vertex projection in Katsura ideal}, we see that
\begin{align*}
\psi\big(as_v^{\Lambda^i}b\big)
=\psi(a)\psi\big(s_v^{\Lambda^i}\big)\psi(b)
=\psi(a)\Big(\sum_{\mu\in v\Lambda^{e_i}}\Theta_{s_\mu^\Lambda,s_\mu^\Lambda}\Big)\psi(b)
&=\sum_{\mu\in v\Lambda^{e_i}}\Theta_{\psi(a)s_\mu^\Lambda,\psi(b^*)s_\mu^\Lambda}\\
&=\sum_{\mu\in v\Lambda^{e_i}}\Theta_{\phi(a)s_\mu^\Lambda,\phi(b^*)s_\mu^\Lambda}.
\end{align*}
Thus, using Proposition~\ref{vertex projection in Katsura ideal} again, we get that
\[
(\iota,\phi)^{(1)}\big(\psi\big(as_v^{\Lambda^i}b\big)\big)=\sum_{\mu\in v\Lambda^{e_i}}\phi(a)s_\mu^\Lambda{s_\mu^\Lambda}^*\phi(b)=\phi(a)s_v^\Lambda \phi(b)=\phi(as_v^{\Lambda^i}b).
\]
This completes the proof of (i). 

Next we check (ii). Fix $a,b\in C^*(\Lambda^i)$ and $E\in \mathrm{FE}(\Lambda^i\setminus \Lambda^i H_{J_X})$ with $E\subseteq \bigcup_{j=1}^k \Lambda^{e_j}$. Suppose that there exists $F\subseteq r(E)\Lambda^{e_i}$ such that $E\cup F\in \mathrm{FE}(\Lambda)$. By Lemma~\ref{extending to finite exhaustive set implies compact}, we have that
\begin{align*}
\psi\big(a\Delta(s^{\Lambda^i})^Eb\big)
=\psi(a)\psi\big(\Delta(s^{\Lambda^i})^E\big)\psi(b)
&=
\sum_{\substack{G\subseteq E\cup F\\G\cap F\neq \emptyset\\\mu\in \mathrm{MCE}(G)}}
(-1)^{(|G|+1)}\Theta_{\psi(a)s_\mu^\Lambda, \psi(b^*)s_\mu^\Lambda}\\
&=
\sum_{\substack{G\subseteq E\cup F\\G\cap F\neq \emptyset\\\mu\in \mathrm{MCE}(G)}}
(-1)^{(|G|+1)}\Theta_{\phi(a)s_\mu^\Lambda, \phi(b^*)s_\mu^\Lambda}.
\end{align*}
Thus, combining equations \eqref{splitting sum} and \eqref{applying lemma 5.12 to E} for the second equality, we see that
\begin{align*}
(\iota,\phi)^{(1)}\big(\psi\big(a\Delta(s^{\Lambda^i})^Eb\big)\big)
&=
\sum_{\substack{G\subseteq E\cup F\\G\cap F\neq \emptyset\\\mu\in \mathrm{MCE}(G)}}
(-1)^{(|G|+1)}\phi(a)s_\mu^\Lambda {s_\mu^\Lambda}^*\phi(b)\\
&=
\phi(a)
\phi\big(\Delta(s^{\Lambda^i})^E\big)
\phi(b)\\
&=
\phi\big(a\Delta(s^{\Lambda^i})^Eb\big).
\end{align*}

We conclude that $(\iota, \phi)$ is a Cuntz--Pimsner covariant Toeplitz representation of $X$ in $C^*(\Lambda)$. Hence, there exists a $*$-homomorphism $\iota \times_\mathcal{O} \phi:\mathcal{O}_X\rightarrow C^*(\Lambda)$ such that $(\iota \times_\mathcal{O} \phi)\circ j_X=\iota$ and $(\iota \times_\mathcal{O} \phi)\circ j_{C^*(\Lambda^i)}=\phi$, where $(j_X, j_{C^*(\Lambda^i)})$ is the universal Cuntz--Pimsner covariant Toeplitz representation of $X$.

Next, we show that the collection $\{u_\lambda:\lambda\in \Lambda\}\subseteq \mathcal{O}_X$ of partial isometries defined by 
$u_\lambda:=j_X^{\otimes d(\lambda)_i}\left(\Omega_{d(\lambda)_i}\left(s_\lambda^\Lambda\right)\right)$
is a Cuntz--Krieger $\Lambda$-family. The same calculations as in the proof of Theorem~\ref{Toeplitz algebra iterated construction} show that $\{u_\lambda:\lambda\in \Lambda\}$ is a Toeplitz--Cuntz--Krieger $\Lambda$-family. Thus, it remains to check that $\{u_\lambda:\lambda\in \Lambda\}$ satisfies (CK). By \cite[Theorem~C.1]{MR2069786}, it suffices to show that if $v\in \Lambda^0$ and $E\subseteq \bigcup_{j=1}^k v\Lambda^{e_j}$ belongs to $v\mathrm{FE}(\Lambda)$, then
$\Delta(u)^E=0$.

Firstly we consider the case where $E\cap \Lambda^{e_i}=\emptyset$. Then $E=E\cap \Lambda^i\in v\mathrm{FE}(\Lambda^i)$, and so
\[
\Delta(u)^E=j_{C^*(\Lambda^i)}\big(\Delta(s^{\Lambda^i})^E\big)=0,
\]
where the last equality comes from applying the Cuntz--Krieger relation in $C^*(\Lambda^i)$.

It remains to consider the situation where $E\cap \Lambda^{e_i}\neq\emptyset$. Using Lemma~\ref{products to sums general}, we get that
\begin{align*}
\Delta(u)^E
&=u_v+\sum_{\substack{\emptyset\neq G\subseteq E\\\lambda\in \mathrm{MCE}(G)}}(-1)^{|G|}u_\lambda u_\lambda^*\\
&=u_v+\sum_{\substack{\emptyset\neq G\subseteq E\cap \Lambda^i\\\lambda\in \mathrm{MCE}(G)}}(-1)^{|G|}u_\lambda u_\lambda^*
+\sum_{\substack{\emptyset\neq G\subseteq E\\ G\cap \Lambda^{e_i}\neq \emptyset \\ \lambda\in \mathrm{MCE}(G)}}(-1)^{|G|}u_\lambda u_\lambda^*\\
&=j_{C^*(\Lambda^i)}\big(\Delta(s^{\Lambda^i})^{E\cap \Lambda^i}\big)
+(j_X, j_{C^*(\Lambda^i)})^{(1)}\Bigg(\sum_{\substack{\emptyset\neq G\subseteq E\\ G\cap \Lambda^{e_i}\neq \emptyset\\\lambda\in \mathrm{MCE}(G)}}(-1)^{|G|}\Theta_{s_\lambda^\Lambda, s_\lambda^\Lambda}\Bigg).
\end{align*}
Since $\emptyset\neq E\cap \Lambda^{e_i}\subseteq v\Lambda^{e_i}$, Proposition~\ref{vertex projection in Katsura ideal} tells us that $s_v^{\Lambda^i}\in \mathrm{ker}(\psi)^\perp$. Thus, 
\[
\Delta(s^{\Lambda^i})^{E\cap \Lambda^i}=s_v^{\Lambda^i}\Delta(s^{\Lambda^i})^{E\cap \Lambda^i}\in \mathrm{ker}(\psi)^\perp
\]
as well. Since the Toeplitz representation $j_X$ is, by definition, Cuntz--Pimsner covariant, to establish that $\Delta(u)^E=0$, we need only verify that
\begin{equation}
\label{last thing to check for main theorem}
\psi\big(\Delta(s^{\Lambda^i})^{E\cap \Lambda^i}\big)
=
\sum_{\substack{\emptyset\neq G\subseteq E\\ G\cap \Lambda^{e_i}\neq \emptyset\\\lambda\in \mathrm{MCE}(G)}}(-1)^{(|G|+1)}\Theta_{s_\lambda^\Lambda, s_\lambda^\Lambda}.
\end{equation}
Again applying Lemma~\ref{products to sums general} (now to the Cuntz--Krieger $\Lambda$-family $\{s_\lambda^\Lambda:\lambda\in \Lambda\}$), and recalling that $E$ is finite and exhaustive in $\Lambda$, we get that
\begin{align*}
0=\Delta(s^\Lambda)^E
=\Delta(s^\Lambda)^{E\cap \Lambda^i}
+\sum_{\substack{\emptyset\neq G\subseteq E\\ G\cap \Lambda^{e_i}\neq \emptyset \\ \lambda\in \mathrm{MCE}(G)}}(-1)^{|G|}s_\lambda^\Lambda {s_\lambda^\Lambda}^*.
\end{align*}
Rearranging, we see that
\begin{align*}
\phi\big(\Delta(s^{\Lambda^i})^{E\cap \Lambda^i}\big)
=\Delta(s^\Lambda)^{E\cap \Lambda^i}
=\sum_{\substack{\emptyset\neq G\subseteq E\\ G\cap \Lambda^{e_i}\neq \emptyset \\ \lambda\in \mathrm{MCE}(G)}}(-1)^{(|G|+1)}s_\lambda^\Lambda {s_\lambda^\Lambda}^*,
\end{align*}
and so \eqref{last thing to check for main theorem} follows. 

Thus, $\{u_\lambda:\lambda\in \Lambda\}$ is a Cuntz--Krieger $\Lambda$-family in $\mathcal{O}_X$. The universal property of $C^*(\Lambda)$ then induces a $*$-homomorphism $\pi_u:C^*(\Lambda)\rightarrow \mathcal{O}_X$ such that $\pi_u\big(s_\lambda^\Lambda\big)=u_\lambda$ for each $\lambda\in \Lambda$. Exactly the same argument as in Theorem~\ref{Toeplitz algebra iterated construction} shows that $\pi_u$ and $\iota \times_{\mathcal{O}}\phi$ are mutually inverse. We conclude that $C^*(\Lambda)\cong \mathcal{O}_X$. 
\end{proof}
\end{thm}

\section{Relationships to other constructions}

\subsection{Graph correspondences}
\label{graph correspondences}

It is well known that if $E=(E^0, E^1,r,s)$ is a directed graph, then the graph algebra $C^*(E)$ can be realised as the Cuntz--Pimsner algebra of a Hilbert $C_0(E^0)$-bimodule. We briefly summarise this procedure and show that it is a special case of our construction when $k=1$. For $a,b\in C_0(E^0)$ and $x,y\in C_c(E^1)$ we define $a\cdot x\cdot b \in C_c(E^1)$ and $\langle x,y\rangle_{C_0(E^0)}\in C_0(E^0)$ by
\[
(a\cdot x\cdot b)(e):=a(r(e))x(e)b(s(e)) \quad \text{and} \quad \langle x,y\rangle_{C_0(E^0)}(v):=\sum_{e\in s^{-1}(v)}\overline{x(e)}y(e)
\] 
for each $e\in E^1$ and $v\in E^0$. Taking the completion of $C_c(E^1)$ with respect to the seminorm induced by $\langle \cdot, \cdot \rangle_{C_0(E^0)}$ gives a Hilbert $C_0(E^0)$-bimodule $X(E)$ (see \cite[Lemma~2.16]{MR1634408} for the details of this procedure), which we call the graph correspondence. For $v\in E^0$ and $e\in E^1$ we write $\delta_v$ and $\delta_e$ for the point masses of $v$ and $e$, which we view as elements of $C_0(E^0)$ and $C_c(E^1)\subseteq X(E)$ respectively. It follows from \cite[Example~8.8]{MR2135030} that  the Katsura ideal of the graph correspondence is $J_{X(E)}=\cspan\{\delta_v:0<|r^{-1}(v)|<\infty\}$. Moreover, \cite[Example~8.13]{MR2135030} tells us that the maps $s_v^E\mapsto j_{C_0(E^0)}(\delta_v)$ and $s_e^E\mapsto j_{X(E)}(\delta_e)$ for $v\in E^0$ and $e\in E^1$ induce an isomorphism from $C^*(E)$ to $\mathcal{O}_{X(E)}$.
 
Let $\Lambda$ be the path category of $E$. Then $\Lambda$ is a locally-convex $1$-graph, and we can apply our procedure from Section~\ref{section: realising as Cuntz--Pimsner algebra} to $\Lambda$. Removing all edges of degree $e_1$ from $\Lambda$ leaves $\Lambda^1=E^0$, and so $C^*(\Lambda^1)\cong C_0(E^0)$ via the isomorphism $s_v^{\Lambda^1}\mapsto \delta_v$. Similarly, $X=\cspan\{s_\lambda^\Lambda {s_\mu^\Lambda}^*:\lambda,\mu\in \Lambda, d(\lambda)_1=1, d(\mu)_1=0\}=\cspan\{s_e^\Lambda:e\in E^1\}$, which (if we identify the respective coefficient algebras) is isomorphic to $X(E)$ as a Hilbert bimodule via the map $s_e^\Lambda\mapsto \delta_e$. Thus, the isomorphism given by Theorem~\ref{main theorem} is the same as that given by \cite[Example~8.13]{MR2135030}. For this reason we like to think of the construction in Section~\ref{section: realising as Cuntz--Pimsner algebra} as a higher-rank graph correspondence. It is also not difficult to see that the description of the Katsura ideal for the graph correspondence given by \cite[Example~8.13]{MR2135030} is just a special case of Proposition~\ref{putting everything together}. Since $\Lambda^1\setminus \Lambda^1 H_{J_X}$ consists of just vertices, it follows that $\mathrm{FE}(\Lambda^1\setminus \Lambda^1 H_{J_X})=\emptyset$, and so $\tilde{B}_{J_X}=\emptyset$. Thus, Proposition~\ref{putting everything together} tells us that $J_X$ is generated as an ideal of $C^*(\Lambda^1)$ by the vertex projections $\{s_v^{\Lambda^1}:0<|vE^1|<\infty\}$, and so $J_X=\cspan\{s_v^{\Lambda^1}:0<|vE^1|<\infty\}$.

\subsection{Iterating the Nica--Toeplitz and Cuntz--Nica--Pimsner construction}
\label{Iterating the Nica--Toeplitz and Cuntz--Nica--Pimsner construction}

Sims and Yeend showed in \cite[Section~5.3]{MR2718947} that the Cuntz--Krieger algebra of a finitely aligned $k$-graph may be realised as the Cuntz--Nica--Pimsner algebra of a compactly aligned product system over $\N^k$. In \cite{2017arXiv170608626F} we showed how the Nica--Toeplitz and Cuntz--Nica--Pimsner algebras of a compactly aligned product system over $\N^k$ can be realised as iterated Toeplitz and iterated Cuntz--Pimsner algebras respectively. We now briefly explain how the results of the current paper can be deduced by combining these two constructions (at least for row finite graphs with no sources). For the relevant background information on product systems and their associated $C^*$-algebras we direct the reader to \cite{MR2718947}. 

Let $\Lambda$ be a finitely aligned $k$-graph. For $n\in \N^k$, $(\Lambda^n,\Lambda^0, r|_{\Lambda^n},s|_{\Lambda^n})$ is a directed graph, and we write $\mathbf{X}(\Lambda)_n$ for the associated graph correspondence. It can be shown that there exists an associative multiplication on $\mathbf{X}(\Lambda):=\bigsqcup_{n\in \N^k}\mathbf{X}(\Lambda)_n$ such that $\delta_\mu\delta_\nu=\delta_{s(\mu),r(\nu)}\delta_{\mu\nu}$ for each $\mu,\nu\in \Lambda$. This multiplication induces a Hilbert $C_0(\Lambda^0)$-bimodule isomorphism from the balanced tensor product $\mathbf{X}(\Lambda)_m\otimes_{C_0(\Lambda^0)} \mathbf{X}(\Lambda)_n$ to $\mathbf{X}(\Lambda)_{m+n}$ for each $m,n\in \N^k$, and so $\mathbf{X}(\Lambda)$ has the structure of a compactly aligned product system over $\N^k$ with coefficient algebra $C_0(\Lambda^0)$. It can then be shown that there is an isomorphism from the Nica--Toeplitz algebra $\mathcal{NT}_{\mathbf{X}(\Lambda)}$ to $\mathcal{T}C^*(\Lambda)$ that maps $i_{\mathbf{X}(\Lambda)}(\delta_\lambda)$ to $t_\lambda^\Lambda$ for each $\lambda\in \Lambda$. Similarly there exists an isomorphism from the Cuntz--Nica--Pimsner algebra $\mathcal{NO}_{\mathbf{X}(\Lambda)}$ to $C^*(\Lambda)$ mapping $j_{\mathbf{X}(\Lambda)}(\delta_\lambda)$ to $s_\lambda^\Lambda$ for each $\lambda\in \Lambda$. 

Fix $i\in\{1,\ldots, k\}$. Then $\mathbf{X}':=\bigsqcup_{\{n\in \N^k:n_i=0\}}\mathbf{X}(\Lambda)_n$ has the structure of a compactly aligned product system over $\N^{(k-1)}$. Clearly, $\mathbf{X}'$ is isomorphic as a product system to $\mathbf{X}(\Lambda^i)$, and so $\mathcal{NT}_{\mathbf{X}'}\cong \mathcal{T}C^*(\Lambda^i)$ and $\mathcal{NO}_{\mathbf{X}'}\cong C^*(\Lambda^i)$. By \cite[Proposition~4.2]{2017arXiv170608626F}, the inclusion $\mathbf{X}'\subseteq\mathbf{X}(\Lambda)$ induces an injective $*$-homomorphism $\phi_{\mathbf{X}'}^\mathcal{NT}:\mathcal{NT}_{\mathbf{X}'}\rightarrow \mathcal{NT}_{\mathbf{X}(\Lambda)}$. Similarly, \cite[Propositions~5.6]{2017arXiv170608626F} says that if $C_0(\Lambda^0)$ acts faithfully on each $\mathbf{X}(\Lambda)_n$ (equivalently, $\Lambda$ has no sources), then there exists an injective $*$- homomorphism $\phi_{\mathbf{X}'}^\mathcal{NO}:\mathcal{NO}_{\mathbf{X}'}\rightarrow \mathcal{NO}_{\mathbf{X}(\Lambda)}$ induced by the inclusion $\mathbf{X}'\subseteq \mathbf{X}(\Lambda)$. It follows from \cite[Propositions~4.3 and 4.6]{2017arXiv170608626F} that the closed subspace 
\[
\mathbf{Y}_1^\mathcal{NT}:=\cspan\{i_{\mathbf{X}(\Lambda)}(\mathbf{X}(\Lambda)_{e_i})\phi_{\mathbf{X}'}^\mathcal{NT}(\mathcal{NT}_{\mathbf{X}'})\}\subseteq \mathcal{NT}_{\mathbf{X}(\Lambda)}
\] 
has the structure of a Hilbert $\mathcal{NT}_{\mathbf{X}'}$-bimodule with operations 
\[
a\cdot y\cdot b=\phi_{\mathbf{X}'}^\mathcal{NT}(a)y\phi_{\mathbf{X}'}^\mathcal{NT}(b)
\quad \text{and} \quad
\langle y,w\rangle_{\mathcal{NT}_{\mathbf{X}'}}=(\phi_{\mathbf{X}'}^\mathcal{NT})^{-1}(y^*w) 
\]
for $y,w\in \mathbf{Y}_1^\mathcal{NT}$ and $a,b\in \mathcal{NT}_{\mathbf{X}'}$. After identifying the coefficient algebras $\mathcal{NT}_{\mathbf{X}'}$ and $\mathcal{T}C^*(\Lambda^i)$, routine calculations show that the map 
\[
i_{\mathbf{X}(\Lambda)}(\delta_\lambda)\phi_{\mathbf{X}'}^\mathcal{NT}\big(i_{\mathbf{X}'}(\delta_\mu)i_{\mathbf{X}'}(\delta_\nu)^*\big)
\mapsto t_{\lambda\mu}^\Lambda {t_\nu^\Lambda}^* 
\]
for $\lambda\in \Lambda^{e_i}$, $\mu,\nu\in \Lambda^i$ with $s(\lambda)=r(\mu)$, $s(\mu)=s(\nu)$ extends to a Hilbert bimodule isomorphism from $\mathbf{Y}_1^\mathcal{NT}$ to the bimodule $X$ constructed in Proposition~\ref{module for TCK algebra}. Using $\phi_{\mathbf{X}'}^\mathcal{NO}$ in place of $\phi_{\mathbf{X}'}^\mathcal{NT}$, we also have a Hilbert $\mathcal{NO}_{\mathbf{X}'}$-bimodule
\[
\mathbf{Y}_1^\mathcal{NO}:=\cspan\{j_{\mathbf{X}(\Lambda)}(\mathbf{X}(\Lambda)_{e_i})\phi_{\mathbf{X}'}^\mathcal{NO}(\mathcal{NO}_{\mathbf{X}'})\}\subseteq \mathcal{NO}_{\mathbf{X}(\Lambda)},
\] 
which we can identify with the bimodule $X$ from Section~\ref{section: realising as Cuntz--Pimsner algebra}. Finally, \cite[Theorems~4.17 and 5.20]{2017arXiv170608626F} imply that the inclusions $\mathbf{Y}_1^\mathcal{NT}\subseteq \mathcal{NT}_{\mathbf{X}(\Lambda)}$ and $\mathbf{Y}_1^\mathcal{NO}\subseteq \mathcal{NO}_{\mathbf{X}(\Lambda)}$ induce isomorphisms $\mathcal{T}_{\mathbf{Y}_1^\mathcal{NT}}\cong \mathcal{NT}_{\mathbf{X}(\Lambda)}$ and $\mathcal{O}_{\mathbf{Y}_1^\mathcal{NO}}\cong \mathcal{NO}_{\mathbf{X}(\Lambda)}$ (the latter assuming that $C_0(\Lambda^0)$ acts faithfully and compactly  on each $\mathbf{X}(\Lambda)_n$).
Consequently, in the situation where $\Lambda$ is row finite and has no sources,  the main result of this paper (Theorem~\ref{main theorem}) can be obtained by combining \cite[Proposition~5.4]{MR2718947} and the results of \cite{2017arXiv170608626F}.

\subsection{Semi-saturated circle actions and generalised crossed products}

Let $B$ be a $C^*$-algebra and $\alpha:\T\rightarrow \mathrm{Aut}(B)$ an action of the circle group. For each $n\in \Z$, we define the $n$th spectral subspace for $\alpha$ to be $B_n:=\{b\in B:\alpha_z(b)=z^nb \text{ for each } z\in \T\}$. It is routine to check that each $B_n$ is a closed subspace of $B$, and $B_0$ (which we call the fixed point algebra of $\alpha$) is also closed under multiplication and taking adjoints. Moreover, $B_n^*=B_{-n}$ and $B_n B_m\subseteq B_{n+m}$ for each $n,m\in \Z$. In particular, since $B_0 B_1 B_0\subseteq B_1$ and $B_1^* B_1\subseteq B_0$, \cite[Lemma~3.2(1)]{MR1426840} tells us that $B_1$ is a Hilbert $B_0$-bimodule with left and right actions given by multiplication and inner product $\langle \xi,\eta \rangle_{B_0}=\xi^*\eta$ for each $\xi,\eta\in B_1$. In fact, since $B_1 B_1^*\subseteq B_0$, $B_1$ also has a left $B_0$-valued inner product given by ${}_{B_0}\langle \xi,\eta \rangle=\xi\eta^*$ for each $\xi,\eta\in B_1$, which gives $B_1$ the structure of a left Hilbert $B_0$-bimodule. It is straightforward to see that the two inner products satisfy the imprimitivity condition 
\begin{equation}
\label{imprim. condition}
\xi\cdot \langle \eta,\mu\rangle_{B_0}={}_{B_0}\langle \xi,\eta \rangle\cdot \mu \quad \text{for $\xi,\eta,\mu\in B_1$.}
\end{equation}
Thus, if $B_1^*B_1=B_0=B_1B_1^*$, then both of these inner products are full, and $B_1$ is a $B_0$--$B_0$ imprimitivity bimodule (see \cite[Definition~3.1]{MR1634408}). If the action $\alpha$ is semi-saturated  in the sense that $B$ is generated as a $C^*$-algebra by the fixed point algebra $B_0$ and the first spectral subspace $B_1$ (see \cite[Definition~4.1]{MR1276163}), then \cite[Theorem~3.1]{MR1467459} says that $B$ can be realised as the generalised crossed product $B_0 \rtimes_{B_1}\Z$. Proposition~3.7 of \cite{MR2029622} tells us that $B_0 \rtimes_{B_1}\Z$ is canonically isomorphic to $\mathcal{O}_{B_1}$, and so we conclude that any $C^*$-algebra with a semi-saturated circle action may be realised as the Cuntz--Pimsner algebra of a Hilbert bimodule whose coefficient algebra is equal to the fixed point algebra of the action. 

Now suppose that $\Lambda$ is some locally-convex finitely aligned $k$-graph and resume the notation of Section~\ref{section: realising as Cuntz--Pimsner algebra}. Let $\gamma_i^\Lambda:\T\rightarrow \mathrm{Aut}(C^*(\Lambda))$ denote the restriction of the gauge action $\gamma^\Lambda:\T^k\rightarrow \mathrm{Aut}(C^*(\Lambda))$ to the $i$th coordinate of $\T^k$. The $n$th spectral subspace for $\gamma_i^\Lambda$ is then $C^*(\Lambda)_n=\cspan\{s_\lambda^\Lambda {s_\mu^\Lambda}^*:\lambda,\mu\in \Lambda, d(\lambda)_i-d(\mu)_i=n\}$. Thus, $\phi(C^*(\Lambda^i))\subseteq C^*(\Lambda)_0$ and $X\subseteq C^*(\Lambda)_1$. Theorem~\ref{main theorem} tells us that $\phi(C^*(\Lambda^i))$ and $X$ generate $C^*(\Lambda)$, and so we see immediately that $\gamma_i^\Lambda$ is a semi-saturated action. Consequently, the discussion in the previous paragraph shows that $C^*(\Lambda)$ can be realised as the Cuntz--Pimsner algebra of the Hilbert $C^*(\Lambda)_0$-bimodule $C^*(\Lambda)_1$. We now explain how this decomposition of $C^*(\Lambda)$ as a Cuntz--Pimsner algebra relates to that given by Theorem~\ref{main theorem}.

The first point to note is that whilst $\phi(C^*(\Lambda^i))$ and $X$ are always subsets of  $C^*(\Lambda)_0$ and $C^*(\Lambda)_1$ respectively, these containments are usually strict. Thus, the descriptions of $C^*(\Lambda)$ given by \cite[Theorem~3.1]{MR1467459} and Theorem~\ref{main theorem} are not the same. For example consider the $1$-graph $\Sigma$ consisting of a single vertex and $n\geq 2$ loops $\{e_1,\ldots, e_n\}$. Then $C^*(\Sigma)$ is the Cuntz algebra $\mathcal{O}_n$, and removing all the edges from $\Sigma$ leaves the graph $\Sigma^1$ consisting of just one vertex. Hence, $\phi(C^*(\Sigma^1))\cong \C$, whilst we see that $s_{e_i}^\Sigma {s_{e_i}^\Sigma}^*\in C^*(\Sigma)_0\setminus \phi(C^*(\Sigma^1))$ for each $i\in \{1,\ldots,n\}$ (in fact $C^*(\Sigma)_0$ is the UHF algebra $M_{n^\infty}$). In general, the bimodules $X$ and $C^*(\Lambda)_1$ (and their respective coefficient algebras $C^*(\Lambda^i)$ and $C^*(\Lambda)_0$) are related by Pimsner's process of extending the scalars (see \cite[Section~2]{MR1426840} and \cite[Section~3.1]{2016arXiv160508593A} for the details): the map 
\[
s_\lambda^\Lambda {s_\mu^\Lambda}^*\mapsto s_{\lambda(0,e_i)}^\Lambda \otimes_{C^*(\Lambda^i)}s_{\lambda(e_i,d(\lambda))}^\Lambda {s_\mu^\Lambda}^* \quad \text{for $\lambda,\mu\in \Lambda$ with $d(\lambda)_i-d(\mu)_i=1$}
\] 
extends by linearity and continuity to an Hilbert $C^*(\Lambda)_0$-bimodule isomorphism from $C^*(\Lambda)_1$ to $X\otimes_{C^*(\Lambda^i)}C^*(\Lambda)_0$. 

Another key difference between our procedure for realising $C^*(\Lambda)$ as a Cuntz--Pimsner algebra and that of \cite{MR1467459} is the existence of a left inner product satisfying the imprimitivity condition. Since $C^*(\Lambda)_1{C^*(\Lambda)_1}^*\subseteq C^*(\Lambda)_0$, the spectral subspace $C^*(\Lambda)_1$ carries a left $C^*(\Lambda)_0$-valued inner product given by ${}_{C^*(\Lambda)_0}\langle \xi,\eta \rangle=\xi\eta^*$, and the left and right inner products on $C^*(\Lambda)_1$ satisfy \eqref{imprim. condition}. On the other hand, it is not true in general that $XX^*\subseteq \phi(C^*(\Lambda^i))$. For example, if we return to the $1$-graph $\Sigma$ discussed above, then $s_{e_i}^\Sigma\in X$ for each $i\in\{1,\ldots,n\}$, but $s_{e_i}^\Sigma{s_{e_i}^\Sigma}^*\not\in \phi(C^*(\Sigma^1))$. It would be interesting to see what the condition $XX^*\subseteq \phi(C^*(\Lambda^i))$ implies about the structure of the graph $\Lambda$: as the next result shows, this condition determines precisely when $X$ also has the structure of a left Hilbert $C^*(\Lambda^i)$-bimodule and the two inner products satisfy the imprimitivity condition. 

\begin{prop}
\label{existence of left IP satisfying imprimitivity condition}
Let $\Lambda$ be a locally-convex $k$-graph and let $X$ be the Hilbert $C^*(\Lambda^i)$-bimodule constructed in Section~\ref{section: realising as Cuntz--Pimsner algebra}. Then there exists a left $C^*(\Lambda^i)$-valued inner product ${}_{C^*(\Lambda^i)}\langle \cdot,\cdot \rangle$ giving $X$ the structure of a left Hilbert $C^*(\Lambda^i)$-bimodule and satisfying the imprimitivity condition 
\[
{}_{C^*(\Lambda^i)}\langle x,y \rangle \cdot z=x\cdot \langle y,z \rangle_{C^*(\Lambda^i)} \quad \text{for each $x,y,z\in X$}
\]
if and only if $XX^*\subseteq \phi(C^*(\Lambda^i))$.
\begin{proof}
It follows from \cite[Proposition~5.18]{MR2102572} that $X$ has left $C^*(\Lambda^i)$-valued inner product with the required properties if and only if $\mathcal{K}_{C^*(\Lambda^i)}(X)\subseteq \psi(J_X)$. By \cite[Lemma~3.2(3)]{MR1426840}, the $*$-homomorphism $(\iota,\phi)^{(1)}:\mathcal{K}_{C^*(\Lambda^i)}(X)\rightarrow XX^*$, which sends $\Theta_{x,y}$ to $\iota(x)\iota(y)^*=xy^*$ for each $x,y\in X$, is an isomorphism. By Theorem~\ref{main theorem}, the Toeplitz representation $(\iota,\phi)$ is Cuntz--Pimsner covariant, and so $(\iota,\phi)^{(1)}\circ \psi$ and $\phi$ agree on $J_X$. Thus, $X$ has a left $C^*(\Lambda^i)$-valued inner product with the required properties if and only if $XX^*\subseteq \phi(J_X)$. Consequently, to prove the result it remains to show that $XX^*\subseteq \phi(C^*(\Lambda^i))$ implies $XX^*\subseteq \phi(J_X)$.

Suppose $XX^*\subseteq \phi(C^*(\Lambda^i))$. Since $XX^*=\cspan\{s_\lambda^\Lambda {s_\mu^\Lambda}^*:\lambda,\mu\in \Lambda, d(\lambda)_i=d(\mu)_i=1\}$, by linearity and continuity it suffices to show that $s_\lambda^\Lambda {s_\mu^\Lambda}^*\in \phi(J_X)$ for each $\lambda,\mu\in \Lambda$ with $d(\lambda)_i=d(\mu)_i=1$ and $s(\lambda)=s(\mu)$. By assumption, $s_\lambda^\Lambda {s_\mu^\Lambda}^*=\phi(a)$ for some $a\in C^*(\Lambda^i)$, and so we need only check that $a\in J_X$. Since $s_\lambda^\Lambda, s_\mu^\Lambda\in X$, we see immediately that $\psi(a)=\Theta_{s_\lambda^\Lambda, s_\mu^\Lambda}$, and so $a\in \psi^{-1}(\mathcal{K}_{C^*(\Lambda^i)}(X))$. Thus, it remains to show that $a\in \mathrm{ker}(\psi)^\perp$. By Proposition~\ref{generators of ker(psi)}, $\mathrm{ker}(\psi)$ is generated as an ideal of $C^*(\Lambda^i)$ by the collection of vertex projections $\{s_v^{\Lambda^i}\in \Lambda^0:v\Lambda^{e_i}=\emptyset\}$. Hence, by linearity and continuity, if $a s_\nu^{\Lambda^i} {s_\eta^{\Lambda^i}}^*s_v^{\Lambda^i}=0$ whenever $\nu,\eta\in \Lambda^i$, $v\in \Lambda^0$ with $v\Lambda^{e_i}=\emptyset$, and $s(\nu)=s(\eta)$, $r(\eta)=v$, then $a\in \mathrm{ker}(\psi)^\perp$ as required. Since $\phi(a s_\nu^{\Lambda^i} {s_\eta^{\Lambda^i}}^*s_v^{\Lambda^i})=s_\lambda^\Lambda {s_\mu^\Lambda}^*s_\nu^\Lambda {s_\eta^\Lambda}^*s_v^\Lambda$ and $\phi$ is injective, we need only show that $r(\mu)\neq r(\nu)$. Looking for a contradiction, suppose that $r(\mu)=r(\nu)$. Since $d(\mu)_i=1$ and $d(\nu)_i=0$, Lemma~\ref{extending local convexity} tells us that $s(\nu)\Lambda^{e_i}\neq \emptyset$. As $s(\nu)=s(\eta)$, the factorisation property implies that $v\Lambda^{e_i}\neq \emptyset$, which is impossible. Hence, $a\in \mathrm{ker}(\psi)^\perp$, and we conclude that $a\in J_X$. Thus, $s_\lambda^\Lambda {s_\mu^\Lambda}^*=\phi(a)\in \phi(J_X)$ as required. 
\end{proof}
\end{prop}

Finally, we note that even when $XX^*$ is not contained in $\phi(C^*(\Lambda^i))$, $X$ can still have the structure of a left $C^*(\Lambda^i)$-bimodule, provided we do not require that the two inner products satisfy the imprimitivity condition. If $E$ is a directed graph, then, as shown in \cite[Proposition~3.8]{2018arXiv180408114R}, 
the graph correspondence $X(E)$ has a left $C_0(E^0)$-valued inner product given by
\[
{}_{C_0(E^0)}\langle f,g\rangle(v)=\sum_{e\in r^{-1}(v)}f(e)\overline{g(e)} \quad \text{for $f,g\in C_c(E^1)$, $v\in \Lambda^0$,}
\]
which gives $X(E)$ the structure of a left Hilbert $C_0(E^0)$-bimodule. It is straightforward to see that the left and right $C_0(E_0)$-valued inner products on $X(E)$ do not satisfy the imprimitivity condition (as predicted by Proposition~\ref{existence of left IP satisfying imprimitivity condition}): if $e,f\in E^1$ are distinct edges with common range, then
\[
{}_{C_0(E^0)}\langle \delta_e,\delta_e\rangle\cdot \delta_f=\delta_{r(e)}\cdot\delta_f=\delta_f\neq 0=\delta_{e,f}\delta_e\cdot\delta_{s(e)}=\delta_e\cdot \langle \delta_e,\delta_f\rangle_{C_0(E^0)}.
\]
However, the left and right inner products on $X(E)$ are compatible in the sense that the right action of $C_0(E^0)$ is adjointable with respect to the left inner product and the left action of $C_0(E^0)$ is adjointable with respect to the right inner product. As shown in \cite[Remark~1.9]{MR1304338}, this compatibility condition is automatic if the two inner products satisfy the imprimitivity condition. Thus, the graph correspondence $X(E)$ is a bi-Hilbertian $C^*(\Lambda^i)$-bimodule in the sense of \cite[Definition~2.1]{MR3622236}. Unfortunately, we have so far been unable to determine whether an analogous left inner product exists for the bimodule $X$ associated to graphs of rank $2$ or more. This is certainly an issue worth exploring further: if $X$ has the structure of a bi-Hilbertian $C^*(\Lambda^i)$-bimodule, then the results of \cite{2016arXiv160508593A, MR3663586, MR3622236, 2018arXiv180408114R} could be applied to higher-rank graph algebras. 

\section*{Acknowledgements}
The results in this article extend the work from the second chapter of my PhD thesis. I would like to thank my supervisors Adam Rennie and Aidan Sims at the University of Wollongong for their advice and encouragement during my PhD and during the writing of this article. 

\appendix
\section{Gauge-invariant ideals of higher-rank graph algebras}

Suppose $\Sigma$ is a finitely aligned $k$-graph and $I$ is a gauge-invariant ideal of $C^*(\Sigma)$. If 
\[
H_I:=\{v\in \Sigma^0:s_v^\Sigma\in I\}
\quad \text{and}\quad
B_{I}:=\big\{E\in \mathrm{FE}(\Sigma\setminus \Sigma H_I):\Delta(s^\Sigma)^E\in I\big\},
\]
then \cite[Theorem~4.6]{MR3262073} tells us that $I$ is generated as an ideal of $C^*(\Sigma)$ by the collection of projections
\[
\big\{s_v^{\Sigma}:v\in H_I\}\cup \{\Delta(s^\Sigma)^E:E\in B_I\big\}.
\]
In this appendix we show that to get a generating set for $I$, we need only consider those finite exhaustive sets in the collection $B_I$ consisting of edges. We make use of this refinement of \cite[Theorem~4.6]{MR3262073} in our proof of Lemma~\ref{extending to finite exhaustive set implies compact}.

To this end, we prove in the next result that if 
\[
\tilde{B}_{I}:=\Big\{E\in B_I: E\subseteq \bigcup_{j=1}^k \Sigma^{e_j}\Big\}
\] 
and $F\in B_I$, then $\Delta(s^\Sigma)^F$ belongs to the ideal of $C^*(\Sigma)$ generated by the collection of projections $\{\Delta(s^\Sigma)^E:E\in \tilde{B}_I\}$. Our proof uses the same techniques as deployed in \cite[Appendix~C]{MR2069786} to show that a Toeplitz--Cuntz--Krieger $\Sigma$-family $\{q_\lambda:\lambda\in \Sigma\}$ satisfies relation (CK) if and only if $\Delta(q)^E=0$ for each $E\in \mathrm{FE}(\Sigma)$ with $E\subseteq \bigcup_{j=1}^k \Sigma^{e_j}$. 

\begin{prop}
\label{sufficient to consider exhaustive sets of edges}
Let $\Sigma$ be a finitely aligned $k$-graph. Suppose $I$ is a gauge-invariant ideal of $C^*(\Sigma)$. Let $J$ denote the ideal of $C^*(\Sigma)$ generated by the projections $\{\Delta(s^\Sigma)^E:E\in \tilde{B}_I\}$. Then 
\begin{equation}
\label{what we need to show}
E\in B_I \implies \Delta(s^\Sigma)^E \in J. 
\end{equation}
\begin{proof}
We will prove the result using induction on $L(E):=\sum_{j=1}^k \max\{d(\lambda)_j:\lambda\in E\}$. If $E\in B_I$ and $L(E)=1$, then $E\subseteq \Sigma^{e_j}$ for some $j\in \{1,\ldots, k\}$, and so $E\in \tilde{B}_I$. Thus, $\Delta(s^\Sigma)^E \in J$ as required.

Now let $n\geq 1$ and suppose that \eqref{what we need to show} holds whenever $L(E)\leq n$. Fix $F\in B_I$ with $L(F)=n+1$. Define
\[
I(F):=\bigcup_{j=1}^k \{\lambda(0,e_j):\lambda\in F, \ d(\lambda)_j\geq 1\}.
\]
Since $F\in \mathrm{FE}(\Sigma\setminus \Sigma H_I)$, \cite[Lemma~C.6]{MR2069786} tells us that $I(F)\in \mathrm{FE}(\Sigma\setminus \Sigma H_I)$. Moreover, since $\Delta(s^\Sigma)^F\in I$, and each element of $F$ extends an element of $I(F)$, we see that
\[
\Delta(s^\Sigma)^{I(F)}=\Delta(s^\Sigma)^{I(F)}\Delta(s^\Sigma)^{F}\in I.
\]
Thus, $I(F)\in \tilde{B}_{I}$.

For each $\mu\in I(F)$, we also define
\[
\mathrm{Ext}_{\Sigma\setminus \Sigma H_I}(\mu;F):=\bigcup_{\lambda\in F}\{\alpha\in s(\mu)(\Sigma\setminus \Sigma H_I):\mu\alpha\in\mathrm{MCE}(\mu,\lambda)\}.
\]
We claim that
\begin{equation}
\label{using the induction hypothesis}
\Delta(s^\Sigma)^{\mathrm{Ext}_{\Sigma\setminus \Sigma H_I}(\mu;F)}\in J.
\end{equation}
If $s(\mu)\in \mathrm{Ext}_{\Sigma\setminus \Sigma H_I}(\mu;F)$, then $\Delta(s^\Sigma)^{\mathrm{Ext}_{\Sigma\setminus \Sigma H_I}(\mu;F)}=0$ which is certainly in $J$, so we may as well suppose that $s(\mu)\not\in \mathrm{Ext}_{\Sigma\setminus \Sigma H_I}(\mu;F)$. Since $F\in \mathrm{FE}(\Sigma\setminus \Sigma H_I)$, it follows from \cite[Lemma~C.5]{MR2069786} that $\mathrm{Ext}_{\Sigma\setminus \Sigma H_I}(\mu;F)\in s(\mu)\mathrm{FE}(\Sigma\setminus \Sigma H_I)$. Since $\Delta(s^\Sigma)^F\in I$, we can use \cite[Lemma~3.7]{MR3262073} to see that
\[
\Delta(s^\Sigma)^{\mathrm{Ext}_{\Sigma\setminus \Sigma H_I}(\mu;F)}={s_\mu^\Sigma}^*\Delta(s^\Sigma)^F s_\mu^{\Sigma}\in I.
\]
Thus, $\mathrm{Ext}_{\Sigma\setminus \Sigma H_I}(\mu;F)\in B_I$. By \cite[Lemma~C.8]{MR2069786}, $L(\mathrm{Ext}_{\Sigma\setminus \Sigma H_I}(\mu;F))\leq n$, and so we may apply the inductive hypothesis to conclude that $\Delta(s^\Sigma)^{\mathrm{Ext}_{\Sigma\setminus \Sigma H_I}(\mu;F)}\in J$. This completes the proof of claim \eqref{using the induction hypothesis}.

Observe that if $\mu\in I(F)$ and $\lambda\in \mathrm{Ext}_{\Sigma\setminus \Sigma H_I}(\mu;F)$, then there exists $\nu\in F$ such that $\mu\lambda\in \mathrm{MCE}(\mu,\nu)$, and so $s_{r(F)}^\Sigma-s_{\nu}^\Sigma{s_{\nu}^\Sigma}^*=(s_{r(F)}^\Sigma-s_{\nu}^\Sigma{s_{\nu}^\Sigma}^*)(s_{r(F)}^\Sigma-s_{\mu\lambda}^\Sigma{s_{\mu\lambda}^\Sigma}^*)$. Hence,
\[
\Delta(s^\Sigma)^F=\Delta(s^\Sigma)^F\prod_{\mu\in I(F)}\prod_{\lambda\in \mathrm{Ext}_{\Sigma\setminus \Sigma H_I}(\mu;F)}(s_{r(F)}^\Sigma-s_{\mu\lambda}^\Sigma{s_{\mu\lambda}^\Sigma}^*).
\]
Thus, $\Delta(s^\Sigma)^F $ will belong to the ideal $J$ provided
\begin{equation}
\label{last thing to check}
\prod_{\mu\in I(F)}\prod_{\lambda\in \mathrm{Ext}_{\Sigma\setminus \Sigma H_I}(\mu;F)}(s_{r(F)}^\Sigma-s_{\mu\lambda}^\Sigma{s_{\mu\lambda}^\Sigma}^*)\in J.
\end{equation}
For $\mu\in I(F)$, we have
\[
s_{r(F)}^\Sigma-s_\mu^\Sigma{s_\mu^\Sigma}^*
=(s_{r(F)}^\Sigma-s_\mu^\Sigma{s_\mu^\Sigma}^*)\Bigg(\prod_{\lambda\in \mathrm{Ext}_{\Sigma\setminus \Sigma H_I}(\mu;F)}(s_{r(F)}^\Sigma-s_{\mu\lambda}^\Sigma{s_{\mu\lambda}^\Sigma}^*)\Bigg).
\]
On the other hand,
\begin{align*}
&(s_{r(F)}^\Sigma-s_\mu^\Sigma{s_\mu^\Sigma}^*)\Bigg(\prod_{\lambda\in \mathrm{Ext}_{\Sigma\setminus \Sigma H_I}(\mu;F)}(s_{r(F)}^\Sigma-s_{\mu\lambda}^\Sigma{s_{\mu\lambda}^\Sigma}^*)\Bigg)\\
&\quad
=s_{r(F)}^\Sigma\Bigg(\prod_{\lambda\in \mathrm{Ext}_{\Sigma\setminus \Sigma H_I}(\mu;F)}(s_{r(F)}^\Sigma-s_{\mu\lambda}^\Sigma{s_{\mu\lambda}^\Sigma}^*)\Bigg)
-s_\mu^\Sigma{s_\mu^\Sigma}^*\Bigg(\prod_{\lambda\in \mathrm{Ext}_{\Sigma\setminus \Sigma H_I}(\mu;F)}(s_{r(F)}^\Sigma-s_{\mu\lambda}^\Sigma{s_{\mu\lambda}^\Sigma}^*)\Bigg)\\
&\quad
=\Bigg(\prod_{\lambda\in \mathrm{Ext}_{\Sigma\setminus \Sigma H_I}(\mu;F)}(s_{r(F)}^\Sigma-s_{\mu\lambda}^\Sigma{s_{\mu\lambda}^\Sigma}^*)\Bigg)
-\Bigg(\prod_{\lambda\in \mathrm{Ext}_{\Sigma\setminus \Sigma H_I}(\mu;F)}(s_\mu^\Sigma{s_\mu^\Sigma}^*-s_{\mu\lambda}^\Sigma{s_{\mu\lambda}^\Sigma}^*)\Bigg)\\
&\quad
=\Bigg(\prod_{\lambda\in \mathrm{Ext}_{\Sigma\setminus \Sigma H_I}(\mu;F)}(s_{r(F)}^\Sigma-s_{\mu\lambda}^\Sigma{s_{\mu\lambda}^\Sigma}^*)\Bigg)
-s_\mu^\Sigma\Delta(s^\Sigma)^{\mathrm{Ext}_{\Sigma\setminus \Sigma H_I}(\mu;F)}{s_\mu^\Sigma}^*.
\end{align*}
Combining the last two calculations, we see that
\[
\prod_{\lambda\in \mathrm{Ext}_{\Sigma\setminus \Sigma H_I}(\mu;F)}(s_{r(F)}^\Sigma-s_{\mu\lambda}^\Sigma{s_{\mu\lambda}^\Sigma}^*)
=(s_{r(F)}^\Sigma-s_\mu^\Sigma{s_\mu^\Sigma}^*)+s_\mu^\Sigma\Delta(s^\Sigma)^{\mathrm{Ext}_{\Sigma\setminus \Sigma H_I}(\mu;F)}{s_\mu^\Sigma}^*.
\]
Taking the product as $\mu$ ranges over $I(F)$ and recalling that $\Delta(s^\Sigma)^{\mathrm{Ext}_{\Sigma\setminus \Sigma H_I}(\mu;F)}\in J$ and $\prod_{\mu\in I(F)}(s_{r(F)}^\Sigma-s_\mu^\Sigma{s_\mu^\Sigma}^*)=\Delta(s^\Sigma)^{I(F)}\in J$, we conclude that \eqref{last thing to check} holds. Thus $\Delta(s^\Sigma)^F \in J$, and the result holds by induction. 
\end{proof}
\end{prop}

\begin{rem}
The results of \cite{MR3262073} in fact deal with the more general situation of twisted relative Cuntz--Krieger algebras associated to higher-rank graphs. In this more general setting the analogous version of Proposition~\ref{sufficient to consider exhaustive sets of edges} still holds using exactly the same argument: all our calculations take place in the diagonal subalgebra (which is the same regardless of the twist) and do not make use of relation (CK). 
\end{rem}

\end{document}